\documentclass[11pt,a4paper]{article}
\usepackage[utf8]{inputenc}
\usepackage[margin=1in]{geometry}
\usepackage{amsmath,amssymb,amsthm}
\usepackage{xcolor}
\usepackage{graphicx}
\usepackage{longtable,booktabs,array}
\usepackage{caption}
\usepackage{calc}
\usepackage{etoolbox}
\usepackage{bm}
\usepackage[T1]{fontenc}
\usepackage{lmodern}
\usepackage[colorlinks=true,linkcolor=blue,citecolor=blue,urlcolor=blue]{hyperref}
\IfFileExists{footnotehyper.sty}{\usepackage{footnotehyper}}{\usepackage{footnote}}
\makesavenoteenv{longtable}
\IfFileExists{xurl.sty}{\usepackage{xurl}}{}
\allowdisplaybreaks
\graphicspath{{figures/}}

\providecommand{\tightlist}{%
  \setlength{\itemsep}{0pt}\setlength{\parskip}{0pt}}
\makeatletter
\patchcmd\longtable{\par}{\if@noskipsec\mbox{}\fi\par}{}{}
\def\fps@figure{htbp}
\makeatother

\title{Edmunds--Evans essential spectra of the generalized
Constantin--Lax--Majda linearization}
\author{Jie Xu\\[2pt]
\small Department of Mechanical and Industrial Engineering,\\
\small University of Illinois Chicago, Chicago, IL, USA\\
\small \texttt{jiexu@uic.edu}\\
\small ORCID: \href{https://orcid.org/0000-0001-5765-7431}{0000-0001-5765-7431}}
\date{}

\hypersetup{
  pdftitle={Edmunds--Evans essential spectra of the generalized Constantin--Lax--Majda linearization},
  pdfauthor={Jie Xu},
  pdfkeywords={essential spectrum; Fredholm index; non-self-adjoint operator; Hilbert transform; self-similar blow-up; generalized Constantin-Lax-Majda equation},
  pdfsubject={Essential spectra of the positive-advection gCLM self-similar linearization}
}

\begin{document}
\maketitle

\begin{abstract}

Fix \(0<a<1\) and a smooth odd self-similar collapse profile
\((\Omega,c_l)\) of the generalized Constantin--Lax--Majda equation,
where \(c_l\) is the focusing exponent. Let \(L_a\) be the
corresponding linearization on the origin-\(H^2\) realization
\(X=\{\varphi\in L^2(0,\infty):\varphi(0)=0,\ \varphi''\in L^2(0,\infty)\}\).
Under explicit fixed-profile transport, regularity, weighted-derivative,
and far-field hypotheses, the first three Edmunds--Evans essential
spectra are exactly the union of two (not necessarily distinct) vertical
indicial lines:
\(\operatorname{Re}\lambda=-1+c_l/2\) at infinity and
\(\operatorname{Re}\lambda=-\tilde c/2\) at the origin, where
\(\tilde c=c_l+a(H\Omega)(0)\) and \(H\) is the Hilbert transform.
Away from these lines, \(L_a-\lambda\) is Fredholm. Its index is zero on
the two exterior components and is \(+1\) or \(-1\) in the open strip
between the lines, according to their orientation. Consequently, the
fourth and fifth Edmunds--Evans spectra and the Browder essential
spectrum are exactly the closed inter-line strip. In the orientation
yielding Fredholm index \(+1\), every point of the open strip is an
eigenvalue.

We further prove that every fixed point of the Huang--Qin--Wang--Wei
positive-advection construction with \(0<a<400/(848-9\pi^2)\), after normalization,
satisfies the complete fixed-profile hypothesis package. For all
sufficiently small positive \(a\), every such fixed point has distinct
endpoint lines in the orientation yielding Fredholm index \(+1\). Thus actual
positive-advection collapse profiles realize the closed Browder band. An
additional kernel-nondegeneracy condition (ND) converts the Fredholm index
into the exact individual kernel and cokernel dimensions. We prove this
condition outside an explicit vertical energy slab and, at sufficiently large
imaginary part, on every closed vertical substrip of the open band separated
from the far-field line.
Any failures in the open band are locally finite and can accumulate at
high frequency only toward the far-field line.

\end{abstract}



\section{1. Introduction}\label{introduction}

\subsection{1.1 The gCLM equation and a fixed self-similar
profile}\label{the-gclm-equation-and-a-fixed-self-similar-profile}

The regularity of three-dimensional incompressible Euler turns on the
competition between vortex stretching and advection. The one-dimensional
caricature studied here is the generalized Constantin--Lax--Majda (gCLM)
equation \[
w_t+a\,u\,w_x=u_xw,\qquad u_x=Hw,
\] where \(H\) is the Hilbert transform (Fourier multiplier
\(-i\,\operatorname{sgn}\xi\) on \(\mathbb R\)); \(a=0\) is the original
Constantin--Lax--Majda model {[}1{]}, \(a=1\) is De Gregorio {[}2{]},
and the interpolation is due to Okamoto, Sakajo, and Wunsch {[}3{]}.

Throughout, the equation is posed on the whole line \(x\in\mathbb R\).
In the odd sector we fix the velocity primitive by \(u(0,t)=0\),
equivalently \[
u(x,t)=\int_0^x Hw(s,t)\,ds.
\] Then \(Hw\) is even, \(u\) is odd, and the odd sector is invariant.
Identifying an odd function with its restriction to \(x>0\) gives the
folded real-line Hilbert operators of Section 2.1. The far-field threshold
analyzed here comes from \(y\to\infty\), so the unbounded similarity
domain is intrinsic to the resulting two-ended Fredholm geometry.

Throughout the theorem, \(0<a<1\) is fixed and so is one smooth odd
self-similar collapse profile \((\Omega,c_l)\) satisfying the hypotheses in Section
2.2. It generates \[
w(x,t)=(T-t)^{-1}\Omega(y),\qquad y=\frac{x}{(T-t)^{c_l}},
\] and solves \[
\Omega+(c_ly+aU)\Omega'=\Omega H\Omega,\qquad U'=H\Omega.
\] Huang, Qin, Wang, and Wei {[}5, Theorem 3.11{]} prove fixed-point
existence at each parameter in their stated range; uniqueness is open
for general \(a\), and parameter continuity is posed in
{[}5, Conjecture 3.12{]}. Appendix
C works with the entire fixed-point set: every fixed point with
\(0<a<400/(848-9\pi^2)\), after normalization, satisfies the full
registry used here, and for all sufficiently small positive \(a\) the
far-field line lies strictly to the left of the origin line
for every such fixed point. The exact \(a=0\) profile
\(\Omega=-y/(y^2+\tfrac14)\), \(H\Omega=\tfrac12/(y^2+\tfrac14)\),
\(c_l=1\) serves as the benchmark.

Linearizing the self-similar dynamics about the fixed profile produces,
on the odd half-line \((0,\infty)\), \[
L_a\varphi=-\varphi-b\varphi'+q\varphi+\Omega H_{\mathrm{odd}}\varphi-aV(\varphi)\Omega',
\qquad b=c_ly+aU,\quad q=H\Omega,
\] where \(U(y)=\int_0^yq(s)\,ds\) and
\(V(\varphi)(y)=\int_0^y(H_{\mathrm{odd}}\varphi)(s)\,ds\). This is a first-order
transport operator with \(b>0\), plus a multiplier and two nonlocal
terms. We write
\(L_{\mathrm{loc}}\varphi:=-\varphi-b\varphi'+q\varphi\)
for its local transport part. To keep the fixed-profile formulas legible, we suppress the
\(a\)-subscript on profile quantities.

\subsection{\texorpdfstring{1.2 The endpoint picture and the
origin-\(H^2\)
realization}{1.2 The endpoint picture and the origin-H\^{}2 realization}}\label{the-endpoint-picture-and-the-origin-h2-realization}

The spectral analysis is realization-dependent. We work on \[
X=\{\varphi\in L^2(0,\infty):\varphi(0)=0,\ \varphi''\in L^2(0,\infty)\},
\qquad \|\varphi\|_X=\bigl(\|\varphi\|_2^2+\|\varphi''\|_2^2\bigr)^{1/2}.
\] The second-derivative control at the collapse point removes the
strip-filling behavior of the maximal plain-\(L^2\) realization. In the
companion paper {[}6{]}, we prove this realization dichotomy and obtain
the rational \(a=0\) benchmarks used here: the single essential line in
\(\{\operatorname{Re}\lambda\ge-1/2\}\) and the full odd point spectrum
\(\{0,1\}\).

For the fixed profile considered here, the two frozen endpoint models
instead select \[
\Gamma_\infty=\{\operatorname{Re}\lambda=-1+c_l/2\},\qquad
\Gamma_0=\{\operatorname{Re}\lambda=\lambda_X\},\qquad
\lambda_X=-1+q(0)-\tfrac32\tilde c=-\tfrac12\tilde c.
\] Section 3 proves directly that both lines belong to the essential
spectrum of the fixed operator. The task of the present paper is to
prove Fredholmness everywhere off them and compute the intervening
index.

Two-ended Fredholm problems are classically organized by limit-operator
theory {[}7{]}, weighted and conic analysis on noncompact spaces
{[}8--11{]}, and singular-integral methods for fixed singularities
{[}12{]}. Three features require a direct argument here: the Hilbert
multiplier is noncompact, its output does not preserve support, and the
origin trace changes with the operator realization. The proof therefore
supplies an explicit endpoint-subordination estimate, routes the trace
channel through a rank-one cancellation, and continues a common-domain
upper-semi-Fredholm path from a directly solved transport anchor.

Nonlinear stability of self-similar blow-up has been established for
related nonlocal transport equations in perturbative regimes {[}14{]}.
The present paper concerns instead the Fredholm and essential-spectrum
geometry of one fixed linearization. The registry below addresses the
smooth one-scale origin-\(H^2\) setting; multiscale CLM blow-up and
singular-profile gCLM regimes are treated in {[}15,16{]}.

\subsection{1.3 Goal of this paper}\label{goal-of-this-paper}

The companion paper's admissibility package \(\mathrm{Adm}(a)\)
{[}6, Definition 4.2{]} collects its profile equation and regularity,
origin compatibility, far-field decay, and nondegenerate transport
assumptions on the range \(0<a<a_c\). Here \(a_c\) denotes the
transition at which \(c_l\) changes sign and the self-similar spatial
extent changes from shrinking (focusing) to expanding. For profiles
satisfying that package, {[}6, Section 4.5{]} asked, under its Browder
convention, whether the essential spectrum on the origin-\(H^2\)
realization consists only of the two endpoint lines. Put \[
\begin{aligned}
F&=-1+c_l/2,\qquad O=-\tilde c/2,\\
S&:=\{\lambda:\min(F,O)<\operatorname{Re}\lambda<\max(F,O)\},\\
\overline S&:=\{\lambda:\min(F,O)\le\operatorname{Re}\lambda\le\max(F,O)\}.
\end{aligned}
\] We call \(S\) the open inter-line band, namely the open vertical strip
between the two endpoint lines, and \(\overline S\) its closed band.
Under the fixed-profile hypotheses of Section 2.2, the present
theorem applies. With the Edmunds--Evans conventions of Definition~2.3,
it proves \[
\sigma_{e4}(L_a|_X)=\sigma_{e5}(L_a|_X)
=\sigma_{\mathrm{B}}(L_a|_X)=\overline S,
\] while the first three Edmunds--Evans spectra are exactly the two
boundary lines. Thus, for every separated-line profile satisfying both
hypothesis packages, the question in {[}6, Section 4.5{]} has a negative
answer under the Browder convention used there: the Browder essential
spectrum is the closed strip \(\overline S\), not only its two boundary
lines. The fixed-profile equalities and the Huang--Qin--Wang--Wei
realizations below use the present package alone.

\begin{quote}
\noindent\textbf{Main result (fixed-profile overview).} Fix \(0<a<1\) and one
smooth odd self-similar collapse profile satisfying the
profile-and-focusing-transport hypotheses \((H_{\mathrm{prof}})\), the
boundedness-and-weighted-derivative package \((K4^+)\), and the
far-field-data hypotheses \((D_\infty)\), all listed in Section~2.2. For
every \(\lambda\) off the two endpoint lines, the operator \(L_a-\lambda\)
on the maximal domain \(D_a\) of Definition~2.1 is Fredholm
and \[
\operatorname{ind}(L_a-\lambda)=p_0+p_\infty-1,
\] where \(p_0,p_\infty\in\{0,1\}\) are the endpoint-admissibility flags
defined in Section~5.2. Consequently
\(\sigma_{e1}=\sigma_{e2}=\sigma_{e3}=\Gamma_0\cup\Gamma_\infty\),
whereas \(\sigma_{e4}=\sigma_{e5}=\overline S\). Every point of the open
band \(S\), when its index is \(+1\), lies in the point spectrum. Condition
(ND) enters only the later
kernel--cokernel multiplicity refinement.
\end{quote}

The proof first establishes a three-zone upper-semi-Fredholm estimate
off both lines. It then uses the common-domain path \[
L_a^t=L_{\mathrm{loc}}+t\big(\Omega H_{\mathrm{odd}}-aV(\,\cdot\,)\Omega'\big),\qquad 0\le t\le1,
\] which is norm-continuous in the common graph norm and remains in the
upper-semi-Fredholm class \(\Phi_+\) of Definition~2.3. The direct first-order anchor for
\(L_{\mathrm{loc}}-\lambda\) supplies Fredholmness and the endpoint
index, and continuation inside \(\Phi_+\) carries both to
\(L_a-\lambda\). Sections 4 and 5 form the Fredholm argument.

The origin profile equation and the gauge \(\Omega'(0)=-4\) give the
exact pointwise identity \[
q(0)-\tilde c=1,\qquad \tilde c:=c_l+a q(0).
\] Solving these two scalar identities gives \[
q(0)=\frac{1+c_l}{1-a},\qquad
\tilde c=\frac{c_l+a}{1-a}. \tag{1.1}
\] Since \(c_l>0\) and \(\tilde c=\lim_{y\downarrow0}b(y)/y>0\), the
fixed-profile hypotheses are inconsistent at \(a=1\) and for \(a>1\).
Thus every profile covered by this paper necessarily has \(0<a<1\), and
\(\lambda_X=-\tilde c/2\). Equations (1.1) therefore restrict the
registry to \(0<a<1\); Propositions 2.7--2.8 and Appendix C supply
positive-\(a\) nonemptiness and separated-line realizations.

For this fixed profile, if \[
g:=\frac{c_l+\tilde c}{2}-1,
\] then (1.1) gives the profile-local orientation criterion \[
g=\frac{(2-a)c_l+3a-2}{2(1-a)},\qquad
F<O\iff c_l<\frac{2-3a}{2-a},\qquad
F=O\iff c_l=\frac{2-3a}{2-a}. \tag{1.2}
\] In particular, every focusing registry profile with \(a\ge2/3\) has
the reverse orientation \(O<F\). Existence in this parameter range is a
separate profile question.

With this notation, \(\lambda_X-(-1+c_l/2)=-g\) and the right essential
edge is \[
\max(-1+c_l/2,\lambda_X)=-1+c_l/2+\max(0,-g).
\] Thus the corresponding distance of the essential spectrum from the
imaginary axis is the fixed-profile quantity \(1-c_l/2-\max(0,-g)\).

\subsection{1.4 Logical layers}\label{logical-layers}

\noindent\textbf{(I) Fixed-profile essential-spectrum theorem.} Under the
profile-and-focusing-transport hypotheses \((H_{\mathrm{prof}})\), the
boundedness-and-weighted-derivative package \((K4^+)\), and the
far-field-data hypotheses \((D_\infty)\) for the chosen profile, the two-line exactness of
\(\sigma_{e1},\sigma_{e2},\sigma_{e3}\), the exact closed-strip formula
for \(\sigma_{e4}=\sigma_{e5}\) and the Browder essential spectrum,
Fredholmness off the lines, and the index \(p_0+p_\infty-1\) are analytic
results. The inclusion of the lines comes from Section 3. The
upper-semi-Fredholm estimate of Section 4, its common-domain \(\Phi_+\)
homotopy, and the direct local anchor and continuation argument of
Section 5 give the reverse inclusion and full Fredholmness. The
construction is fixed-profile and primal; the merged-line obstruction
uses a Banach-adjoint singular sequence.

\noindent\textbf{(II) Actual Huang--Qin--Wang--Wei realizations.} Proposition 2.7 and Appendix C
prove the complete fixed-profile hypothesis registry for every
Huang--Qin--Wang--Wei fixed point with \(0<a<400/(848-9\pi^2)\).
Proposition 2.8 proves that all such fixed points have \(F<O\) for every
sufficiently small positive \(a\). Together with {[}5, Theorem 3.11{]},
this supplies actual collapse profiles with a nonempty band of Fredholm index \(+1\)
throughout the fixed-point set, with a qualitative small-\(a\)
threshold.

The realizations proved here are the small-\(a\) family realizing the band
of Fredholm index \(+1\) and
the exact merged-line witness at \(a=1/2\). Reverse-orientation
conclusions apply to any supplied registry profile with \(O<F\).

\noindent\textbf{(III) Spectral-type refinement.} The index does not by itself
determine the individual kernel and cokernel dimensions. Positive index
already forces point spectrum. Lemma 5.1 proves condition (ND) outside
the energy slab, Proposition 5.6 proves it on every high-frequency
closed inter-line sub-band a positive distance from the far-field line, and
Corollary 5.7 makes its possible failure set locally finite in the open
band. This is an interior statement: it neither excludes accumulation at
a boundary line nor resolves the simultaneous limit
\(\operatorname{Re}\lambda\to F\) and
\(|\operatorname{Im}\lambda|\to\infty\). At that corner the far-end
transport rate \(|F-\operatorname{Re}\lambda|\) tends to zero, so the
inverse bounds used to absorb the \(O(|\operatorname{Im}\lambda|^{-1})\)
Hilbert--Volterra remainder lose uniformity. Condition (ND) is therefore
an independent kernel--cokernel multiplicity refinement.

\section{2. Setting and realization}\label{setting-and-realization}

\subsection{2.1 The odd sector and the Hilbert
pair}\label{the-odd-sector-and-the-hilbert-pair}

On the half-line we use the folded Hilbert transforms \[
H_{\mathrm{odd}}\varphi(y) = \frac1\pi\,\mathrm{p.v.}\!\int_0^\infty\!\varphi(t)\Big[\frac{1}{y-t}-\frac{1}{y+t}\Big]dt
= \frac2\pi\,\mathrm{p.v.}\!\int_0^\infty\!\frac{t\,\varphi(t)}{y^2-t^2}\,dt,
\] the restriction to \((0,\infty)\) of \(H\) applied to the odd
extension, and \(H_{\mathrm{ev}}\) with the \(+\) kernel for even
functions. Both are bounded on \(L^2(0,\infty)\) with norm \(\le 1\);
\(H_{\mathrm{odd}}\) maps odd to even, \(H_{\mathrm{ev}}\) even to odd,
and \(H\) commutes with \(d/dy\). Two consequences are used throughout:
if \(\varphi(0)=0\) and \(\varphi'\in L^2_{\mathrm{loc}}\) then
\((H_{\mathrm{odd}}\varphi)' = H_{\mathrm{ev}}(\varphi')\), and if in
addition \(\varphi''\in L^2\) then
\((H_{\mathrm{odd}}\varphi)'' = H_{\mathrm{odd}}(\varphi'')\in L^2_{\mathrm{loc}}\).
The antiderivative
\(V(\varphi)(y)=\int_0^y(H_{\mathrm{odd}}\varphi)(s)\,ds\)
satisfies \(V(0)=0\), \(V' = H_{\mathrm{odd}}\varphi\),
\(|V(y)|\le y^{1/2}\|\varphi\|_2\); it is an order \(-1\) (smoothing)
operator.

\subsection{2.2 The fixed-profile data
registry}\label{the-fixed-profile-data-registry}

Fix \(0<a<1\) and one profile \((\Omega,c_l)\). All quantities below
belong to this fixed profile. Since \(a\) and the profile are fixed, all
profile quantities below are written without an \(a\)-subscript.
Registry A.4 provides a tabular summary and use map. Their finiteness is
a standing hypothesis of the theorem. Our notation differs slightly
from {[}6{]}: the transport weight \(c_ly+aU\), denoted by \(q\) in
{[}6, Definition 4.2(H4){]}, is denoted by \(b\) here, while
\(q:=H\Omega\).

\noindent\textbf{\((H_{\mathrm{prof}})\) (profile and focusing transport).} The
pair is a smooth odd self-similar collapse solution of \[
\Omega+b\Omega'=\Omega q,\qquad q=H\Omega,\qquad b(y)=c_ly+aU(y),\qquad U(y)=\int_0^yq(s)\,ds,
\] with \(0<c_l<2\), \(\Omega\in L^2\), the normalization
\(\Omega'(0)=-4\), and \[
0<\beta_-\le \frac{b(y)}{y}\le\beta_+<\infty\qquad(y>0).
\] The required origin expansion is available on some radius
\(r_{\mathrm{org}}>0\) and
has finite constants \[
M_1=\sup_{(0,r_{\mathrm{org}})}|q'|,\quad
M_2=\sup_{(0,r_{\mathrm{org}})}|\Omega''|,\quad
M_3=\sup_{(0,r_{\mathrm{org}})}\frac{|\Omega(y)-\Omega'(0)y|}{y^3},\quad
M_4=\sup_{(0,r_{\mathrm{org}})}|q''|.
\] We also assume the exact origin identity \(q(0)-\tilde c=1\),
where \(\tilde c=c_l+aq(0)\). Consequently (1.1) holds, the standing
hypotheses force \(0<a<1\), and \(\lambda_X=-\tilde c/2\).

\noindent\textbf{\((K4^+)\) (boundedness and weighted derivatives).} Here \(K4^+\)
labels the following global boundedness, integrability, and
weighted-derivative finiteness conditions: \[
Q_0=\|q\|_\infty,\quad Q_1=\|q'\|_\infty,\quad
W_j=\|\Omega^{(j)}\|_\infty\ (j=0,1,2),\quad I_q=\int_0^\infty|q(y)|\,dy,
\] \[
J_k=\left(\int_0^\infty y|\Omega^{(k)}(y)|^2\,dy\right)^{1/2}<\infty
\quad(k=1,2,3),
\] together with the weighted tails implicit in these finite integrals.
The tail moduli used below are defined explicitly after \((D_\infty)\);
Registry A.4 records where each enters.

\noindent\textbf{\((D_\infty)\) (far-field data).} For some \(C_q,C_\Omega>0\)
and some decay exponent \(\alpha_{\rm dec}>1\), \[
|q(y)|\le C_q(1+y)^{-\min(\alpha_{\rm dec},2)},\qquad
|\Omega(y)|\le C_\Omega(1+y)^{-\alpha_{\rm dec}},\qquad U(y)\to U_\infty,
\] with \(b(y)=c_ly+O(1)\). Together with \((K4^+)\), define \[
\varepsilon_\Omega(R):=\sup_{y\ge R}\bigl(|\Omega(y)|+|\Omega'(y)|+|\Omega''(y)|\bigr),\qquad
\varepsilon_q(R):=\sup_{y\ge R}|q(y)|,
\] \[
\varepsilon_V(R):=
\left(\int_R^\infty y\bigl(|\Omega'(y)|^2+|\Omega'''(y)|^2\bigr)\,dy\right)^{1/2}.
\] All three tend to zero. For the derivative suprema this follows, for
\(f=\Omega'\) and \(f=\Omega''\), from
\(\sup_{y\ge R}|f(y)|^2\le2\|f\|_{L^2(R,\infty)}\|f'\|_{L^2(R,\infty)}\)
and the weighted \(J_1,J_2,J_3\) tails. Hence the displayed \((K4^+)\)
and \((D_\infty)\) bounds already imply all differentiated tail
estimates used below. The fixed-profile theorem uses these bounds for
some \(\alpha_{\rm dec}>1\). For the Huang--Qin--Wang--Wei profiles in
Appendix C, the verified bookkeeping exponent is \(\alpha_{\rm reg}\)
of Section C.3; the raw asymptotic exponent
\(s_{\mathrm{prof}}:=1/c_l\) is recorded separately in
(C.5).

The registry remains a fixed-profile bookkeeping device. Appendix C
verifies \((H_{\mathrm{prof}})+(K4^+)+(D_\infty)\) for every
Huang--Qin--Wang--Wei fixed point with
\(0<a<\underline a:=400/(848-9\pi^2)\), after the normalization used
here; {[}5, Theorem 3.11{]} supplies at least one such fixed point at
every parameter in that range. Other profiles require individual
verification.

The companion's four-part admissibility package \(\mathrm{Adm}(a)\)
{[}6, Definition 4.2{]} covers the profile equation and regularity,
origin compatibility, far-field decay, and nondegenerate transport on
its range \(0<a<a_c\). It is distinct from
\((H_{\mathrm{prof}})+(K4^+)+(D_\infty)\). Comparisons with {[}6{]} are
made on the intersection of the two packages. The explicit \(a=1/2\)
profile in Section C.5 belongs to that intersection, but its endpoint
lines coincide. For general Huang--Qin--Wang--Wei fixed points, Appendix C proves
the zeroth-order far-field rate, all-order Sobolev regularity, weighted
\(L^2\) derivative control, and derivative tails without a pointwise
rate. It does not prove the sharper pointwise rates
\[
|\Omega^{(j)}(y)|\le C y^{-1/c_l-j},\qquad 1\le j\le3,
\]
required by {[}6, Definition 4.2(H3){]}. Thus nonemptiness of the
separated-line part of the intersection is not asserted here. Conditional
on a separated-line profile in that intersection, Theorem 7.1 gives the
negative answer under the Browder convention of {[}6{]}; independently,
the present fixed-profile theorem and its actual small-\(a\)
Huang--Qin--Wang--Wei realizations are nonvacuous under the present
registry.

\subsection{2.3 The operator and its
domain}\label{the-operator-and-its-domain}

\noindent\textbf{Definition 2.1.} \(L_a\) acts on \(X\) by the formula of Section
1.1 with maximal domain
\(D(L_a)=D_a:=\{\varphi\in X:L_a\varphi\in X\}\). The
zeroth-order terms are bounded on \(L^2\)
(\(\|q\varphi\|_2\le Q_0\|\varphi\|_2\),
\(\|\Omega H_{\mathrm{odd}}\varphi\|_2\le W_0\|\varphi\|_2\),
\(\|aV(\varphi)\Omega'\|_2\le aJ_1\|\varphi\|_2\)), so
\(L_a\varphi\in L^2\) iff \(y\varphi'\in L^2\).

\noindent\textbf{Proposition 2.2 (closed realization and trace splitting).} The
maximal realization \(L_a:D_a\subset X\to X\) is closed. The trace
\(\tau\varphi:=\varphi'(0)\) is bounded on \(X\). With a fixed smooth
cutoff \(\chi_0\) (\(\equiv1\) on \([0,1]\), supported in \([0,2]\)) and
\(e(y):=y\chi_0(y)\), one has \(e\in D_a\), \(\tau e=1\), and \[
X=X_w\oplus\mathbb C e,\qquad D_a=(D_a\cap X_w)\oplus\mathbb C e,\qquad X_w=X_{00}:=\ker\tau,
\] with bounded projections. Moreover
\(C_c^\infty(0,\infty)\subset D_a\) is dense in \(X_w\); hence
\(C_c^\infty(0,\infty)\oplus\mathbb C e\) is dense in \(X\) and \(D_a\)
is densely defined. We do \textbf{not} assert that
\(C_c^\infty(0,\infty)\) alone is dense in \(X\) or is a graph core for
\(D_a\). The \(y\)-mode is the rank-one bookkeeping piece used below.
For \(e=y\chi_0\), the folded transform has
\(H_{\mathrm{odd}}e=O(y^{-2})\) and \(V(e)=O(1)\) at infinity; the
\((K4^+)\) bounds then give \(L_ae\in X\).

\emph{Proof.} The norm in Section 1.2 is equivalent on \(X\) to the
Hilbert norm \((\|\varphi\|_2^2+\|\varphi''\|_2^2)^{1/2}\). Indeed, \[
\|\varphi'\|_2^2\le \|\varphi\|_2\|\varphi''\|_2, \tag{2.2a}
\] obtained by applying the Fourier interpolation inequality to the odd
extension of \(\varphi\) (the condition \(\varphi(0)=0\) creates no
boundary mass). Thus \(X\) is the closed trace-zero subspace of
\(H^2(0,\infty)\), and the trace theorem gives \[
|\tau\varphi|\le C_{\mathrm{tr}}(\|\varphi\|_2+\|\varphi''\|_2). \tag{2.2b}
\] Write \[
B\varphi:=(-1+q)\varphi+\Omega H_{\mathrm{odd}}\varphi-aV(\varphi)\Omega',
\qquad L_a\varphi=-b\varphi'+B\varphi.
\] The registry and \(|V(\varphi)(y)|\le y^{1/2}\|\varphi\|_2\) give \[
\|B\varphi\|_2\le(1+Q_0+W_0+aJ_1)\|\varphi\|_2. \tag{2.2c}
\] Suppose \(\varphi_n\in D_a\), \(\varphi_n\to\varphi\) in \(X\), and
\(L_a\varphi_n\to g\) in \(X\). Then (2.2c) gives \[
b\varphi_n'=B\varphi_n-L_a\varphi_n\longrightarrow B\varphi-g
\quad\hbox{in }L^2.
\] On every compact subinterval, (2.2a) gives \(\varphi_n'\to\varphi'\)
in \(L^2\) and \(b\) is bounded, so the same sequence converges
distributionally to \(b\varphi'\). Hence
\(b\varphi'=B\varphi-g\in L^2\). Since \(\beta_-\le b(y)/y\le\beta_+\),
one has \(y\varphi'\in L^2\), and therefore
\(L_a\varphi=-b\varphi'+B\varphi=g\in X\). Thus \(\varphi\in D_a\) and
the maximal realization is closed.

For \(e=y\chi_0\), the local and transport terms belong to \(X\). Put
\(h=H_{\mathrm{odd}}e\). Boundedness of the folded transforms and
commutation with derivatives give \(h,h',h''\in L^2\); compact support
of \(e\) gives \(h^{(j)}(y)=O(y^{-2-j})\) for \(j=0,1,2\), and hence
\(V(e)=O(1)\) at infinity and \(V(e)=O(y)\) at zero. Consequently \[
(\Omega h)''=\Omega''h+2\Omega'h'+\Omega h'',
\] \[
(V(e)\Omega')''=h'\Omega'+2h\Omega''+V(e)\Omega'''
\] belong to \(L^2\): use \(W_0,W_1,W_2\) for the bounded factors and
\(J_3\) for the last term. Also \(q''=H\Omega''\in L^2\), so the
remaining compactly supported product terms lie in \(X\). Hence
\(e\in D_a\) and \(\tau e=1\).

The bounded projection \(P\varphi=(\tau\varphi)e\) now gives
\(X=\ker\tau\oplus\mathbb Ce\). Since \(D_a\) is linear and contains
\(e\), its restriction gives \(D_a=(D_a\cap\ker\tau)\oplus\mathbb Ce\).
Both projections are bounded on \(X\); on \(D_a\) with its graph norm
this follows from (2.2b) and \(e\in D_a\). The same nonlocal calculation
shows \(C_c^\infty(0,\infty)\subset D_a\). Lemma A.1 proves that this
class is dense in \(\ker\tau\); adding \(\mathbb Ce\) proves density in
\(X\) and completes the proof. \(\square\)

\emph{Remark (realization).} \(X\) is a Hilbert space with
\(\varphi\in H^2_{\mathrm{loc}}([0,\infty))\), so \(\varphi,\varphi'\)
are continuous and \(\varphi'(0)\) is well defined; the only genuine
origin condition is \(\varphi(0)=0\). The domain requirement
``\(\varphi\sim y\)'' of the companion means exactly ``\(\varphi'(0)\)
exists,'' which a closed operator domain, being a linear subspace,
includes the case \(\varphi'(0)=0\).

\noindent\textbf{Definition 2.3 (essential-spectrum conventions).} For
\(T_\lambda:=L_a-\lambda:D_a\to X\), viewed as a bounded operator from
the graph space \(D_a\) to \(X\), write \[
\Phi_+:=\{T:\operatorname{ran}T\text{ is closed},\ \dim\ker T<\infty\},
\] \[
\Phi_-:=\{T:\operatorname{ran}T\text{ is closed},\quad
\operatorname{codim}\operatorname{ran}T<\infty\},\qquad
\Phi:=\Phi_+\cap\Phi_-.
\] We use the Edmunds--Evans ordering {[}13, Chapter I, Section 4{]} \[
\sigma_{e1}(L_a):=\{\lambda:T_\lambda\notin\Phi_+\cup\Phi_-\},\qquad
\sigma_{e2}(L_a):=\{\lambda:T_\lambda\notin\Phi_+\},
\] \[
\sigma_{e3}(L_a):=\{\lambda:T_\lambda\notin\Phi\},\qquad
\sigma_{e4}(L_a):=\{\lambda:T_\lambda\notin\Phi\ \text{or}\ \operatorname{ind}T_\lambda\ne0\}.
\] If \(\lambda_0\) is isolated in \(\sigma(L_a)\), let \[
P_{\lambda_0}:=\frac{1}{2\pi i}\int_\gamma(\zeta-L_a)^{-1}\,d\zeta,
\] where \(\gamma\subset\rho(L_a)\) is a positively oriented circle
enclosing only \(\lambda_0\). Define
{[}17, Chapter III, Sections 6.4--6.5{]} \[
\sigma_{\mathrm{disc}}(L_a):=
\{\lambda_0\in\sigma(L_a):\lambda_0\text{ is isolated and }
0<\operatorname{rank}P_{\lambda_0}<\infty\},
\] and call \(\operatorname{rank}P_{\lambda_0}\) its algebraic
multiplicity. Set \[
\sigma_{e5}(L_a):=\sigma(L_a)\setminus\sigma_{\mathrm{disc}}(L_a).
\] We also use the Browder essential spectrum in its operator-theoretic
sense: \[
\sigma_{\mathrm{B}}(L_a):=
\{\lambda\in\mathbb C:T_\lambda\text{ is not Fredholm with finite ascent and finite descent}\}.
\] Here the ascent is the first index at which the kernels of the powers
of \(T_\lambda\) stabilize, and the descent is defined analogously from
the ranges. For a closed operator, finite ascent and descent at a Fredholm
spectral point are equivalent to a resolvent pole with finite-rank Riesz
projection. Thus the spectral points outside \(\sigma_{\mathrm{B}}\) are
exactly the isolated eigenvalues with finite-rank Riesz projection
{[}17, Chapter III, Sections 6.4--6.5{]}.
Consequently, \[
\sigma_{\mathrm{B}}(L_a)
=\sigma(L_a)\setminus\sigma_{\mathrm{disc}}(L_a)
=\sigma_{e5}(L_a). \tag{2.3a}
\] Thus the symbol \(\sigma_{\mathrm{ess}}\) used in {[}6{]} corresponds
to \(\sigma_{\mathrm B}=\sigma_{e5}\) in the present notation. Moreover,
\[
\sigma_{e1}\subseteq\sigma_{e2}\subseteq\sigma_{e3}
\subseteq\sigma_{e4}\subseteq\sigma_{e5}.
\] The last definition is the equivalent discrete-spectrum
characterization of the fifth Edmunds--Evans essential spectrum. In the
original component formulation, \(\rho_{e5}(L_a)\) is the union of those
connected components of the semi-Fredholm set
\(\rho_{e1}(L_a):=\mathbb C\setminus\sigma_{e1}(L_a)\) that meet the
resolvent \(\rho(L_a)\); for closed operators this is equivalent to
\(\sigma_{e5}(L_a)=\sigma(L_a)\setminus\sigma_{\mathrm{disc}}(L_a)\)
{[}13, Chapter I, Section 4{]}. These definitions fix the numbering used
below. The inter-line band is an open subset of the spectrum whenever
the two lines are distinct, because its Fredholm index is nonzero; hence
it belongs to both \(\sigma_{e4}\) and \(\sigma_{e5}\) even though it is
disjoint from \(\sigma_{e3}\). Theorem 7.1(4), using Sections 5.3 and 6,
makes the component picture explicit: both exterior components meet the
resolvent, whereas a nonempty inter-line component has nonzero index and
cannot do so.

\subsection{2.4 The first-order graph
norm}\label{the-first-order-graph-norm}

The estimates of Section 4 are run in a weighted first-order graph norm
adapted to the two endpoints. With the weight \(w(y) = y^{-2}\) for
\(y\le 1\), \(w=1\) for \(y\ge 1\), and the decomposition
\(\varphi = \psi + l\,y\chi_0\), \(l=\varphi'(0)\), so that
\(\psi(0)=\psi'(0)=0\), set \[
N(\varphi) := \|w\psi\|_2 + \|w\,b\psi'\|_2 + |l|.
\] Hardy at the origin and the bounds on \(b/y\) at infinity give, for
\(\varphi\in D_a\), \[
N(\varphi)\le C_{\mathrm{int}}
\bigl(\|\varphi\|_X+\|y\varphi'\|_2\bigr). \tag{2.4}
\] The far-field derivative term in (2.4) is essential and is not
controlled by \(\|\varphi\|_X\) alone. Conversely,
\(\|\varphi\|_2\le C_N\,N(\varphi)\) and
\(\|\psi'\|_2\le C_b\|wb\psi'\|_2\), with constants from the registry
(\(\beta_\pm\), \(c_{\min}=\inf_{0<y\le1}b/y>0\),
\(\kappa_\infty=\inf_{y\ge1}b/(c_ly)>0\)). For every
\(\lambda\in\mathbb C\) and \(\varphi\in D_a\), the undifferentiated
equation and \(b/y\ge\beta_-\) give
\[
\|y\varphi'\|_2\le\beta_-^{-1}
\bigl[(1+|\lambda|+Q_0+W_0+aJ_1)\|\varphi\|_2
+\|(L_a-\lambda)\varphi\|_2\bigr].
\]
Thus the weighted first derivative is graph-bounded uniformly on compact
spectral sets and closes (2.4). Separation from the origin line enters
only in Theorem 2.6, when the second derivative is recovered.

\subsection{2.5 The bootstrap lemma}\label{the-bootstrap-lemma}

The graph norm \(N\) controls only the first derivative; the
second-derivative content of \(\|\varphi\|_X\) is recovered from it by a
one-shot bootstrap that needs no relation to the far-field line.

\noindent\textbf{Lemma 2.5 (compatibility cancellation).} Let \(\varphi\in D_a\)
and \((L_a-\lambda)\varphi=f\). Put \[
z=\varphi',\qquad p=z(0),\qquad
h=H_{\mathrm{odd}}\varphi,\qquad h_0=h(0),
\] and \[
\kappa=1+\lambda+b'-q,
\] \[
\mathfrak r=q'\varphi+(1-a)\Omega'h+\Omega H_{\mathrm{ev}}z
-aV(\varphi)\Omega''-f'.
\] Then \(\mathfrak r\) has the normalized local \(L^2\) trace \[
(1-a)\Omega'(0)h_0-f'(0)=\lambda p=\tilde c\,\nu p,
\qquad \nu=\lambda/\tilde c.
\] Here and below, \(\mathfrak r(0)\) denotes this normalized local
\(L^2\) trace, not a separately assumed pointwise value. In particular,
this is the compatibility value that removes the apparent \(1/y\)
singularity in the equation for the trace-free part of \(z\).

\emph{Proof.} We say that \(g\) has normalized local \(L^2\) trace \(c\)
if \[
\delta^{-1/2}\|g-c\|_{L^2(0,\delta)}\longrightarrow0
\qquad(\delta\downarrow0).
\] Since \(\varphi,f\in X\), the functions \(z\) and \(f'\) are locally
\(H^1\) and have traces \(p\) and \(f'(0)\). Moreover,
\(h'=H_{\mathrm{ev}}z\in L^2\), so \(h\) has trace \(h_0\). The Taylor
bounds at the origin now show, term by term, that \(\mathfrak r\) has
normalized trace \((1-a)\Omega'(0)h_0-f'(0)\). Indeed, \(q'\varphi\) has trace
zero, \(\Omega'h\) has trace \(\Omega'(0)h_0\), and \[
\delta^{-1/2}\|\Omega H_{\mathrm{ev}}z\|_{L^2(0,\delta)}
\le C\delta^{1/2}\|H_{\mathrm{ev}}z\|_{L^2(0,\delta)}\longrightarrow0.
\] Also \(V(\varphi)(0)=0\) and
\(|V(\varphi)(y)|\le y^{1/2}\|\varphi\|_2\), so the
\(V(\varphi)\Omega''\) term has trace zero.

Differentiating the equation in distributions gives \[
bz'+\kappa z=\mathfrak r.
\] Because \(z'=\varphi''\in L^2\) and \(b(y)=O(y)\), \[
\delta^{-1/2}\|bz'\|_{L^2(0,\delta)}
\le C\delta^{1/2}\|z'\|_{L^2(0,\delta)}\longrightarrow0.
\] The normalized trace of \(\kappa z\) is \(\kappa(0)p\). Comparison in
the last displayed equation therefore gives
\((1-a)\Omega'(0)h_0-f'(0)=\kappa(0)p\). Finally, \(b'(0)=\tilde c\) and
\(q(0)-\tilde c=1\), hence
\(\kappa(0)=1+\lambda+\tilde c-q(0)=\lambda\). \(\square\)

\noindent\textbf{Theorem 2.6 (bootstrap, fixed profile).} Let
\(K\subset\mathbb C\) be compact with
\(d_X(\lambda):=|\operatorname{Re}\lambda-\lambda_X|\ge d_0>0\) on
\(K\). There is \(C_K<\infty\), depending only on \(K\), \(d_0\), and
the fixed-profile registry, such that
\[
\|\varphi''\|_2 + \|y\varphi'\|_2 \ \le\ C_K\big(\|(L_a-\lambda)\varphi\|_X + \|\varphi\|_2\big)
\qquad(\lambda\in K,\ \varphi\in D_a).
\] No condition relative to the far-field line is needed. The constant must
deteriorate as \(d_0\downarrow0\): already the Euler-model factor below
is \(1+|\nu|/|\operatorname{Re}\nu+1/2|=1+|\lambda|/d_X(\lambda)\), and
the log-widening origin Weyl packets at exponent \(3/2\) (Section 3.3)
witness this loss. We make no sharper asymptotic claim for the full
registry-dependent constant. The weighted first-derivative part is the
unconditional estimate of Section 2.4; separation from the origin line
is needed only for the second-derivative estimate.

\emph{Proof.} Put \[
f=(L_a-\lambda)\varphi,\qquad z=\varphi',\qquad
p=z(0),\qquad h=H_{\mathrm{odd}}\varphi,\qquad
\nu=\lambda/\tilde c .
\] The undifferentiated equation and its derivative are \[
bz=-(1+\lambda)\varphi+q\varphi+\Omega h-aV(\varphi)\Omega'-f, \tag{2.5a}
\] \[
bz'+\kappa z=\mathfrak r,\qquad
\kappa:=1+\lambda+b'-q=1+\lambda+c_l+(a-1)q, \tag{2.5b}
\] where \[
\mathfrak r:=q'\varphi+(1-a)\Omega'h+\Omega H_{\mathrm{ev}}z
-aV(\varphi)\Omega''-f'. \tag{2.5c}
\] At the origin \(\kappa(0)=\lambda\), and Lemma 2.5 identifies the
normalized local \(L^2\) trace and gives \[
\mathfrak r(0)=(1-a)\Omega'(0)h(0)-f'(0)=\lambda p. \tag{2.5d}
\]

We use one scalar estimate. If \(u(0)=0\) and \(u'\in L^2(0,\infty)\),
then \[
\|u'\|_2\le C_\nu\left\|\frac{yu'+\nu u}{y}\right\|_2,\qquad
C_\nu:=1+\frac{|\nu|}{|\operatorname{Re}\nu+1/2|}. \tag{2.5e}
\] To see this, put \(v(t)=e^{-t/2}u(e^t)\) and take the Fourier
transform in \(t\); comparison of the multipliers \(i\xi+1/2\) and
\(i\xi+\nu+1/2\) gives (2.5e), first for a dense smooth class and then
by closure. Since
\(|\operatorname{Re}\nu+1/2|=d_X(\lambda)/\tilde c\), the number
\(C_E:=\sup_{\lambda\in K}C_\nu\) is finite.

Choose \(\delta\le r_{\mathrm{org}}/2\) so that the origin Taylor bounds from
\((H_{\mathrm{prof}})\) hold and the coefficient errors below, after
multiplication by \(C_E\), are at most \(1/16\). Take \(\chi=1\) on
\([0,\delta]\), supported in \([0,2\delta]\), with
\(|\chi'|\le2/\delta\), and set \[
\rho:=\min\{\delta,\delta/(64C_E^2)\},\qquad u=\chi(z-p).
\] Subtracting the origin identity from (2.5b) gives exactly \[
\frac{yu'+\nu u}{y}
=\chi'(z-p)+\frac{\chi(\mathfrak r-\mathfrak r(0))}{\tilde c\,y}
-\frac{\chi(\kappa-\lambda)z}{\tilde c\,y}
-\frac{\chi(b-\tilde c y)z'}{\tilde c\,y}. \tag{2.5f}
\] Every quotient is regular. Hardy's inequality, the \(L^2\) bounds of
the folded Hilbert transforms, and the registry give \[
\left\|\frac{\mathfrak r-\mathfrak r(0)}{y}\right\|_{L^2(0,2\delta)}
\le C_{\mathrm{reg}}\bigl(\|\varphi\|_2+\|z\|_2+\|f''\|_2\bigr). \tag{2.5g}
\] For completeness, decompose
\(\Omega'h-\Omega'(0)h(0)=\Omega'(0)(h-h(0))+(\Omega'-\Omega'(0))h\) and use
\(h-h(0)=\int_0^yH_{\mathrm{ev}}z(s)\,ds\). The other four differences
obey \[
\begin{aligned}
\left\|\frac{q'\varphi}{y}\right\|_2&\le2M_1\|z\|_2,
&\left\|\frac{\Omega H_{\mathrm{ev}}z}{y}\right\|_2&\le W_1\|z\|_2,\\
\left\|\frac{V(\varphi)\Omega''}{y}\right\|_2&\le2M_2\|\varphi\|_2,
&\left\|\frac{f'-f'(0)}{y}\right\|_2&\le2\|f''\|_2.
\end{aligned}
\] on \((0,2\delta)\), with harmless registry factors absorbed into
\(C_{\mathrm{reg}}\). Evenness of \(q\) and the \(M_4\) bound also give
\[
|\kappa(y)-\lambda|\le\frac{1-a}{2}M_4y^2,\qquad
|b(y)-\tilde c y|\le\frac{a}{6}M_4y^3. \tag{2.5h}
\]

To control the transition term, average
\(p=z(y)-\int_0^y z'(s)\,ds\) on \((0,\rho)\) to obtain \[
|p|\le \rho^{-1/2}\|z\|_{L^2(0,\rho)}
+\rho^{1/2}\|z'\|_{L^2(0,\rho)}. \tag{2.5i}
\] Consequently the \(p\chi'\) part of (2.5f), after (2.5e), contributes
at most \(2C_E(\rho/\delta)^{1/2}\|z'\|_{L^2(0,\rho)}
\le\tfrac14\|z'\|_{L^2(0,\rho)}\), plus \(C_{K,\delta}\|z\|_2\).

On \([\delta,\infty)\), \(b\) is bounded below by a positive multiple of
\(\min(1,y)\), using \(c_{\min}\) and \(\kappa_\infty\). Directly from
(2.5b)--(2.5c), \[
\|z'\|_{L^2(\delta,\infty)}
\le C_{K,\delta,\mathrm{reg}}
\bigl(\|\varphi\|_2+\|z\|_2+\|f'\|_2\bigr). \tag{2.5j}
\] Apply (2.5e) to (2.5f), use (2.5g)--(2.5j), and first choose
\(\delta\) as above. The \(z'\) coefficient from (2.5h) and the trace
coefficient from (2.5i) are absorbed, while the annulus
\((\delta,2\delta)\) is controlled by (2.5j). Thus \[
\|z'\|_2\le C_K\bigl(\|f\|_X+\|\varphi\|_2+\|z\|_2\bigr). \tag{2.5k}
\] The interpolation inequality
\(\|z\|_2\le\varepsilon\|z'\|_2+C_\varepsilon\|\varphi\|_2\) absorbs the
last term. The first-derivative estimate of Section 2.4 supplies the
remaining \(\|yz\|_2\) bound. All identities above hold distributionally
for \(\varphi\in D_a\);
since the scalar localization belongs to the local \(H^1\) class, the
density passage is internal to (2.5e). This proves the theorem.
\(\square\)

As a corollary, on \(D_a\) the norm
\(\|\varphi\|_X + \|(L_a-\lambda)\varphi\|_X\) is equivalent to
\(\|\varphi\|_2 + \|\varphi''\|_2 + \|y\varphi'\|_2 + \|(L_a-\lambda)\varphi\|_X\),
so every localization and commutator error of the form
\(\|\varphi\|_2\), \(\|\varphi'\|_2\),
\(\|\chi\varphi\|_{H^1(\mathrm{cpt})}\) in Section 4 is graph-bounded
with explicit constants.

\subsection{2.6 Actual Huang--Qin--Wang--Wei
realizations}\label{actual-hqww-realizations}

Let \(f\in\mathbb D\) be a fixed point of the Huang--Qin--Wang--Wei map
\(R_a\), with raw profile \(\omega_0(x)=-xf(x)\) and raw similarity
constants \(c_{\omega,0},c_{l,0}\). The class \(\mathbb D\), the map
\(R_a\), and all scalar functionals used here are defined explicitly in
Appendix C, equations (C.0a)--(C.0i). Define \(q_0=H\omega_0\) and
\(U_0'=q_0\) with \(U_0(0)=0\). Set the sufficient focusing threshold
from {[}5, Corollary 4.6{]} \[
\underline a:=\frac{400}{848-9\pi^2}.
\tag{2.6a}
\] This threshold has a structural form. Put
\[
\mu_{\mathrm H}:=\frac{9\pi^2}{64}-\frac34,\qquad
\underline\mu:=\frac{9\pi^2}{400}-\frac3{25}.
\]
Put \(\mu=\mu(f):=2Q(f)/b(f)^2\), as in~(C.0i). Then
\(\underline a=1/(2-\underline\mu)\), and
\[
c_l=\frac{1-a(2-\mu)}{1-a\mu}>0
\quad\Longleftrightarrow\quad
a<\frac1{2-\mu}.
\]
For \(0<a<1/2\), the bound \(\mu\ge0\) gives
\((2-\mu)^{-1}\ge1/2>a\). For
\(1/2\le a<\underline a\), {[}5, Theorem 4.3{]} gives
\[
\mu\ge
\frac{\mu_{\mathrm H}}{9(1-a/3)^2}
\ge \underline\mu,
\]
so \((2-\mu)^{-1}\ge\underline a>a\). This is the two-case
focusing argument behind the threshold.

\noindent\textbf{Proposition 2.7 (full fixed-profile hypothesis registry for every Huang--Qin--Wang--Wei fixed point).} For
every \(0<a<\underline a\) and every fixed point \(f\in\mathbb D\) of
\(R_a\), the positive rescaling \[
\begin{aligned}
\alpha&=-\frac1{c_{\omega,0}},
&\qquad \beta&=\frac4\alpha,\\
\Omega(y)&=\alpha\omega_0(\beta y),
& q(y)&=\alpha q_0(\beta y),\\
U(y)&=\frac{\alpha}{\beta}U_0(\beta y),
& c_l&=\alpha c_{l,0}.
\end{aligned}
\tag{2.6b}
\] normalizes \(c_\omega=-1\) and \(\Omega'(0)=-4\). The resulting
profile satisfies the complete fixed-profile package \[
(H_{\mathrm{prof}})+(K4^+)+(D_\infty).
\tag{2.6c}
\] In particular, {[}5, Theorem 3.11{]} makes this package nonempty at
every parameter \(0<a<\underline a\).

\emph{Proof.} Appendix C.1--C.4 transfers the raw profile and similarity
constants from {[}5{]} through (2.6b), proves the weighted derivative
bounds by dyadic interpolation, derives the velocity and far-field
estimates, and verifies every registry entry. No profile selection is
used. \(\square\)

\noindent\textbf{Proposition 2.8 (uniform small-\(a\) line separation).} There
exists a qualitative threshold \(a_{\mathrm{sep}}\in(0,1/2]\) such that,
for every \(0<a<a_{\mathrm{sep}}\) and every fixed point
\(f\in\mathbb D\) of \(R_a\), \[
\mu(f)<\frac12,
\qquad
F-O=\frac{a(2\mu(f)-1)}{1-a\mu(f)}<0,
\tag{2.6d}
\] where \(F=-1+c_l/2\) and \(O=-\tilde c/2\). Thus \(F<O\): the
open inter-line band has Fredholm index \(+1\), every point in it is an
eigenvalue, and the band lies in the Browder essential spectrum. By
{[}5, Theorem 3.11{]}, actual positive-advection collapse profiles
realize this conclusion at every \(0<a<a_{\mathrm{sep}}\).

\emph{Proof.} Appendix C.5 proves the all-fixed-point estimate for
\(\mu\) and shows that one may choose \(a_{\mathrm{sep}}\le1/2\);
moreover, no all-fixed-point separation interval \((0,A)\) can have
\(A>1/2\). The spectral conclusions follow from Theorem 7.1 and are
collected in Corollary 7.3. The argument yields a qualitative, possibly
nonoptimal threshold \(a_{\mathrm{sep}}\), with the sharp upper
obstruction \(a_{\mathrm{sep}}\le1/2\). \(\square\)

\section{3. The frozen endpoint models and the inclusion of both
lines}\label{the-frozen-endpoint-models-and-the-inclusion-of-both-lines}

The two-line inclusion was previously proved in {[}6, Proposition 1{]}
under \(\mathrm{Adm}(a)\). We reproduce it under the present
fixed-profile hypotheses as the starting point for the reverse inclusion
and exactness argument.

\subsection{3.1 The far model}\label{the-far-model}

As \(y\to\infty\), the lower-order profile coefficients vanish and the
bounded transport correction is small relative to \(c_ly\partial_y\) on
far log-localized packets, so \(L_a\) approaches the pure dilation \[
A_\infty = -1 - c_l\,y\,\partial_y \qquad\text{on } L^2(dy).
\] The unitary \(u(t) = e^{t/2}\varphi(e^t)\) conjugates \(A_\infty\) to
the Fourier multiplier \((-1+c_l/2) - i c_l\tau\), so
\(A_\infty-\lambda\) is invertible iff \(\lambda\) is off the far-field line,
with the exact resolvent norm \[
\|(A_\infty-\lambda)^{-1}\|_{L^2\to L^2} = \frac{1}{d_{\mathrm{far}}(\lambda)}, \qquad d_{\mathrm{far}} = |\operatorname{Re}\lambda - (-1+c_l/2)|.
\] The kernels are one-sided powers, causal from \(0\) right of the line
and anti-causal from \(\infty\) left of it (this orientation is the sign
data for the index count), and the first-order gain is
\(\|y\partial_y(A_\infty-\lambda)^{-1}\|\le c_l^{-1} + (1+|\lambda|)/(c_l\,d_{\mathrm{far}})\).
Since the Mellin modes \(y^{-1/2+i\xi}\) are not in \(L^2\),
\(A_\infty\) has no eigenvalues and its spectrum is purely essential,
equal to the far-field line.

\subsection{3.2 The origin model on the second-derivative-shifted
line}\label{the-origin-model-on-the-second-derivative-shifted-line}

Near \(y=0\) the leading operator is the shifted dilation \[
L_a^0\varphi=(-1+q(0))\varphi-\tilde c\,y\varphi'.
\] Here the homogeneous origin space is defined intrinsically by \[
X^0:=\overline{C_c^\infty(0,\infty)}^{\,\|\partial_y^2(\cdot)\|_2}
=\{\varphi\in H^2_{\mathrm{loc}}([0,\infty)):\varphi(0)=\varphi'(0)=0,\ \varphi''\in L^2(0,\infty)\},
\] where the set on the right denotes the canonical representative and
\(\|\varphi\|_{X^0}:=\|\varphi''\|_2\). The two traces remove the affine
nullspace. The map \(\partial_y^2:X^0\to L^2(0,\infty)\) is an isometric
isomorphism: its inverse is \[
g\longmapsto\int_0^y(y-t)g(t)\,dt.
\] Indeed, an \(L^2\) function orthogonal to
\(\{\eta'':\eta\in C_c^\infty(0,\infty)\}\) has distributional second
derivative zero and hence is an \(L^2\) affine function, therefore zero.
This also proves completeness. The second-order Hardy inequality gives
\[
\|y^{-2}\varphi\|_2\le\frac43\|\varphi''\|_2. \tag{3.0}
\] The plain-\(L^2\) Mellin line of \(L_a^0\) is
\(\{\operatorname{Re}s=-\tfrac12\}\),
i.e.~\(\{\operatorname{Re}\lambda = -1+q(0)+\tilde c/2\}\); its modes
\(\varphi\sim y^{-1/2-i\tau}\) are not in \(L^2\) near \(0\) and are
excluded by the domain, which is why they are invisible in \(X\). The
marginal exponent for the \(X\)-norm is instead \(y^{3/2}\) (the
borderline power for \(\varphi''\in L^2\) near \(0\)), and at that
exponent the same dilation traces the origin line
\(\operatorname{Re}\lambda = \lambda_X = -1 + q(0) - \tfrac32\tilde c = -\tilde c/2\).
The exact conjugation identity \[
\partial_y^2\big((\tilde c\,y\partial_y + \kappa)\varphi\big) = (\tilde c\,y\partial_y + \kappa + 2\tilde c)\,\varphi''
\] turns \((L_a^0-\lambda)\varphi = f\) into
\((\tilde L-\lambda)\varphi'' = f''\) with \(\tilde L\) a dilation on
plain \(L^2\), giving the \textbf{exact} origin resolvent bound \[
\|\varphi''\|_2 \le \frac{\|f''\|_2}{d_X(\lambda)}, \qquad
|\varphi(y)|\le \frac{y^{3/2}\|\varphi''\|_{L^2(0,y)}}{\sqrt3}, \quad |\varphi'(y)|\le y^{1/2}\|\varphi''\|_{L^2(0,y)}.
\] In the parallel weighted realization on \(L^2(y^{-4}dy)\), the
unitary \(u(t) = e^{-3t/2}\varphi(e^t)\) conjugates \(L_a^0\) to
\(\lambda_X - i\tilde c\tau\), with resolvent norm \(1/d_X\). This
settles which norm carries the origin line: the weighted \(L^2\) metric
with marginal Mellin exponent \(3/2\), no anisotropic space. The
construction requires \(f'(0)=0\), which is exactly why the parametrix
must feed the origin inverse the \(X_{00}\)-type data carried by the
rank-one split of Proposition 2.2. Attempting instead to invert
\(A_\infty\) in the \(H^2\) seminorm on far supports produces a spurious
line \(\operatorname{Re}\lambda = -1-\tfrac32 c_l\) that is \textbf{not}
essential spectrum (no Weyl sequence: the tilt \(e^{4t}\) concentrates
packet mass in an \(O(1)\) log-window); the assembly therefore never
inverts the far model in the \(H^2\) seminorm, handling the far face
purely in plain \(L^2\) and recovering \(\|\varphi''\|\) a posteriori by
Theorem 2.6.

\subsection{3.3 Inclusion of both lines}\label{inclusion-of-both-lines}

Both lines lie in the essential spectrum, by explicit singular sequences
at the two endpoints.

\noindent\textbf{Graph-space transfer used below.} Give \(D_a\) the graph norm \[
\|u\|_{D_a}^2=\|u\|_X^2+\|L_au\|_X^2.
\] Suppose \(u_n\in D_a\), \(\|u_n\|_X=1\), \(u_n\rightharpoonup0\) in
\(X\), and \(\|(L_a-\lambda)u_n\|_X\to0\). Then
\(L_au_n=\lambda u_n+o_X(1)\), so \[
\|u_n\|_{D_a}\longrightarrow(1+|\lambda|^2)^{1/2}
\] and \(u_n\rightharpoonup0\) in the graph space. Indeed, the second
component of the graph inner product tends to zero because
\(L_au_n=\lambda u_n+o_X(1)\). After graph normalization, call the
sequence \(v_n\). If \(T_\lambda=L_a-\lambda\) belonged to \(\Phi_+\),
closed range and finite-dimensional kernel would give \[
\operatorname{dist}_{D_a}(v,\ker T_\lambda)
\le C\|T_\lambda v\|_X,\qquad v\in D_a.
\] The orthogonal projections of the weakly null \(v_n\) onto the
finite-dimensional kernel tend to zero, whereas the displayed estimate
would force their distance from that kernel to tend to zero. This
contradiction proves \(T_\lambda\notin\Phi_+\). Thus every weakly null
\(X\)-normalized approximate eigenvector constructed below places
\(\lambda\) in \(\sigma_{e2}(L_a)\), not merely in the approximate point
spectrum.

\emph{Far line.} Fix \(\lambda=-1+c_l/2-i c_l\xi\) and let \[
\varphi_n(y)=c_n\,\eta\!\Big(\frac{\log y-T_n}{L_n}\Big)y^{-1/2+i\xi},
\qquad L_n\to\infty,\qquad R_n:=e^{T_n-L_n}\to\infty,
\] where \(\eta\in C_c^\infty(-1,1)\) and \(c_n\) is chosen so that
\(\|\varphi_n\|_2=1\). Thus \(c_n\asymp L_n^{-1/2}\) and
\(\operatorname{supp}\varphi_n\subset[R_n,e^{2L_n}R_n]\). Direct
differentiation gives, for \(j=1,2,3\), \[
\|\partial_y^j\varphi_n\|_2\le C_jR_n^{-j},\qquad
g_n:=(A_\infty-\lambda)\varphi_n
=-\frac{c_lc_n}{L_n}\eta'\!\Big(\frac{\log y-T_n}{L_n}\Big)y^{-1/2+i\xi},
\] and hence \[
\|g_n\|_2=O(L_n^{-1}),\qquad
\|g_n''\|_2=O(L_n^{-1}R_n^{-2}). \tag{3.1a}
\] It follows that \(\|\varphi_n\|_X=1+o(1)\) and \(\|g_n\|_X=o(1)\).
After the harmless \(X\)-renormalization, the packets are
\(X\)-normalized and weakly null in \(X\).

It remains to check the full perturbation in the realization norm. Write
\[
(L_a-A_\infty)\varphi
=-aU\varphi'+q\varphi+\Omega H_{\mathrm{odd}}\varphi
-aV(\varphi)\Omega'. \tag{3.2}
\] The local terms are immediate from the derivative bounds: \[
\|(U\varphi_n')\|_X=o(1),\qquad
\|q\varphi_n\|_X=o(1). \tag{3.3}
\] Indeed \(U\) is bounded, \(U'=q\), \(U''=q'\), while
\((q\varphi_n)''=q\varphi_n''+2q'\varphi_n'+q''\varphi_n\); here
\(q''=H(\Omega'')\in L^2\), and
\(\|\varphi_n\|_\infty\lesssim L_n^{-1/2}R_n^{-1/2}\).

For the nonlocal terms use \[
(\Omega H_{\mathrm{odd}}\varphi_n)''
=\Omega''H_{\mathrm{odd}}\varphi_n
+2\Omega'H_{\mathrm{ev}}\varphi_n'
+\Omega H_{\mathrm{odd}}\varphi_n'', \tag{3.4}
\] \[
(V(\varphi_n)\Omega')''
=(H_{\mathrm{ev}}\varphi_n')\Omega'
+2(H_{\mathrm{odd}}\varphi_n)\Omega''
+V(\varphi_n)\Omega'''. \tag{3.5}
\] The \(L^2\) boundedness of \(H\) and (3.1a) make every term
containing \(\varphi_n'\) or \(\varphi_n''\) tend to zero. The remaining
terms require the fact that the Hilbert output is not supported with the
input. For \(y<R_n/2\), the separated-support kernels give \[
|H_{\mathrm{odd}}\varphi_n(y)|
\le C\int_{R_n}^\infty\frac{|\varphi_n(t)|}{t}\,dt
\le Cc_nR_n^{-1/2}
\le CL_n^{-1/2}R_n^{-1/2}. \tag{3.6}
\] On \(y\ge R_n/2\), the tails of \(\Omega,\Omega',\Omega''\)
tend uniformly to zero. Explicitly, the decay of \(\Omega\) is in
\((D_\infty)\), while for \(f=\Omega'\) and \(f=\Omega''\) the tail
estimate
\(\sup_{y\ge R}|f(y)|^2\le2\|f\|_{L^2(R,\infty)}\|f'\|_{L^2(R,\infty)}\)
follows from \(J_1,J_2,J_3\). Splitting at \(R_n/2\) in (3.4)--(3.5),
(3.6), the global \(L^2\) bounds on \(\Omega\) and \(\Omega''\), and the profile tail
bounds therefore control all terms except possibly those containing
\(V(\varphi_n)\). For those use the exact kernel \[
V(\varphi_n)(y)
=\frac1\pi\int_0^\infty\varphi_n(t)
\log\left|\frac{t-y}{t+y}\right|dt,
\] which gives
\(|V(\varphi_n)(y)|\le Cc_nyR_n^{-1/2}
\le CL_n^{-1/2}yR_n^{-1/2}\) for \(y<R_n/2\),
together with the global estimate
\(|V(\varphi_n)(y)|\le y^{1/2}\|\varphi_n\|_2\). First fix \(M\), use
the separated estimate on \(y\le M\), and then use the tails of
\(\int_M^\infty y|\Omega'|^2\) and \(\int_M^\infty y|\Omega'''|^2\).
Letting first \(n\to\infty\) and then \(M\to\infty\) yields \[
\|\Omega H_{\mathrm{odd}}\varphi_n\|_X
+\|V(\varphi_n)\Omega'\|_X\longrightarrow0. \tag{3.7}
\] Equation (3.1a) and equations (3.2)--(3.7) show
\(\|(L_a-\lambda)\varphi_n\|_X\to0\). The graph-space transfer above
therefore puts the far-field line in \(\sigma_{e2}(L_a|_X)\).

\emph{Origin line.} Fix \(\lambda=\lambda_X-i\tilde c\,\tau\), so
\(\tau=-\operatorname{Im}\lambda/\tilde c\), and set \[
\varphi_n(y)=c_n\,\eta\!\Big(\frac{\log y-T_n}{L_n}\Big)y^{3/2+i\tau},
\qquad L_n\to\infty,\qquad T_n+L_n\to-\infty,
\] with \(c_n\) chosen so that \(\|\varphi_n''\|_2=1\). Then
\(\|\varphi_n\|_2\to0\), hence \(\|\varphi_n\|_X=1+o(1)\), and the
\(X\)-renormalized packets are weakly null. The exact differentiated
dilation identity of Section 3.2 gives \[
\|(L_a^0-\lambda)\varphi_n\|_X=O(L_n^{-1}).
\] The Taylor remainders in \(b-\tilde c\,y\) and \(q-q(0)\)
give \(o(1)\) in \(X\) on the shrinking supports. The only nonlocal term
not settled by this order count is
\(\Omega H_{\mathrm{odd}}(\varphi_n'')\); Lemma 3.1 below proves its
full product contribution is \(o(1)\), and the same lemma treats the
\(V\) channel. Thus \(\varphi_n\in D_a\) and
\(\|(L_a-\lambda)\varphi_n\|_X\to0\), proving that the origin line lies
in \(\sigma_{e2}(L_a|_X)\) for the fixed profile by the same graph-space
transfer.

\subsection{3.4 Subordination of the nonlocal terms at both
endpoints}\label{subordinacy-of-the-nonlocal-terms-at-both-endpoints}

Neither frozen model contains a Hilbert term, and the reason is
structural, not a compactness accident.

\begin{itemize}
\tightlist
\item
  The \textbf{local multiplier} \(q\) has \(q(0)\ne 0\): it is
  degree-preserving and \emph{belongs} to the origin model. It freezes
  to \(q(0)\) and, together with the \(-1\) term, contributes the
  constant \(-1+q(0)\). Its non-vanishing is part of the model, not an
  obstruction.
\item
  The \textbf{nonlocal} \(\Omega H_{\mathrm{odd}}\) has coefficient
  \(\Omega\) with \(\Omega(0)=0\), \(\Omega\sim \Omega'(0) y\): on a
  homogeneous mode \(y^\mu\) it produces \(O(y^{\mu+1})\), raising the
  origin homogeneity by one, and at infinity \(\Omega\to 0\). It
  enters neither indicial family.
\item
  The \textbf{nonlocal} \(aV(\varphi)\Omega'\): \(V\) raises the power
  by one exactly (order \(-1\)), \(\Omega'\) bounded; subordinate at
  both faces by order count.
\end{itemize}

\noindent\textbf{Lemma 3.1 (nonlocal \(X\)-norm subordination on origin
packets).} For the origin packets \(\varphi_n\) of Section 3.3, \[
\|\Omega H_{\mathrm{odd}}\varphi_n\|_X
+\|V(\varphi_n)\Omega'\|_X\longrightarrow0. \tag{3.8}
\]

\emph{Proof.} Write \[
\alpha=\frac32+i\tau,\qquad
\eta_n(y)=\eta\!\left(\frac{\log y-T_n}{L_n}\right),\qquad
\varphi_n(y)=c_ny^\alpha\eta_n(y),
\] and let \([u_n,s_n]\) be its support. Then \(s_n\to0\) and direct
differentiation gives \[
\varphi_n''(y)=c_ny^{\alpha-2}
\left[
\alpha(\alpha-1)\eta_n(y)
+\frac{2\alpha-1}{L_n}\eta'\!\left(\frac{\log y-T_n}{L_n}\right)
+\frac1{L_n^2}\eta''\!\left(\frac{\log y-T_n}{L_n}\right)
\right]. \tag{3.9}
\] The normalization gives \(c_n\asymp L_n^{-1/2}\) and
\(\|\varphi_n''\|_2=1\). Since \(\varphi_n\) and \(\varphi_n'\) vanish
to the left of \(u_n\), two integrations and Cauchy--Schwarz give \[
\|\varphi_n'\|_2\le s_n,\qquad
\|\varphi_n\|_2\le s_n^2. \tag{3.10}
\] Equation (3.9) also gives the moment bound \[
M_n:=\int_{u_n}^{s_n}t|\varphi_n''(t)|\,dt
\le Cc_ns_n^{3/2}. \tag{3.11}
\]

The only term not covered immediately by \(L^2\) boundedness of the
folded Hilbert transforms is \(\Omega H_{\mathrm{odd}}(\varphi_n'')\).
For \(y\ge2s_n\) the odd kernel has no principal value singularity and
\[
|H_{\mathrm{odd}}(\varphi_n'')(y)|
\le\frac{C}{y^2}\int_{u_n}^{s_n}t|\varphi_n''(t)|\,dt
\le\frac{CM_n}{y^2}. \tag{3.12}
\] On \((0,2s_n)\), the origin bound \(|\Omega(y)|\le Cy\) and the
\(L^2\) boundedness of \(H_{\mathrm{odd}}\) give \[
\|\Omega H_{\mathrm{odd}}(\varphi_n'')\|_{L^2(0,2s_n)}
\le Cs_n.
\] For \(n\) sufficiently large that \(2s_n<r_{\mathrm{org}}\), on
\((2s_n,r_{\mathrm{org}})\)
combine \(|\Omega(y)|\le Cy\) with (3.12); the resulting norm is at
most \(CM_ns_n^{-1/2}\le Cc_ns_n\). On \((r_{\mathrm{org}},\infty)\), boundedness of
\(\Omega\) and (3.12) give the bound \(CM_n\le Cc_ns_n^{3/2}\). Hence
\[
\|\Omega H_{\mathrm{odd}}(\varphi_n'')\|_2\longrightarrow0. \tag{3.13}
\]

Now differentiate the first nonlocal product twice: \[
(\Omega H_{\mathrm{odd}}\varphi_n)''
=\Omega''H_{\mathrm{odd}}\varphi_n
+2\Omega'H_{\mathrm{ev}}\varphi_n'
+\Omega H_{\mathrm{odd}}\varphi_n''.
\] The first two terms tend to zero by (3.10), the global coefficient
bounds, and \(L^2\) boundedness of \(H\); the last one is (3.13). For
the Volterra channel, \[
(V(\varphi_n)\Omega')''
=(H_{\mathrm{ev}}\varphi_n')\Omega'
+2(H_{\mathrm{odd}}\varphi_n)\Omega''
+V(\varphi_n)\Omega'''.
\] The first two terms again tend to zero. The pointwise estimate
\(|V(\varphi_n)(y)|\le y^{1/2}\|\varphi_n\|_2\) and the weighted bound
\(J_3<\infty\) give \[
\|V(\varphi_n)\Omega'''\|_2
\le J_3\|\varphi_n\|_2\longrightarrow0. \tag{3.14}
\] Finally,
\(\|\Omega H_{\mathrm{odd}}\varphi_n\|_2\le W_0\|\varphi_n\|_2\) and
\(\|V(\varphi_n)\Omega'\|_2\le J_1\|\varphi_n\|_2\). Together with
(3.10)--(3.14), these estimates prove (3.8). \(\square\)

The essential point is that \(\Omega H_{\mathrm{odd}}\) is
\textbf{noncompact}, even after \(y=e^t\): there it is a vanishing
coefficient times a Mellin multiplier that does not vanish at high
frequency. Endpoint subordination supplies the needed replacement. At
infinity, translated coefficient tails tend strongly to zero, and the
far model is used only in plain \(L^2\). At the origin,
\(|\Omega(y)|\lesssim y\), together with the separated-support kernel
gain in (3.11)--(3.13), supplies the extra power of \(y\) required by the
\(X\)-estimate. For the reverse estimate, Lemmas 4.1--4.2 first remove
the rank-one trace channel before weighted inversion. Sections 3.1--3.3
turn this endpoint information into the line inclusion. Section 4 gives
the reverse inclusion at the
level of \(\sigma_{e2}\), and the common-domain homotopy and local
anchor of Sections 4.4 and 5.2 upgrade the result to Fredholmness and
\(\sigma_{e3}\).

\section{4. The semi-Fredholm
estimate}\label{the-semi-fredholm-estimate}

This is the technical heart. We prove a three-zone a priori estimate
that gives, for every \(\lambda\) off the two lines, finite kernel and
closed range for \(L_a-\lambda\), uniformly on compact \(\lambda\)-sets
for the fixed profile.

\subsection{4.1 Setup and the scalar
cancellation}\label{setup-and-the-scalar-cancellation}

Fix a compact \(K\subset\mathbb{C}\) with
\(\mathrm{dist}(K,\text{far-field line})\ge d\) and
\(\mathrm{dist}(K,\text{origin line})\ge d\). Fix
\(\delta_0\in(0,\min(r_{\mathrm{org}},1)/4]\), rescale \(\chi_0\) to \(\equiv 1\) on
\([0,\delta_0/2]\), supported in \([0,\delta_0]\), and write
\(\varphi = \psi + l\,y\chi_0\), \(l=\varphi'(0)\), so
\(\psi(0)=\psi'(0)=0\) and \(\psi=\varphi\) on \(y\ge\delta_0\). The
equation \((L_a-\lambda)\varphi=f\) becomes \[
(L_a-\lambda)\psi = f - l\,g_\star, \qquad g_\star := (L_a-\lambda)(y\chi_0),
\] with \(g_\star\) smooth and odd near \(0\) (hence
\(g_\star(y)-g_\star'(0)y=O(y^3)\)),
\(g_\star'(0) = -\lambda + (1-a)\Omega'(0) h_1\)
(\(h_1 = H_{\mathrm{odd}}(y\chi_0)(0)\); the constant
\(-1+q(0)-\tilde c\) vanished by the exact identity), and
\(C_g(\lambda):=\|w(g_\star-g_\star'(0)y\chi_0)\|_2<\infty\), bounded
above by an affine function of \(|\lambda|\) and therefore uniformly
bounded for \(\lambda\in K\). Measuring data in
\(\|f\|_\# := \|w f_w\|_2 + |f'(0)|\), \(f_w = f - f'(0)y\chi_0\), one
has \(\|f\|_\#\le C_\#\|f\|_X\).

The one place the estimate could silently be infinite is the
\(y\)-direction: \(wy = y^{-1}\) is not locally \(L^2\), so any linear
part in the data would make \(\|w(\cdot)\|_2\) diverge. All three
\(y\)-direction linear parts (\(f'(0)\), \(l g_\star'(0)\), and the
nonlocal output \(\Omega'(0)(1-a)(H_{\mathrm{odd}}\psi)(0)\)) must be routed
through a single identity before any weighted inversion.

\noindent\textbf{Lemma 4.1 (scalar cancellation).} For \(\varphi\in D_a\),
\((L_a-\lambda)\varphi=f\), \(\psi=\varphi-l\,y\chi_0\), \[
f'(0) - l\,g_\star'(0) = (1-a)\,\Omega'(0)\,(H_{\mathrm{odd}}\psi)(0).
\] (Lemma 2.5 applied to \(\psi\), whose \(p\)-value is \(0\).) The
three linear parts cancel inside the \(\mathrm{span}\{y\chi_0\}\)
component of the triangular split \(X = \mathbb{C}e\oplus X_w\) and never
enter the weighted component. The off-diagonal entries are the three
bounded scalar functionals \(l(\varphi)\),
\((H_{\mathrm{odd}}\psi)(0)\), \(f'(0)\), each rank-one hence compact;
they sit in the compact remainder.

\subsection{4.2 The origin nonlocal
machinery}\label{the-origin-nonlocal-machinery}

The odd kernel gives
\(H_{\mathrm{odd}}\psi(0) = -\tfrac2\pi\int_0^\infty\psi(t)/t\,dt\), a
bounded rank-one functional,
\(|H_{\mathrm{odd}}\psi(0)|\le C N(\varphi)\). The absorption estimate,
with the linear part removed, is:

\noindent\textbf{Lemma 4.2 (weighted origin bounds for both nonlocal channels).}
Put \[
h:=H_{\mathrm{odd}}\psi,\qquad m:=h(0),
\] and let \(0<\delta\le\min(\delta_0,r_{\mathrm{org}}/4,1/4)\). Then \[
\|\mathbf1_{(0,\delta)}w[\Omega h-\Omega'(0)my]\|_2
\le\varepsilon_H(\delta)N(\varphi), \tag{4.0a}
\] \[
\|\mathbf1_{(0,\delta)}w[V(\psi)\Omega'-\Omega'(0)my]\|_2
\le\varepsilon_{V,0}(\delta)N(\varphi), \tag{4.0b}
\] where \[
\varepsilon_H(\delta)=\delta\left[M_3+|\Omega'(0)|\left(\frac8{c_{\min}}+32+C_{H,f}\right)\right],
\qquad C_{H,f}:=\frac{64}{9\pi}\left(2^{-1/2}+5^{-1/2}\right),
\] and \(\varepsilon_{V,0}(\delta)=C_{V,0}\delta\), with \[
C_{V,0}:=W_1\left[\frac43\left(\frac4{c_{\min}}+16\right)+\frac{64}{9\pi}\right]
+\frac{C_mC_{\Omega',0}}{\sqrt3},
\] \[
C_m:=\frac2\pi(1+3^{-1/2}),\qquad
C_{\Omega',0}:=\frac{M_3(1+Q_0)+|\Omega'(0)|M_4(\tfrac12+\tfrac a6)}{c_{\min}}.
\] The two far-piece constants are deliberately estimated separately:
the Hilbert bound retains both subinterval contributions in (4.0f),
whereas the extra Volterra integration and \(\delta\le1/4\) permit the
smaller coefficient in (4.0k). Only their common \(O(\delta)\) decay is
used below. In particular, \[
\|\mathbf1_{(0,\delta)}w[P_a\psi-(1-a)\Omega'(0)my]\|_2
\le[\varepsilon_H(\delta)+a\varepsilon_{V,0}(\delta)]N(\varphi). \tag{4.0c}
\]

\emph{Proof.} Write \[
\Omega h-\Omega'(0)my=(\Omega-\Omega'(0)y)h+\Omega'(0)y(h-m).
\] The first term is at most \(M_3\delta N(\varphi)\) after
multiplication by \(w\), because \(|\Omega-\Omega'(0)y|\le M_3y^3\),
\(\|h\|_2\le\|\psi\|_2\), and
\(\|\psi\|_2\le\|w\psi\|_2\le N(\varphi)\), since \(w\ge1\). Thus no
cutoff-dependent constant enters this modulus. Choose a
smooth cutoff \(\zeta=1\) on \([0,2\delta]\), supported in
\([0,4\delta]\), with \(|\zeta'|\le\delta^{-1}\), and split
\(\psi=\psi_n+\psi_f\), where \(\psi_n=\zeta\psi\). For
\(h_n=H_{\mathrm{odd}}\psi_n\) and \(m_n=h_n(0)\), the first-order Hardy
inequality gives \[
\left\|\frac{h_n-m_n}{y}\right\|_2
\le2\|h_n'\|_2\le2\|\psi_n'\|_2
\le\delta\left(\frac8{c_{\min}}+32\right)N(\varphi). \tag{4.0d}
\] Here we used \(b/y\ge c_{\min}\) on \([0,4\delta]\) and \[
\|\psi_n'\|_2\le\frac{4\delta}{c_{\min}}\|wb\psi'\|_2+16\delta\|w\psi\|_2. \tag{4.0e}
\] For the far part, \(\psi_f=0\) on
\([0,2\delta]\), and the folded kernel obeys \[
\left|\partial_y\frac{2t}{y^2-t^2}\right|\le\frac{64}{9}\delta t^{-3}
\qquad(0<y\le\delta,\ t\ge2\delta).
\] Consequently \[
\sup_{0<y<\delta}|h_f'(y)|
\le\frac{64}{9\pi}\delta\int_{2\delta}^\infty|\psi(t)|t^{-3}\,dt,
\] while Cauchy--Schwarz gives \[
\int_{2\delta}^\infty|\psi(t)|t^{-3}\,dt
\le[(2\delta)^{-1/2}+5^{-1/2}]\|w\psi\|_2. \tag{4.0f}
\] Since \(h_f(y)-m_f=\int_0^yh_f'(s)\,ds\), (4.0f) yields \[
\left\|\frac{h_f-m_f}{y}\right\|_{L^2(0,\delta)}
\le C_{H,f}\delta N(\varphi).
\] This proves (4.0a).

For the Volterra channel, first note that \[
|m|\le\frac2\pi\left(\|t\|_{L^2(0,1)}+\|t^{-1}\|_{L^2(1,\infty)}\right)\|w\psi\|_2
=C_m\|w\psi\|_2. \tag{4.0g}
\] The quadratic derivative remainder follows from the profile equation
and introduces no new hypothesis. Evenness gives \(q'(0)=0\), and hence
\[
|q(y)-q(0)|\le\tfrac12M_4y^2,\qquad
|b(y)-\tilde c y|\le\tfrac a6M_4y^3.
\] Using \(q(0)-1=\tilde c\) and \(b\Omega'=\Omega(q-1)\), \[
b(\Omega'-\Omega'(0))=(\Omega-\Omega'(0)y)(q-1)+\Omega'(0)y(q-q(0))-\Omega'(0)(b-\tilde c y).
\] Since \(b\ge c_{\min}y\), \[
|\Omega'(y)-\Omega'(0)|\le C_{\Omega',0}y^2. \tag{4.0h}
\] Now write \[
V(\psi)\Omega'-\Omega'(0)my=\Omega'[V(\psi)-my]+my(\Omega'-\Omega'(0)). \tag{4.0i}
\] After multiplication by \(w\), the last term in (4.0i) is bounded by
\((C_mC_{\Omega',0}/\sqrt3)\delta^{3/2}N(\varphi)\), hence by the same
constant times \(\delta N(\varphi)\).

For \(j=n,f\), put \(h_j=H_{\mathrm{odd}}\psi_j\), \(m_j=h_j(0)\), and
\[
G_j(y):=\int_0^y(h_j(s)-m_j)\,ds.
\] Then \(V(\psi)-my=G_n+G_f\). For the near part, \(G_n(0)=G_n'(0)=0\)
and \(G_n''=H_{\mathrm{ev}}(\psi_n')\). The second-order Hardy
inequality and (4.0e) give \[
\|y^{-2}G_n\|_2\le\frac43\|G_n''\|_2
\le\frac43\|\psi_n'\|_2
\le\frac43\delta\left(\frac4{c_{\min}}+16\right)N(\varphi). \tag{4.0j}
\] For the far part, \[
G_f(y)=\int_0^y(y-s)h_f'(s)\,ds.
\] The kernel bound and (4.0f) therefore imply \[
\|y^{-2}G_f\|_{L^2(0,\delta)}
\le\frac{64}{9\pi}\delta N(\varphi). \tag{4.0k}
\] Combining (4.0h)--(4.0k) with \(|\Omega'|\le W_1\) proves (4.0b), and
(4.0c) follows by subtraction. \(\square\)

Thus on \((0,\delta)\) the two nonlocal terms equal \(\Omega'(0)(1-a)my\) (the
rank-one output routed through Lemma 4.1) plus a weighted-\(O(\delta)\)
error. This proves the required endpoint subordination.

The freezing errors at both faces are elementary from the registry: at
the origin \(|(b-\tilde c y)/y|\le\tfrac{aM_1}{2}y\),
\(\sup|q-q(0)|\le 2M_1\delta\), and the nonlocal terms on
\((0,2\delta]\) are \(O(\delta)\|\varphi\|_2\) in the unweighted
\(L^2\) norm; at infinity
\(|(b-c_l y)/y|\le aI_q/R = O(1/R)\) (the coefficient error \(a|U|\)
does not vanish, only the relative error does, so the far face is always
run against the Euler model in log variables),
\(\sup_{y\ge R}|q|\le 2C_q/R\), and the nonlocal terms are
small-coefficient times the global bounded \(H\), with no kernel
splitting needed. Transition-cutoff terms are supported in compact
annuli and enter the Rellich remainder; the restriction
\(X\to L^2(\delta,R)\) is compact.

\subsection{4.3 The three-zone estimate}\label{the-three-zone-estimate}

Take cutoffs \(\eta_0\equiv 1\) on \((0,1/R]\) supported \((0,2/R]\),
\(\eta_\infty\equiv 1\) on \([R,\infty)\) supported \([R/2,\infty)\),
and require \(4/R\le\delta_0\) so the origin zone sits inside
\(\{\chi_0\equiv 1\}\).

\noindent\textbf{Step A (weighted zeroth-order).} We write out the origin
inversion, including its \(X_{00}\) compatibility. Set \[
e=y\chi_0,\qquad m(\psi):=(H_{\mathrm{odd}}\psi)(0),\qquad
c(\psi):=(1-a)\Omega'(0)m(\psi),
\] and \[
P_a\psi:=\Omega H_{\mathrm{odd}}\psi-aV(\psi)\Omega',\qquad
F_w:=f-lg_\star-c(\psi)e,\qquad
P_{a,w}\psi:=P_a\psi-c(\psi)e.
\] Lemma 4.1 says \(c(\psi)=f'(0)-lg_\star'(0)\). Therefore \[
F_w=[f-f'(0)e]-l[g_\star-g_\star'(0)e],\qquad
F_w'(0)=(P_{a,w}\psi)'(0)=0,
\] and, before any weighted inversion, the equation is exactly \[
(L_{\mathrm{loc}}-\lambda)\psi+P_{a,w}\psi=F_w. \tag{4.1a}
\] Thus both data terms presented to the origin model lie in \(X_{00}\).

Let \[
M_0:=-1+q(0)-\tilde c\,y\partial_y.
\] On the support of \(\eta_0\) one has \(e=y\), and direct localization
of (4.1a) gives \[
\begin{aligned}
(M_0-\lambda)(\eta_0\psi)
={}&\eta_0F_w-\eta_0P_{a,w}\psi
+\eta_0(b-\tilde c y)\psi'\\
&+\eta_0(q(0)-q)\psi-\tilde c\,y\eta_0'\psi .
\end{aligned} \tag{4.2a}
\] The left side and every term on the right have zero first trace. The
exact Euler inverse on \(L^2(y^{-4}dy)\) therefore yields \[
\|w\eta_0\psi\|_2
\le d_X(\lambda)^{-1}
\left\|w\left[\text{right side of (4.2a)}\right]\right\|_2. \tag{4.3a}
\] The data term is at most \(\|f\|_\#+|l|C_g(\lambda)\). Let
\(\varepsilon_P(\delta)=\varepsilon_H(\delta)+a\varepsilon_{V,0}(\delta)\)
denote the sum of the two \(O(\delta)\) moduli in Lemma 4.2 for the
Hilbert and Volterra channels. Then \[
\|w\eta_0P_{a,w}\psi\|_2\le\varepsilon_P(2/R)N(\varphi).
\] The claimed \(O(R^{-2})\) freezing uses evenness and \(M_4\), not
merely a Lipschitz bound: \[
|q(y)-q(0)|\le\tfrac12M_4y^2,\qquad
|b(y)-\tilde c y|\le\tfrac16aM_4y^3.
\] Since \(y\le2/R\) and \(b/y\ge c_{\min}\) there, \[
\|w\eta_0(q(0)-q)\psi\|_2
\le\frac{2M_4}{R^2}\|w\psi\|_2,
\] \[
\|w\eta_0(b-\tilde c y)\psi'\|_2
\le\frac{2aM_4}{3c_{\min}R^2}\|wb\psi'\|_2. \tag{4.4a}
\] Set
\(\mathrm{REM}:=\|\mathbf{1}_{[1/R,2R]}\varphi\|_2+|l|+|m(\psi)|\).
The interval \([1/R,2R]\) contains both transition annuli
\([1/R,2/R]\) and \([R,2R]\). The origin commutator is supported in
\([1/R,2/R]\) and enters
\(C(R)\mathrm{REM}\).

For completeness, with \(M_\infty=-1-c_ly\partial_y\), the far
localization obeys \[
(M_\infty-\lambda)(\eta_\infty\psi)
=\eta_\infty(f-lg_\star-P_a\psi)
+\eta_\infty(b-c_ly)\psi'-\eta_\infty q\psi-c_ly\eta_\infty'\psi.
\] Its exact \(L^2\) Euler inverse has norm at most the reciprocal
distance to the far-field line. The registry gives respectively \(O(R^{-1})\),
\(O(R^{-1})\), and \(O(R^{-1})+\varepsilon_V(R/2)\) for the two freezing
errors and the nonlocal channels. More explicitly, with constants
depending only on the fixed registry, \[
\|\eta_\infty(b-c_ly)\psi'\|_2\le \frac{C}{R}\|b\psi'\|_{L^2(R/2,\infty)},
\] \[
\|\eta_\infty q\psi\|_2\le \varepsilon_q(R/2)N(\varphi),\qquad
\|\eta_\infty P_a\psi\|_2
\le[\varepsilon_\Omega(R/2)+a\varepsilon_V(R/2)]N(\varphi).
\] The cutoff commutator again enters the remainder. Combining the two
inversions gives \[
\|w\psi\|_2 \le \frac{2}{d}\|f\|_\# + \varepsilon_A(R)\,N(\varphi) + C(R)\,\mathrm{REM},
\qquad \varepsilon_A(R)\to 0, \tag{4.5}
\] The interior \([1/R,R]\) carries no inversion and lands in
\(\mathrm{REM}\) by Rellich.

\noindent\textbf{Step B (graph reconstruction, \(R\)-independent).} The
trace-free regrouping (4.1a) gives the checkable identity \[
b\psi'=-(1+\lambda-q)\psi-F_w+P_{a,w}\psi. \tag{4.6}
\] Equivalently, before the rank-one terms are combined, \[
b\psi'=-(f-lg_\star)-(1+\lambda-q)\psi
+\Omega H_{\mathrm{odd}}\psi-aV(\psi)\Omega'.
\] Multiply (4.6) by \(w\) and take \(L^2\). The local terms are bounded
by \(\|f\|_\#+|l|C_g(\lambda)
+(1+|\lambda|+Q_0)\|w\psi\|_2\). On \((0,\delta_0/2)\) one has \(e=y\),
so Lemma 4.2 applies directly to
\(P_{a,w}\psi=P_a\psi-(1-a)\Omega'(0)m(\psi)y\).

It remains to record the outer estimate without introducing the
nonintegrable function \(y-e\). Put \(E=[\delta_0/2,\infty)\). Since
\(w\) is bounded on \(E\), the global \(L^2\) bounds give \[
\begin{aligned}
\|\mathbf1_EwP_{a,w}\psi\|_2
&\le C(\delta_0)
\left(\|P_a\psi\|_2+|(1-a)\Omega'(0)m(\psi)|\,\|e\|_2\right)\\
&\le C(\delta_0)
\left[(W_0+aJ_1)\|w\psi\|_2+|m(\psi)|\right].
\end{aligned} \tag{4.6a}
\] Thus the outer nonlocal contribution enters the existing
\(\|w\psi\|_2\) and \(\mathrm{REM}\) terms. As in Step A, the three
\(y\)-direction rank-one parts, namely \(f'(0)\), \(l\,g_\star'(0)\),
and the nonlocal output \(\Omega'(0)(1-a)(H_{\mathrm{odd}}\psi)(0)\), cancel
through Lemma 4.1 before weighted inversion. Hence \[
\|w\,b\psi'\|_2 \le \|f\|_\# + C_1\|w\psi\|_2 + \varepsilon_P(\delta_0)\,N(\varphi) + C\,\mathrm{REM},
\tag{4.7}
\] with \(C_1 = C_1(\delta_0,K,\text{registry})\) independent of \(R\).

\noindent\textbf{Step C (closure).} Finalize \(\delta_0\) with
\(\varepsilon_P(\delta_0)\le\tfrac14\), then \(R=R_0\ge 4/\delta_0\)
with \((2+C_1)\varepsilon_A(R_0)\le\tfrac14\). Adding Steps A and B,
absorbing the \(N\)-terms (total coefficient \(\le\tfrac12\)) and
\(|l|\le\mathrm{REM}\):

\noindent\textbf{Theorem 4.3 (semi-Fredholm estimate, fixed profile).} For
\(\lambda\in K\) there are
\(C_*=C_*(K,d,\text{fixed-profile registry})\) and \(R_0\) with, for all
\(\varphi\in D_a\), \[
N(\varphi) \le C_*\big(\|(L_a-\lambda)\varphi\|_\# + \|\mathbf{1}_{[1/R_0,2R_0]}\varphi\|_2 + |\varphi'(0)| + |(H_{\mathrm{odd}}\psi)(0)|\big).
\] With the interface bounds of Section 2.4 and the bootstrap Theorem
2.6, \[
\|\varphi\|_X \le C\big(\|(L_a-\lambda)\varphi\|_X + \|\mathcal{K}\varphi\|_Y\big),
\] where \[
Y=L^2([1/R_0,2R_0])\oplus\mathbb C^2,
\] \[
\mathcal K\varphi=
\left(
\mathbf1_{[1/R_0,2R_0]}\varphi,\,
\varphi'(0),\,
(H_{\mathrm{odd}}(\varphi-\varphi'(0)e))(0)
\right). \tag{4.7a}
\] This map is compact from \(D_a\) with its graph norm to \(Y\).
Indeed, graph-bounded sets are bounded in \(H^2\) on the displayed
compact interval by Theorem 2.6, so the first component is compact by
Rellich; the other two components are bounded scalar functionals and
hence finite-rank. Therefore, for the fixed profile, \(L_a-\lambda\) has
finite-dimensional kernel and closed range for every \(\lambda\) off the
two lines, with constants uniform on \(K\).

Combined with the inclusion of Section 3.3, this gives
\(\sigma_{e2}(L_a|_X)\) equal to the two-line union for the fixed
profile.

\subsection{4.4 From semi-Fredholm to
Fredholm}\label{from-semi-fredholm-to-fredholm}

\noindent\textbf{Corollary 4.4 (common-domain upper-semi-Fredholm homotopy).} Set
\[
P_a\varphi:=\Omega H_{\mathrm{odd}}\varphi-aV(\varphi)\Omega',\qquad L_a^t:=L_{\mathrm{loc}}+tP_a,\qquad 0\le t\le1.
\] Under \((K4^+)\), \(P_a\in\mathcal B(X)\). Indeed, \[
(\Omega H_{\mathrm{odd}}\varphi)''=\Omega''H_{\mathrm{odd}}\varphi+2\Omega'H_{\mathrm{ev}}\varphi'+\Omega H_{\mathrm{odd}}\varphi'',
\] \[
(V(\varphi)\Omega')''=(H_{\mathrm{ev}}\varphi')\Omega'+2(H_{\mathrm{odd}}\varphi)\Omega''+V(\varphi)\Omega''',
\] so interpolation \(\|\varphi'\|_2\le2^{-1/2}\|\varphi\|_X\) and
\(|V(\varphi)(y)|\le y^{1/2}\|\varphi\|_2\) give \[
\|P_a\varphi\|_X\le C_P\|\varphi\|_X,\qquad
C_P:=2W_0+\sqrt2W_1+W_2
+a(J_1+2^{-1/2}W_1+2W_2+J_3).
\] Consequently the maximal realizations have the common domain
\(D:=D(L_{\mathrm{loc}})=D(L_a^t)\). The bounded perturbation of the
closed operator \(L_a\) shows that \(L_{\mathrm{loc}}\) is closed, so
\(D\) with its \(L_{\mathrm{loc}}\) graph norm is a Hilbert space, and
\[
\|(L_a^t-L_a^s)\varphi\|_X\le |t-s|C_P\|\varphi\|_X\le |t-s|C_P\|\varphi\|_D. \tag{4.8}
\] Thus \(t\mapsto L_a^t-\lambda\) is norm-continuous in
\(\mathcal B(D,X)\). We verify every \(t\)-dependent input to the
three-zone estimate. For \(g_{\star,t}:=(L_a^t-\lambda)e\) and
\(\varphi=le+\psi\), \[
g_{\star,t}'(0)=-\lambda+t(1-a)\Omega'(0)(H_{\mathrm{odd}}e)(0),\qquad
f'(0)-lg_{\star,t}'(0)=t(1-a)\Omega'(0)(H_{\mathrm{odd}}\psi)(0). \tag{4.9}
\] This is Lemmas 2.5 and 4.1 with the common factor \(t\). If
\(c_t(\psi)=t(1-a)\Omega'(0)(H_{\mathrm{odd}}\psi)(0)\), then \[
F_{w,t}=f-lg_{\star,t}-c_t(\psi)e,\qquad
P_{w,t}\psi=tP_a\psi-c_t(\psi)e
\] both have zero first trace and satisfy
\((L_{\mathrm{loc}}-\lambda)\psi+P_{w,t}\psi=F_{w,t}\). Consequently the
exact \(X_{00}\) origin inversion (4.2a)--(4.4a) is unchanged, while the
two moduli of Lemma 4.2 are multiplied by \(t\le1\). The far nonlocal
errors and every Step-B reconstruction term likewise acquire only this
harmless factor. Hence Steps A--C use the same choices of
\(\delta_0,R_0\) for all \(t\in[0,1]\).

The bootstrap is uniform for the same explicit reason. Its
differentiated right side becomes \[
\mathfrak r_t=q'\varphi+t(1-a)\Omega'h+t\Omega H_{\mathrm{ev}}\varphi'
-taV(\varphi)\Omega''-f',
\qquad \mathfrak r_t(0)=\lambda\varphi'(0).
\] The Euler denominator and both indicial lines are independent of
\(t\), and all new coefficients are bounded by their \(t=1\) values.
Thus, on every compact \(K\) separated from the origin line, \[
\|\varphi''\|_2+\|y\varphi'\|_2
\le C_K\bigl(\|(L_a^t-\lambda)\varphi\|_X+\|\varphi\|_2\bigr),
\qquad 0\le t\le1, \tag{4.10}
\] with one constant \(C_K\). Combining (2.4), (4.9), (4.10), and the
uniform versions of Lemmas 4.1--4.2, the proof of Theorem 4.3 applies
uniformly to \(tP_a\). Hence \[
L_a^t-\lambda\in\Phi_+(D,X)\qquad(0\le t\le1)
\] for every \(\lambda\) off the two indicial lines. Here \(\Phi_+\)
denotes the upper-semi-Fredholm operators (closed range and
finite-dimensional kernel). \(\square\)

Theorem 4.3 and Corollary 4.4 give a norm-continuous path entirely
inside \(\Phi_+\). Section 5 uses variation of constants, compactly
supported data, and the closed-range estimate to prove that the local
anchor \(L_{\mathrm{loc}}-\lambda\) is Fredholm and to compute its
index. The continuation lemma then propagates finite cokernel and the
index along the whole \(\Phi_+\) path to \(L_a-\lambda\).

\section{5. Index bookkeeping and the band
classification}\label{index-bookkeeping-and-the-band-classification}

\subsection{5.1 Fixed-profile real-part
localization}\label{fixed-profile-real-part-localization}

\noindent\textbf{Lemma 5.1 (Hilbert-cancelled fixed-profile real-part
enclosure).} Put \[
q_-:=\inf_{y>0}q(y),\qquad q_+:=\sup_{y>0}q(y),
\] and \[
E_-:=-1+\frac{c_l}{2}+\left(1+\frac a2\right)q_- -aJ_1,
\qquad
E_+:=-1+\frac{c_l}{2}+\left(1+\frac a2\right)q_+ +aJ_1.
\] Then every \(\varphi\in D_a\) and every \(\lambda\) with
\(\operatorname{Re}\lambda\notin[E_-,E_+]\) satisfy \[
\|(L_a-\lambda)\varphi\|_2
\ge \operatorname{dist}\!\bigl(\operatorname{Re}\lambda,[E_-,E_+]\bigr)
\|\varphi\|_2. \tag{5.0}
\] In particular \(L_a-\lambda\) is injective there, and every
eigenvalue of the maximal \(X\)-realization lies in the closed vertical
slab \(E_-\le\operatorname{Re}\lambda\le E_+\).

\emph{Proof.} Write \(h=H_{\mathrm{odd}}\varphi\). We first record the
half-line Tricomi identity \[
\operatorname{Re}\langle q\varphi+\Omega h,\varphi\rangle
=\frac12\int_0^\infty q\bigl(|\varphi|^2+|h|^2\bigr). \tag{5.0a}
\] For a real smooth compactly supported odd function \(f\) on the full
line, the Tricomi identity and skew-adjointness of \(H\) give \[
H(fHf)=\frac12\bigl[(Hf)^2-f^2\bigr],\qquad H^*=-H,
\] and hence \[
\int_{\mathbb R}\Omega fHf
=\frac12\int_{\mathbb R}(H\Omega)\bigl[(Hf)^2-f^2\bigr].
\] Adding \(\int_{\mathbb R}qf^2\), with \(q=H\Omega\), and
restricting the resulting even integrands to \((0,\infty)\) proves
(5.0a) for real \(f\). Both sides are continuous quadratic forms on
\(L^2\) because \(q,\Omega\in L^\infty\) and \(H\) is an \(L^2\)
isometry, so density extends the identity to real \(L^2\) data. Applying
it separately to the real and imaginary parts proves the complex case.
The same full-line isometry gives \(\|h\|_2=\|\varphi\|_2\) on the
half-line.

For completeness, the transport integration by parts has no hidden
endpoint term. With \(y=e^t\) and \(u(t):=e^{t/2}\varphi(e^t)\), the
facts \(\varphi,y\varphi'\in L^2(0,\infty)\) give
\(u\in H^1(\mathbb R)\) and hence \(u(t)\to0\) at both ends. Since
\(b(y)/y\) is bounded,
\(b(y)|\varphi(y)|^2=(b(y)/y)|u(\log y)|^2\to0\) as \(y\downarrow0\)
and \(y\to\infty\). Therefore integration by parts, together with
\(b'=c_l+a q\), yields \[
\begin{aligned}
\operatorname{Re}\langle(L_a-\lambda)\varphi,\varphi\rangle
={}&\left(-1+\frac{c_l}{2}-\operatorname{Re}\lambda\right)
\|\varphi\|_2^2\\
&+\frac{1+a}{2}\int_0^\infty q|\varphi|^2
+\frac12\int_0^\infty q|h|^2
-a\operatorname{Re}\langle V(\varphi)\Omega',\varphi\rangle .
\end{aligned} \tag{5.0b}
\] Finally, \[
|\langle V(\varphi)\Omega',\varphi\rangle|
\le\|\varphi\|_2\int_0^\infty y^{1/2}|\Omega'||\varphi|
\le J_1\|\varphi\|_2^2.
\] Because \((1+a)/2\) and \(1/2\) are positive, the two \(q\) terms
in (5.0b) lie between \((1+a/2)q_-\|\varphi\|_2^2\) and
\((1+a/2)q_+\|\varphi\|_2^2\). Thus the real part is at most
\(-(\operatorname{Re}\lambda-E_+)\|\varphi\|_2^2\) to the right of the
slab and at least \((E_--\operatorname{Re}\lambda)\|\varphi\|_2^2\) to
its left. Cauchy--Schwarz gives (5.0). \(\square\)

\emph{Remark (the optimal unweighted \(L^2\) energy slab).} The argument
uses the bounded composition \[
A_a:=M_{\Omega'}V\in\mathcal B\bigl(L^2(0,\infty)\bigr),
\qquad \|A_a\|\le J_1;
\] and then sets \[
B_a^{\mathrm H}:=
\frac{1+a}{2}M_{q}
+\frac12H_{\mathrm{odd}}^*M_{q}H_{\mathrm{odd}}
-\frac a2\left(A_a+A_a^*\right).
\] Equation (5.0b) says exactly that \[
\operatorname{Re}\langle L_a\varphi,\varphi\rangle
=\left(-1+\frac{c_l}{2}\right)\|\varphi\|_2^2
+\langle B_a^{\mathrm H}\varphi,\varphi\rangle.
\] Thus the optimal endpoints for this unweighted energy identity are
\(-1+c_l/2+\inf\sigma(B_a^{\mathrm H})\) and
\(-1+c_l/2+\sup\sigma(B_a^{\mathrm H})\). The explicit values \(E_\pm\)
follow from \(\|H_{\mathrm{odd}}\varphi\|_2=\|\varphi\|_2\) and
\(\|M_{\Omega'}V\|\le J_1\). A rigorous self-adjoint enclosure of
\(B_a^{\mathrm H}\) could therefore narrow the unresolved slab without
estimating a nonnormal resolvent.

After Corollary 5.4, (5.0) and the index formula show that, off the two
indicial lines, \(L_a-\lambda\) is invertible outside the slab when
\(p_0+p_\infty=1\); when \(p_0=p_\infty=0\) it is injective with one-dimensional
cokernel. The case \(p_0=p_\infty=1\) cannot occur outside the slab, because
its Fredholm index is \(+1\).

\subsection{5.2 The local Fredholm anchor and upper-semi-Fredholm
continuation}\label{the-local-fredholm-anchor-and-upper-semi-fredholm-continuation}

Let \(T_0:=L_{\mathrm{loc}}-\lambda\) on the common domain \(D\) of
Corollary 4.4. Its homogeneous solution \[
h_\lambda(y)=\exp\int^y\frac{q-1-\lambda}{b}\,dt
\] has endpoint exponents \[
\mu_0(\lambda)=1-\lambda/\tilde c\quad(y\to0),\qquad
\mu_\infty(\lambda)=-(1+\lambda)/c_l\quad(y\to\infty).
\] Off the two indicial lines define the generic endpoint flags \[
p_0:=\mathbf 1_{\{\operatorname{Re}\lambda<\lambda_X\}},\qquad
p_\infty:=\mathbf 1_{\{\operatorname{Re}\lambda>-1+c_l/2\}}.
\] Thus \(p_0\) records the generic \(H^2\)-admissibility of
\(h_\lambda\) at the origin and \(p_\infty\) its \(L^2\)-admissibility at
infinity. The exceptional value \(\lambda=0\), where the leading origin
power is exactly linear and its second derivative cancels, is treated
separately below.

\noindent\textbf{Lemma 5.2 (local Fredholm anchor, including the trace
resonance).} For every \(\lambda\) off the two indicial lines,
\(T_0:D\to X\) is Fredholm and \[
\dim\ker T_0=p_0p_\infty+\mathbf 1_{\{\lambda=0\}},\qquad
\operatorname{codim}\operatorname{ran}T_0=(1-p_0)(1-p_\infty)+\mathbf 1_{\{\lambda=0\}}.
\] In particular, \[
\operatorname{ind}T_0=p_0+p_\infty-1. \tag{5.1a}
\]

\emph{Proof.} First restrict to \(X_w=\ker\tau\), which \(T_0\)
preserves because the exact trace identity is \[
\tau(T_0\varphi)=-\lambda\,\tau\varphi. \tag{5.2a}
\] Put \[
A_\lambda:=T_0|_{D\cap X_w}:D\cap X_w\longrightarrow X_w.
\] Because \(\tau\in X^*\), \(D\cap X_w\) is closed in the graph space
\(D\). Restricting the compact-remainder estimate of Corollary 4.4 at
\(t=0\) and using
\(L_{\mathrm{loc}}\varphi=A_\lambda\varphi+\lambda\varphi\) gives \[
\begin{aligned}
\|\varphi\|_D
&=\|\varphi\|_X+\|L_{\mathrm{loc}}\varphi\|_X\\
&\le(1+|\lambda|)\|\varphi\|_X+\|A_\lambda\varphi\|_X\\
&\le C\bigl(\|A_\lambda\varphi\|_X+\|\mathcal K\varphi\|_Y\bigr).
\end{aligned}
\] The restriction \(\mathcal K:D\cap X_w\to Y\) remains compact. Peetre's lemma
therefore gives \(A_\lambda\in\Phi_+\); in particular its range is
closed. We now identify that closed range using only compactly supported
data. For \(f\in C_c^\infty(0,\infty)\), variation of constants gives \[
\varphi(y)=h_\lambda(y)\left(C-\int_{y_*}^y\frac{f(s)}{b(s)h_\lambda(s)}\,ds\right). \tag{5.3v}
\] Near either endpoint this solution is either identically zero or a
constant multiple of \(h_\lambda\). At every non-admissible endpoint
choose \(C\) so that it is zero there. Near zero, the smooth coefficient
expansion gives \[
|h_\lambda^{(j)}(y)|\lesssim
y^{\operatorname{Re}\mu_0-j},
\qquad j=0,1,2.
\] Thus \(\operatorname{Re}\mu_0>3/2\) makes \(h_\lambda\) locally
admissible in \(X_w\). Conversely, \(q\) is even, so \(q'(0)=0\).
The \(M_4\) bound and \(b(y)=c_ly+a\int_0^yq(s)\,ds\) give the
quantified expansions \[
|q(y)-q(0)|\le\frac{M_4}{2}y^2,
\qquad
|b(y)-\tilde c\,y|\le\frac{aM_4}{6}y^3.
\] Writing \(m_\lambda=(q-1-\lambda)/b\), these estimates give
\[
m_\lambda(y)=\frac{\mu_0(\lambda)}y+O(y),\qquad
m_\lambda'(y)=-\frac{\mu_0(\lambda)}{y^2}+O(1),
\]
and hence
\[
h_\lambda(y)=C_\lambda y^{\mu_0(\lambda)}(1+O(y^2)),
\quad C_\lambda\ne0.
\tag{5.3o}
\]
The exact identity \(h_\lambda''=(m_\lambda'+m_\lambda^2)h_\lambda\)
therefore yields, if \(\mu_0\notin\{0,1\}\),
\[
h_\lambda''(y)
=C_\lambda\mu_0(\mu_0-1)y^{\mu_0-2}(1+o(1)).
\]
Thus \(h_\lambda''\notin L^2(0,1)\) when
\(\operatorname{Re}\mu_0<3/2\). At \(\mu_0=0\) (equivalently
\(\lambda=\tilde c\)), \(h_\lambda(0)\ne0\). At \(\mu_0=1\)
(equivalently \(\lambda=0\)), \(h_\lambda'(0)\ne0\), so the
mode is not in \(X_w\). The full trace resonance at \(\lambda=0\) is
treated below. Consequently, off \(\Gamma_0\), the homogeneous origin
mode is admissible in \(X_w\) exactly when \(p_0=1\).

At infinity, \(q\in L^1\), \(b=c_ly+O(1)\), and
\(b(y)/y\ge\beta_->0\) give \[
|h_\lambda(y)|\asymp y^{\operatorname{Re}\mu_\infty},
\qquad
|h_\lambda'(y)|\lesssim y^{\operatorname{Re}\mu_\infty-1}.
\] Using the same \(m_\lambda=(q-1-\lambda)/b\) and differentiating the
homogeneous equation, \[
h_\lambda''
=
\left[
\frac{q'}{b}
-\frac{(q-1-\lambda)b'}{b^2}
+\frac{(q-1-\lambda)^2}{b^2}
\right]h_\lambda.
\tag{5.3a}
\] The registry gives \(q'=H\Omega'\in L^2\): local boundedness
controls finite intervals, \(J_1\) controls the tail of the even
whole-line extension of \(\Omega'\), and differentiation commutes with
the Hilbert transform. It also gives \(b'=c_l+aq\in L^\infty\) and
\(b(y)\gtrsim y\). If \(\operatorname{Re}\mu_\infty<-1/2\), the first
term in (5.3a) is controlled in \(L^2\) by
\(q'(y)y^{\operatorname{Re}\mu_\infty-1}\), while the remaining terms
are \(O(y^{\operatorname{Re}\mu_\infty-2})\). Hence
\(h_\lambda,h_\lambda''\in L^2\) near infinity without any pointwise
differentiated far-tail assumption.

Thus the resulting \(\varphi\) lies in \(X_w\), and since
\(T_0\varphi=f\in X\), maximality gives \(\varphi\in D\cap X_w\). If at
least one endpoint is admissible this construction puts every test
function in \(\operatorname{ran}A_\lambda\). Lemma A.1 and the
closedness of the range then give \(\operatorname{ran}A_\lambda=X_w\).

It remains to treat the case \(p_0=p_\infty=0\). The two endpoint constants
agree exactly when the compatibility functional \[
\ell_\lambda(f):=\int_0^\infty\frac{f(s)}{b(s)h_\lambda(s)}\,ds
\] vanishes. This functional extends boundedly and nontrivially to
\(X_w\). Indeed, near zero \[
|f(y)|\le 3^{-1/2}y^{3/2}\|f''\|_{L^2(0,y)},\qquad
|(bh_\lambda)^{-1}|\lesssim y^{-1-\operatorname{Re}\mu_0},
\] and \(\operatorname{Re}\mu_0<3/2\); at infinity
\((bh_\lambda)^{-1}\in L^2(1,\infty)\) precisely because
\(\operatorname{Re}\mu_\infty>-1/2\). Hence \[
|\ell_\lambda(f)|\le C_\lambda\bigl(\|f''\|_{L^2(0,1)}+\|f\|_{L^2(1,\infty)}\bigr). \tag{5.4a}
\] An interior test function can be chosen with
\(\ell_\lambda(\eta)=1\), so \(C_c^\infty\cap\ker\ell_\lambda\) is dense
in \(\ker\ell_\lambda\): if \(f_n\to f\) in \(X_w\), replace \(f_n\) by
\(f_n-\ell_\lambda(f_n)\eta\). Formula (5.3v) therefore shows
\(\ker\ell_\lambda\subseteq\operatorname{ran}A_\lambda\) after taking
the closed-range limit.

Conversely, if \(f=A_\lambda\varphi\), then \[
\left(\frac{\varphi}{h_\lambda}\right)'=-\frac{f}{bh_\lambda}.
\] The right-hand side is integrable at both ends by (5.4a). At zero,
Lemma A.1 gives \(\varphi=o(y^{3/2})\), so \(\varphi/h_\lambda\to0\)
because \(\operatorname{Re}\mu_0<3/2\). At infinity the ratio also tends
to zero: a nonzero limit would imply \(\varphi\sim Ch_\lambda\),
contradicting \(\varphi\in L^2\) when
\(\operatorname{Re}\mu_\infty>-1/2\). Integration over \((0,\infty)\)
gives \(\ell_\lambda(f)=0\). Thus
\(\operatorname{ran}A_\lambda=\ker\ell_\lambda\). The homogeneous
equation contributes \(h_\lambda\) exactly when both endpoints admit it.
We have proved \[
\dim\ker A_\lambda=p_0p_\infty,\qquad
\operatorname{codim}_{X_w}\operatorname{ran}A_\lambda=(1-p_0)(1-p_\infty). \tag{5.5a}
\]

Now use \(D=(D\cap X_w)\oplus\mathbb Ce\) and \(X=X_w\oplus\mathbb Ce\).
With \(B_\lambda:=(I-e\tau)T_0e\), identity (5.2a) gives the triangular
block \[
T_0=\begin{pmatrix}A_\lambda&B_\lambda\\0&-\lambda\end{pmatrix}. \tag{5.6a}
\] For \(\lambda\ne0\) the scalar block is invertible and block
elimination reduces the full kernel and cokernel to those of
\(A_\lambda\).

At \(\lambda=0\) one has \((p_0,p_\infty)=(0,1)\), using \(\tilde c>0\)
and the standing \(c_l<2\). The hypotheses also imply \(\Omega\in X\):
\(\Omega\in L^2\) by \((H_{\mathrm{prof}})\), local boundedness of
\(\Omega''\) follows from smoothness (and is recorded by \(M_2\)), and
\(J_2<\infty\) gives \(\Omega''\in L^2(1,\infty)\). Normalize the
integrating factor so that \(h_0=\Omega\). For a compactly supported
\(f\), the origin-cancelled
solution \[
\varphi(y)=-h_0(y)\int_0^y\frac{f(s)}{b(s)h_0(s)}\,ds
\] vanishes near zero and is a constant multiple of \(h_0=\Omega\)
near infinity, hence belongs to \(D\cap X_w\). Density and the closed
range make \(A_0\) onto; its homogeneous space is \(\mathbb C\Omega\),
which has zero intersection with \(X_w\) because
\(\Omega'(0)\ne0\). Thus \(A_0\) is invertible. The block (5.6a)
now gives \[
\ker T_0=\operatorname{span}\{e-A_0^{-1}B_0\}=\operatorname{span}\{\Omega\},\qquad
\operatorname{ran}T_0=X_w.
\] This adds one to both counts and leaves (5.1a) unchanged. The
\(\Omega\)-mode here belongs to the local anchor \(T_0\); it is distinct
from the \(\lambda=0\) scaling mode \(y\Omega'\) of the full operator
\(L_a\). \(\square\)

\noindent\textbf{Lemma 5.3 (continuation inside \(\Phi_+\)).} Let \(D,X\) be
Hilbert spaces and let \(t\mapsto T_t\in\mathcal B(D,X)\) be
norm-continuous on \([0,1]\). Assume \(T_t\in\Phi_+\) for every \(t\).
If \(T_0\) is Fredholm, then every \(T_t\) is Fredholm and
\(\operatorname{ind}T_t=\operatorname{ind}T_0\).

\emph{Proof.} Fix \(T\in\Phi_+\) and split \(D=\ker T\oplus M\),
\(X=\operatorname{ran}T\oplus Z\). The block
\(T|_M:M\to\operatorname{ran}T\) is an isomorphism. For every
sufficiently small perturbation \(S\), its corresponding block remains
invertible; bounded block elimination reduces \(S\) to this invertible
block plus a map \(G:\ker T\to Z\). Since \(\ker T\) is
finite-dimensional, \(G\) has finite rank. If \(\dim Z<\infty\), then
\(S\) is Fredholm and \[
\operatorname{ind}S=\dim\ker T-\dim Z=\operatorname{ind}T.
\] If \(\dim Z=\infty\), then \(\operatorname{ran}G\) has infinite
codimension in \(Z\), so \(S\) still has infinite defect. Thus the
finite-defect and infinite-defect parts of \(\Phi_+\) are relatively
open and the extended index is locally constant. Connectedness of
\([0,1]\) proves the lemma. \(\square\)

\noindent\textbf{Corollary 5.4 (Fredholmness and index of the full operator).}
Corollary 4.4 gives \(L_a^t-\lambda\in\Phi_+(D,X)\) for all
\(t\in[0,1]\), and Lemma 5.2 gives a Fredholm anchor at \(t=0\). Lemma
5.3 therefore yields \[
L_a-\lambda\ \text{Fredholm},\qquad
\boxed{\ \operatorname{ind}(L_a-\lambda)=p_0+p_\infty-1.\ } \tag{5.7a}
\] The individual kernel and cokernel dimensions may jump by equal
amounts along the homotopy. If condition (ND) of Section 5.3 is imposed,
off the symmetry eigenvalues \(\{0,1\}\), \[
\dim\ker(L_a-\lambda)=p_0p_\infty,\qquad
\dim\operatorname{coker}(L_a-\lambda)=(1-p_0)(1-p_\infty),
\] where the second equality follows from the index formula (5.7a).

\subsection{5.3 Kernel--cokernel refinement by the nondegeneracy
condition (ND)}\label{kernelcokernel-refinement-by-condition-nd}

Set \[
\Lambda_{\mathrm{ND}}:=\mathbb C\setminus
\bigl(\Gamma_0\cup\Gamma_\infty\cup\{0,1\}\bigr).
\] For \(\lambda\in\Lambda_{\mathrm{ND}}\), define the fixed-profile nondegeneracy
condition \[
\mathbf{(ND)}_\lambda\qquad
\dim\ker(L_a-\lambda:D_a\to X)\le p_0p_\infty. \tag{ND}
\] For \(S\subseteq\Lambda_{\mathrm{ND}}\), condition \((ND)_S\) means that
\((ND)_\lambda\) holds for every \(\lambda\in S\). We write (ND) without a
subscript only for \(S=\Lambda_{\mathrm{ND}}\). This operational formulation is the
exact additional input needed below: when \(p_0p_\infty=0\) it excludes a
nonzero global kernel mode invisible to endpoint power counting, and
when \(p_0p_\infty=1\) it asserts geometric simplicity.

\noindent\textbf{An unconditional fixed-profile range.} Lemma 5.1 proves
\((ND)_\lambda\) without any Birman--Schwinger construction whenever
\(\lambda\in\Lambda_{\mathrm{ND}}\) and \(\operatorname{Re}\lambda\notin[E_-,E_+]\):
there the kernel is zero. If \(p_0+p_\infty=1\), Corollary 5.4 further makes
\(L_a-\lambda\) invertible; if \(p_0=p_\infty=0\), its cokernel is
one-dimensional. Thus condition (ND) can remain unresolved only inside
the explicit vertical slab above, not on the whole of \(\Lambda_{\mathrm{ND}}\). To
compare the slab with the inter-line band, write \(F=-1+c_l/2\),
\(O=-\tilde c/2\), and \(g=F-O\). The origin identities give
\(q(0)=1+\tilde c>1\) and \(g=(1-a/2)q(0)-2\). Since
\(q(y)\to0\) by \((D_\infty)\), one has \(q_-\le0\); moreover
\(q_+\ge q(0)\) and \(J_1>0\). Therefore \[
E_-<F,\qquad E_+-O\ge2\tilde c+aJ_1>0.
\] Hence the entire band lies in the slab whenever \(F\le O\). If
\(O<F\), full containment is equivalent to \[
aJ_1-\left(1+\frac a2\right)q_-\ge g.
\] Absent this additional fixed-profile inequality, the unresolved set
is the slab--band intersection; any portion with
\(\operatorname{Re}\lambda<E_-\) is already covered by Lemma 5.1. The
orientation remains a fixed-profile scalar datum.

On each of the two unbounded index-zero components \(C\) of
\(\mathbb C\setminus(\Gamma_0\cup\Gamma_\infty)\), Lemma 5.1 together
with Corollary 5.4 supplies a point \(\lambda_*\in C\) for which
\(L_a-\lambda_*\) is invertible. Equip the common domain \(D_a\) with a
fixed graph norm and put \[
A(\lambda):=L_a-\lambda\in\mathcal B(D_a,X).
\] Then \(A\) is a holomorphic Fredholm family of index zero on \(C\).
The analytic Fredholm theorem for holomorphic operator families
{[}17, Chapter VII, Section 1, Supplementary Note 3{]} therefore gives a
discrete set \(Z_C\subset C\) such that \(A(\lambda)^{-1}\) is
holomorphic on \(C\setminus Z_C\), while at every \(\lambda_0\in Z_C\)
\[
A(\lambda)^{-1}=\sum_{j=1}^{m}(\lambda-\lambda_0)^{-j}B_{-j}+B_0(\lambda),
\qquad \operatorname{rank}B_{-j}<\infty,
\] with \(B_0\) holomorphic. Every \(\lambda_0\in Z_C\) is therefore an
isolated spectral point. If \(\gamma\subset C\setminus Z_C\) is a
positively oriented circle enclosing only \(\lambda_0\) (so
\(\gamma\subset\rho(L_a)\)), its Riesz projection is \[
P_{\lambda_0}=\frac{1}{2\pi i}\int_\gamma(\zeta-L_a)^{-1}\,d\zeta
=-\iota B_{-1},
\] where \(\iota:D_a\hookrightarrow X\) is the continuous inclusion.
Consequently \(P_{\lambda_0}\) is a nonzero finite-rank projection, and
hence \(\lambda_0\in\sigma_{\mathrm{disc}}(L_a)\) with finite algebraic
multiplicity. Conversely, Lemma 5.1 and the index-zero formula make
\(L_a-\lambda\) invertible outside the slab
\(E_-\le\operatorname{Re}\lambda\le E_+\). Thus the exterior
non-invertible set is discrete and contained in that slab. Condition
\((ND)_\lambda\) fails there only at those points. The symmetry eigenvalues
\(0,1\) are two known non-invertible points on the right component and
are excluded from \(\Lambda_{\mathrm{ND}}\).

\noindent\textbf{Corollary 5.5 (dimension split under \((ND)_\lambda\)).} Let
\(\lambda\in\Lambda_{\mathrm{ND}}\) and assume \((ND)_\lambda\). Then \[
\dim\ker(L_a-\lambda)=p_0p_\infty,\qquad
\dim\operatorname{coker}(L_a-\lambda)=(1-p_0)(1-p_\infty). \tag{5.8}
\]

\emph{Proof.} If \(p_0=p_\infty=1\), the index is \(1\), so \(\dim\ker\ge1\);
\((ND)_\lambda\) forces \(\dim\ker=1\) and the cokernel is zero. If
\(p_0p_\infty=0\), \((ND)_\lambda\) forces the kernel to vanish. Substitution
into the index formula (5.7a) gives cokernel dimension zero when exactly
one flag is open and one when both are closed. \(\square\)

The two possible band orientations now have a fixed-profile description
requiring no branch parameter. Let \(S\) be the open strip between the
two lines. If \[
-1+c_l/2<\lambda_X,
\] then \((p_0,p_\infty)=(1,1)\) on \(S\) and the index is \(+1\). Every
point of \(S\) is therefore an eigenvalue even without \((ND)_S\); under
\((ND)_S\) it is geometrically simple and the cokernel vanishes. If
\(\lambda_X<-1+c_l/2\), then \((p_0,p_\infty)=(0,0)\) on \(S\) and the index
is \(-1\); if \((ND)_S\) holds, the kernel vanishes and the cokernel is
one-dimensional, so \(S\) is residual spectrum. Without \((ND)_S\), the
nonzero cokernel is forced by the index but injectivity and hence the
pure residual classification are not determined.

The symmetry eigenvalues \(\{0,1\}\) are separate from (ND). These all-\(a\) symmetry
identities were derived in {[}6, Section 3.2{]} under
\(\mathrm{Adm}(a)\); we record the short computation under the present
hypotheses. Direct differentiation of the fixed profile equation gives
\[
L_a(y\Omega')=0,\qquad L_a\Omega=\Omega+c_ly\Omega'. \tag{5.9}
\] Their membership in the maximal domain follows from the same profile
data: writing \[
y\Omega'=\frac{y(q-1)}{b}\,\Omega,
\] the bounds on \(b/y,q,q'\), the identity \(q''=H(\Omega'')\in L^2\),
and the fixed-profile derivative tails show \(y\Omega'\in X\). Hence
the scaling mode \(\varphi_0=y\Omega'\) and the time-shift mode
\(\varphi_1=\Omega+c_ly\Omega'\) lie in \(X\), and (5.9) puts them in
\(D_a\) with \(L_a\varphi_0=0\) and \(L_a\varphi_1=\varphi_1\).

Moreover \(\tilde c=\lim_{y\downarrow0}b(y)/y>0\) and \(0<c_l<2\),
so both indicial lines lie strictly in the left half-plane. Thus at
\(\lambda=0,1\) one has \((p_0,p_\infty)=(0,1)\) and index zero, while the
nonzero symmetry modes give nontrivial kernels. This is exactly why
\(\{0,1\}\) is excluded from \(\Lambda_{\mathrm{ND}}\).

\subsection{5.4 High-frequency collars and bandwise local
finiteness}\label{high-frequency-collars-and-bandwise-local-finiteness}

Put \[
F:=-1+\frac{c_l}{2},\qquad O:=-\frac{\tilde c}{2},
\] and let \(S\) be the open strip between the two lines.

\pagebreak[3]
\noindent\textbf{Proposition 5.6 (fixed-profile interior-collar vertical tail).}
Assume \(F\ne O\).

\begin{enumerate}
\def\labelenumi{\arabic{enumi}.}
\tightlist
\item
  If \(F<O\), then for every \(0<\delta<O-F\) there is
  \(T_\delta<\infty\) such that \[
  F+\delta\le\operatorname{Re}\lambda<O,\qquad
  |\operatorname{Im}\lambda|\ge T_\delta
  \] implies \[
  \dim\ker_X(L_a-\lambda)=1,\qquad
  \dim\operatorname{coker}_X(L_a-\lambda)=0. \tag{5.10}
  \]
\item
  If \(O<F\), then for every \(0<\delta<F-O\) there is
  \(T_\delta<\infty\) such that \[
  O<\operatorname{Re}\lambda\le F-\delta,\qquad
  |\operatorname{Im}\lambda|\ge T_\delta
  \] implies \[
  \dim\ker_X(L_a-\lambda)=0,\qquad
  \dim\operatorname{coker}_X(L_a-\lambda)=1. \tag{5.11}
  \]
\end{enumerate}

\emph{Proof.} Appendix B proves the corresponding kernel count on the
maximal plain-\(L^2\) transport domain by a Hilbert--Volterra phase
gauge and a uniform first-order minimum-modulus estimate. Every
\(X\)-kernel vector belongs to that plain-\(L^2\) maximal domain. In the
first orientation this gives \(\dim\ker_X\le1\), while the
already-proved \(X\)-index \(+1\) forces equality and zero cokernel. In
the second orientation the plain-\(L^2\) operator is invertible, hence
the \(X\)-kernel vanishes; the \(X\)-index \(-1\) gives the stated
cokernel. \(\square\)

\noindent\textbf{Corollary 5.7 (locally finite failure set in the full band).}
Let \[
\mathcal F_{\mathrm{ND}}:=\{\lambda\in S:(ND)_\lambda\text{ fails}\}.
\] Then \(\mathcal F_{\mathrm{ND}}\) is locally finite in \(S\). More precisely,
there is a locally finite set \(E\subset S\) with
\(\mathcal F_{\mathrm{ND}}\subseteq E\) and condition (ND) holding on
\(S\setminus E\). If \(\lambda_n\in\mathcal F_{\mathrm{ND}}\) and
\(|\operatorname{Im}\lambda_n|\to\infty\), then \[
\operatorname{Re}\lambda_n\longrightarrow F. \tag{5.12}
\]

\emph{Proof.} Equip the common maximal domain \(D_a\) with a fixed graph
norm and set \(T(\lambda)=L_a-\lambda\in\mathcal B(D_a,X)\). This is a
holomorphic Fredholm family of constant index \(+1\) or \(-1\) on the
connected strip \(S\). Proposition 5.6 supplies an anchor
\(\lambda_*\in S\) where (ND) holds.

In the case of Fredholm index \(+1\), choose \(\ell\in D_a^*\) nonzero on the
one-dimensional kernel of \(T(\lambda_*)\) and set \[
\mathcal A(\lambda)u=(T(\lambda)u,\ell u):D_a\to X\oplus\mathbb C.
\] This is an index-zero holomorphic Fredholm family and is invertible
at \(\lambda_*\). In the case of Fredholm index \(-1\), choose
\(v\notin\operatorname{ran}T(\lambda_*)\) and set \[
\mathcal A(\lambda)(u,c)=T(\lambda)u+cv:D_a\oplus\mathbb C\to X;
\] this is again index zero and invertible at the anchor. The analytic
Fredholm theorem, or equivalently the local finite-dimensional Schur
reduction, makes \[
E:=\{\lambda\in S:\mathcal A(\lambda)\text{ is not invertible}\}
\] locally finite. Outside \(E\), the first stabilization makes
\(T(\lambda)\) onto and the second makes it injective; the fixed index
gives exactly (ND). Thus \(\mathcal F_{\mathrm{ND}}\subseteq E\). The inclusion may
be strict because \(E\) depends on the chosen \(\ell\) or \(v\).
Finally, Proposition 5.6 excludes actual failure points at large
imaginary part on every closed sub-band separated from \(F\), which
proves (5.12). \(\square\)

Local finiteness is an interior statement. It does not exclude
accumulation at either boundary line, including at bounded imaginary
part. The only unresolved high-frequency limit is the simultaneous
far-field high-frequency corner \(\operatorname{Re}\lambda\to F\),
\(|\operatorname{Im}\lambda|\to\infty\).

The qualifier ``a positive distance from \(F\)'' is structural in this
argument. Under the logarithmic conjugation in Appendix B, the far-end
coefficient is \(r_\sigma(+\infty)=F-\sigma\), so the frozen far
multiplier has minimum-modulus scale \(|F-\sigma|\) and the inverse and
right-inverse constants in (B.11)--(B.12) deteriorate as
\(\sigma\to F\). The Hilbert--Volterra remainder is \(O(|\tau|^{-1})\),
and the Neumann steps (B.16)--(B.17) are uniform on each fixed collar
\(|\sigma-F|\ge\delta\), while the joint limit \(\sigma\to F\),
\(|\tau|\to\infty\) requires a separate threshold analysis. This corner
is an unresolved uniform estimate for condition (ND); it concerns only
the kernel--cokernel multiplicity refinement, while the Fredholm index,
essential-spectrum formulas, Browder strip, and existence of the
eigenvalue band of Fredholm index \(+1\) remain unchanged.

\subsection{5.5 Status of condition (ND) and spectral
geometry}\label{status-of-condition-nd-and-spectral-geometry}

For the rational \(a=0\) benchmark, set
\[
\Omega_0(y)=-\frac{y}{y^2+\tfrac14}.
\]
Then the identity \[
(L_0-\lambda)(y\Omega_0'+\lambda\Omega_0)
=-\lambda(\lambda-1)\Omega_0
\] and the Hardy analysis of {[}6, Sections 3.3 and 5{]} give the exact
\(a=0\) kernel count. The positive-\(a\) results are those of
Proposition 5.6 and Corollary 5.7.

For \(0<a<1\), Lemma 5.1, Proposition 5.6, and Corollary 5.7 give the
positive-\(a\) energy-slab, high-frequency collar, and local-finiteness
conclusions summarized above.

\noindent\textbf{Fixed-profile spectral geometry.} Figure 1 is a qualitative
schematic of the proved and unresolved regions: the two line
orientations, the nonzero Fredholm index in the open band, the energy
slab, the proved high-frequency closed sub-bands a positive distance
from \(F\), and the unresolved far-field high-frequency corner.

\begin{figure}
\centering
\includegraphics[width=\linewidth,height=0.9\textheight,keepaspectratio]{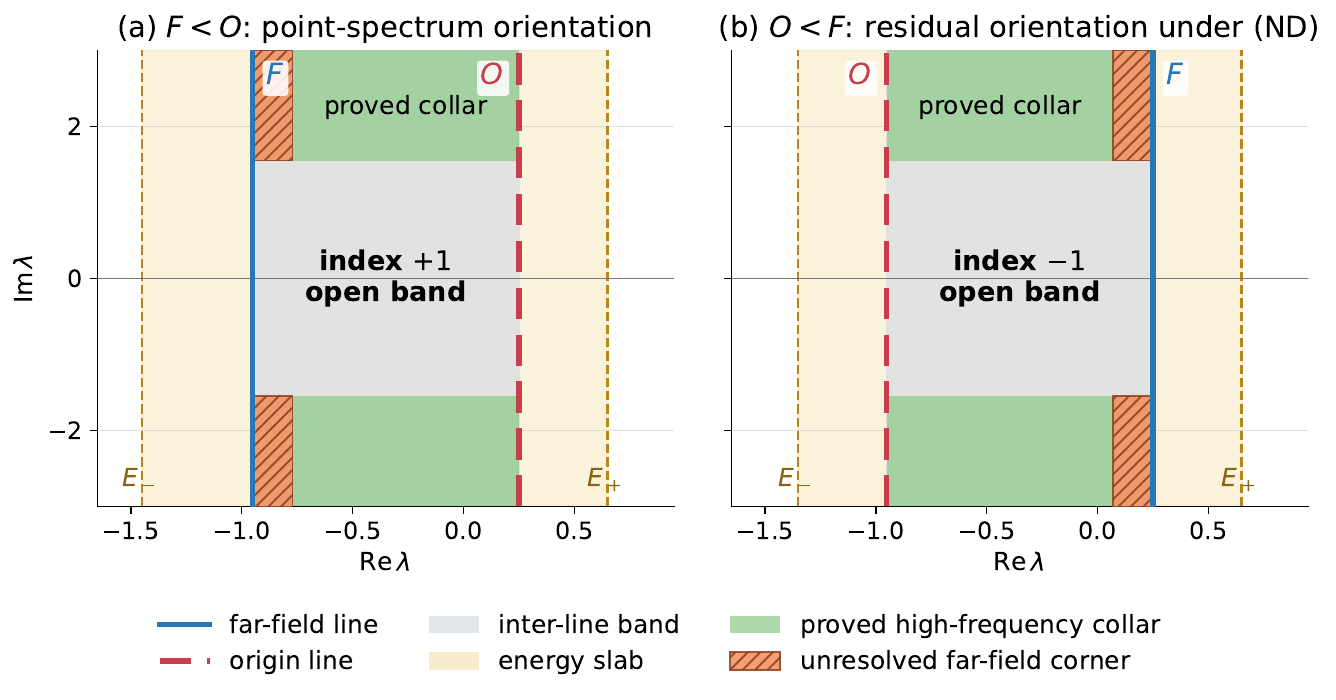}
\caption{Qualitative schematic of the fixed-profile spectral geometry;
no numerical continuation data are shown. In each panel the
solid blue and dashed red vertical lines are the far-field and origin
lines \(F\) and \(O\); the gray strip is the open inter-line band. Panel (a) has index
\(+1\), so every band point is an eigenvalue. Panel (b) has index
\(-1\); injectivity holds on the proved green high-frequency closed
sub-bands. The pale gold vertical region represents the fixed-profile
energy slab, and the orange hatched caps adjacent to \(F\) mark the unresolved
far-field high-frequency corner.}
\end{figure}

\section{6. Semi-Fredholm obstruction on the indicial
union}\label{semi-fredholm-obstruction-on-the-indicial-union}

\noindent\textbf{Lemma 6.1.} Every point of \(\Gamma_0\cup\Gamma_\infty\) belongs
to \(\sigma_{e1}(L_a|_X)\).

\emph{Proof.} Suppose first that the two lines are distinct. On the
common graph domain, \[
\lambda\longmapsto L_a-\lambda\in\mathcal B(D_a,X)
\] is norm-continuous, and the extended index is locally constant on
each of \(\Phi_+\) and \(\Phi_-\). For \(\Phi_+\), this is the block
perturbation argument of Lemma 5.3. For \(\Phi_-\), take Banach
adjoints: \(T\in\Phi_-\) if and only if \(T'\in\Phi_+\), adjunction is
norm-continuous, and the extended indices satisfy
\(\operatorname{ind}T=-\operatorname{ind}T'\), including the infinite
values. The formula (5.7a) changes by one
across each line. If the operator at a line point belonged to either
semi-Fredholm class, nearby Fredholm operators on the two sides would
have the same index, a contradiction.

Now suppose the lines coincide. The primal log-widening packets of
Section 3 are \(X\)-normalized, weakly null, and have
\((L_a-\lambda)\)-defect tending to zero. Hence
\(L_a-\lambda\notin\Phi_+\). To rule out \(\Phi_-\), put \[
\tau=-\frac{\operatorname{Im}\lambda}{c_l},\qquad
R_n=e^{T_n-L_n},
\] and choose \[
u_n(y)=c_n\,\eta\big((\log y-T_n)/L_n\big)y^{-1/2+i\tau},
\qquad L_n\to\infty,\quad R_n\to\infty,
\] with \(\|u_n\|_2=1\). Let \(J_X:L^2(0,\infty)\to X^*\) and
\(J_D:L^2(0,\infty)\to D_a^*\) be the canonical embeddings \[
(J_Xu)(\varphi):=\langle\varphi,u\rangle_{L^2}\quad(\varphi\in X),
\qquad
(J_Du)(\varphi):=\langle\varphi,u\rangle_{L^2}\quad(\varphi\in D_a).
\] Since \(\|u_n''\|_2=O(R_n^{-2})\), \[
1-o(1)=\frac{\|u_n\|_2^2}{\|u_n\|_X}
\le\|J_Xu_n\|_{X^*}\le\|u_n\|_2=1,
\] and \(J_Xu_n\rightharpoonup0\) in \(X^*\) because \(X\) is reflexive
and the supports escape. For smooth compactly supported \(u_n\),
integration by parts in the local terms, \(L^2\) Hilbert-transform
duality, and Fubini in the \(V\) channel give \[
(L_a-\lambda)'J_Xu_n
=J_D\big((L_a^\dagger-\bar\lambda)u_n\big). \tag{6.1a}
\] Here \(L_a^\dagger\) denotes the \(L^2\)-formal adjoint on compactly
supported smooth functions, while the prime is the Banach adjoint
\(X^*\to D_a^*\) of the graph-space operator
\(L_a-\lambda:D_a\to X\). This identity holds on
the whole maximal domain: \(u_n\) is supported away from both endpoints,
and \(W(u_n)\in L^2\) by the estimates below, so every term in these
identities is integrable. No graph-core assertion is used.

The frozen adjoint is \[
A_\infty^\dagger=(-1+c_l)+c_ly\partial_y,
\] and the choice \(\tau=-\operatorname{Im}\lambda/c_l\) gives
\(\|(A_\infty^\dagger-\bar\lambda)u_n\|_2=O(L_n^{-1})\) from the
log-cutoff derivative. The local perturbations are \(o(1)\) because
\(q,\Omega\to0\), \(U\) is bounded, and
\(\|u_n'\|_2=O(R_n^{-1})\). For the remaining \(W\) channel, \[
W(u_n)(y)=\int_y^\infty\Omega'(t)u_n(t)\,dt,
\] Cauchy--Schwarz on \((0,R_n)\) and the dual Hardy inequality on
\((R_n,\infty)\) give \[
\|W(u_n)\|_{L^2(0,R_n)}
\le\left(\int_{R_n}^\infty y|\Omega'(y)|^2\,dy\right)^{1/2}=o(1),
\] \[
\|W(u_n)\|_{L^2(R_n,\infty)}
\le2\|y\Omega'u_n\|_2
\le C L_n^{-1/2}
\left(\int_{R_n}^\infty y|\Omega'(y)|^2\,dy\right)^{1/2}=o(1),
\] where \(|u_n(y)|\le C L_n^{-1/2}y^{-1/2}\). Hence
\(H_{\mathrm{ev}}W(u_n)=o(1)\) in \(L^2\) as well. Consequently \[
\|(L_a-\lambda)'J_Xu_n\|_{D_a^*}
\le \|(L_a^\dagger-\bar\lambda)u_n\|_2=o(1).
\] After normalization, \(J_Xu_n\) is a weakly-null singular sequence
for the Banach adjoint, so the adjoint is not upper-semi-Fredholm.
Closed-range duality gives \(L_a-\lambda\notin\Phi_-\). Thus the merged
line also lies in \(\sigma_{e1}\). \(\square\)

\section{7. Main theorem for a fixed
profile}\label{main-theorem-for-a-fixed-profile}

Fix \(0<a<1\) and one smooth odd self-similar collapse profile \((\Omega,c_l)\)
satisfying \((H_{\mathrm{prof}})\), \((K4^+)\), and \((D_\infty)\) of
Section 2.2. Let \(L_a:D_a\subset X\to X\) be the maximal origin-\(H^2\)
realization, and set \[
\Gamma_\infty:=\{\operatorname{Re}\lambda=-1+c_l/2\},\qquad
\Gamma_0:=\{\operatorname{Re}\lambda=\lambda_X\},\qquad
\lambda_X:=-\tilde c/2,
\] where \(\tilde c=c_l+aq(0)=c_l+a(H\Omega)(0)\). Set
\[
F:=-1+\frac{c_l}{2},\qquad O:=\lambda_X=-\frac{\tilde c}{2}.
\]
Let \(p_0,p_\infty\in\{0,1\}\) be the
origin and infinity admissibility flags defined in Section 5.2.

\noindent\textbf{Theorem 7.1 (two-line Edmunds--Evans exactness, band index, and
Browder classification).} Under the fixed-profile hypotheses above:

\begin{enumerate}
\def\labelenumi{\arabic{enumi}.}
\item
  Both endpoint lines lie in \(\sigma_{e2}(L_a|_X)\). This follows from
  the log-widening singular sequences and the nonlocal \(X\)-norm
  subordination of Section 3.
\item
  For every \(\lambda\notin\Gamma_0\cup\Gamma_\infty\), \(L_a-\lambda\)
  is Fredholm. On compact subsets of the complement, the semi-Fredholm
  and graph estimates are uniform for this fixed operator.
\item
  Off the lines, \[
  \operatorname{ind}(L_a-\lambda)=p_0+p_\infty-1.
  \] Consequently \[
  \sigma_{e1}(L_a|_X)=\sigma_{e2}(L_a|_X)=\sigma_{e3}(L_a|_X)
  =\Gamma_0\cup\Gamma_\infty.
  \] On any component where \(\operatorname{ind}(L_a-\lambda)=+1\),
  Fredholmness forces \(\dim\ker(L_a-\lambda)\ge1\); hence every point of
  that component is an eigenvalue of \(L_a\). Lemma 6.1 supplies the
  semi-Fredholm obstruction, by the index jump
  for separated lines and by primal/dual singular sequences for a merged
  line.
\item
  Put \[
  \overline S:=\{\lambda\in\mathbb C:\min(F,O)\le\operatorname{Re}\lambda\le\max(F,O)\}.
  \] Then, including the merged-line case, \[
  \sigma_{e4}(L_a|_X)=\sigma_{e5}(L_a|_X)
  =\sigma_{\mathrm{B}}(L_a|_X)=\overline S.
  \] Indeed, the boundary line or lines lie in \(\sigma_{e2}\). If
  \(F\ne O\), the open band has index \(\pm1\), so it lies in
  \(\sigma_{e4}\) and is a non-isolated open subset of the spectrum,
  hence lies in \(\sigma_{e5}\). On either exterior component the index
  is zero, and the analytic-Fredholm argument of Section 5.3 puts every
  exterior spectral point in \(\sigma_{\mathrm{disc}}(L_a)\) as defined
  in Section 2.3. For profiles satisfying both the present hypotheses
  and \(\mathrm{Adm}(a)\) of {[}6{]}, the first three Edmunds--Evans
  spectra are the two boundary lines. The question in
  {[}6, Section 4.5{]}, formulated under the Browder convention used
  there, has the following conditional negative answer: for every such
  profile with \(F\ne O\), the Browder set is the strictly larger closed
  strip.
\item
  The right essential edge and its distance from the imaginary axis are
  the fixed-profile quantities \[
  \sup\operatorname{Re}(\Gamma_0\cup\Gamma_\infty)
  =\max(-1+c_l/2,\lambda_X)
  =-1+c_l/2+\max(0,-g),
  \] \[
  \operatorname{gap}_{\mathrm{ess}}=1-c_l/2-\max(0,-g),
  \qquad g=\frac{c_l+\tilde c}{2}-1.
  \] \emph{Proof architecture.} Sections 3--4 give the line inclusions
  and an upper-semi-Fredholm estimate. Corollary 4.4 places the
  common-domain path in \(\Phi_+\). Lemma 5.2 computes the Fredholm
  defect and index of the local endpoint problem, and Lemma 5.3 carries
  that data to the full operator. Section 6 treats a merged line. Item 4
  combines the line inclusion and band index with the finite-rank
  meromorphic exterior resolvent proved in Section 5.3. Proposition 7.2
  adds the kernel--cokernel refinements. \(\square\)
\end{enumerate}

\noindent\textbf{Proposition 7.2 (condition (ND) refinements).} Under the
fixed-profile hypotheses of Theorem 7.1:

\begin{enumerate}
\def\labelenumi{\arabic{enumi}.}
\item
  If condition (ND) of Section 5.3 holds at \(\lambda\in\Lambda_{\mathrm{ND}}\), then
  \[
  \dim\ker(L_a-\lambda)=p_0p_\infty,
  \qquad
  \dim\operatorname{coker}(L_a-\lambda)=(1-p_0)(1-p_\infty).
  \] Thus the band of Fredholm index \(+1\) is geometrically simple point
  spectrum, and the band of Fredholm index \(-1\) is pure residual
  spectrum at those points.
  Without (ND), positive index still forces a nontrivial kernel, but
  geometric multiplicity and cokernel are not determined.
\item
  Condition (ND) is automatic for
  \(\operatorname{Re}\lambda\notin[E_-,E_+]\) by Lemma 5.1. If the
  far-field line lies to the left of the origin line, the entire
  band of Fredholm index \(+1\)
  lies in this slab; every point is nevertheless an eigenvalue by
  Theorem 7.1. This full-band slab containment is only an energy
  localization and does not itself prove (ND); the actual in-band
  refinement is item 3. In the reverse orientation, the part not settled
  by Lemma 5.1 is precisely the slab--band intersection. On every
  portion with \(\operatorname{Re}\lambda<E_-\), the exact dimension
  statement above is unconditional.
\item
  If \(F\ne O\), Proposition 5.6 proves condition (ND) at sufficiently
  large \(|\operatorname{Im}\lambda|\) on every closed inter-line collar
  a positive distance from \(F\). Corollary 5.7 shows that the actual
  failure set is locally finite in the open band and that every failure
  sequence with unbounded imaginary part approaches \(F\). Accumulation
  at either boundary line, including at bounded height, is not excluded.
\end{enumerate}

\emph{Proof.} The dimension statement is the fixed-profile
nondegeneracy argument of Section 5.3 combined with the index formula.
Lemma 5.1 gives the energy-slab conclusion. Proposition 5.6 and
Corollary 5.7 give the collar and analytic-stabilization conclusions.
\(\square\)

\noindent\textbf{Corollary 7.3 (actual positive-advection realizations).} For
every \(0<a<\underline a=400/(848-9\pi^2)\) and every
Huang--Qin--Wang--Wei fixed point \(f\in\mathbb D\), the normalization
in Proposition 2.7 produces a profile to which Theorem 7.1 applies.
There is \(a_{\mathrm{sep}}\in(0,1/2]\) such that, for every
\(0<a<a_{\mathrm{sep}}\) and every such fixed point, \[
F<O,
\] so the open band has index \(+1\), lies in the Browder essential
spectrum, and consists of eigenvalues. The fixed-point existence theorem
{[}5, Theorem 3.11{]} supplies at least one actual profile at every
parameter in these ranges. Theorem 7.1 gives the eigenvalue-band
conclusion, while Proposition 7.2 adds multiplicity information. For any
such profile that also satisfies \(\mathrm{Adm}(a)\) of {[}6{]}, the
conclusion gives the corresponding conditional answer to the question
posed in {[}6, Section 4.5{]}: under the Browder convention used there,
the essential spectrum is the closed strip \(\overline S\), while
\(\sigma_{e1}\), \(\sigma_{e2}\), and \(\sigma_{e3}\) consist only of
its boundary lines.

\emph{Proof.} Proposition 2.7 supplies the full fixed-profile hypothesis registry, Proposition
2.8 supplies the uniform all-fixed-point inequality \(F<O\) for
sufficiently small \(a\), and Theorem 7.1 gives the spectral
conclusions. \(\square\)

\noindent\textbf{Remaining (ND) problem.} Sections 5.4--5.5 reduce the
unresolved cases to threshold analysis near the far-field line and a locally
finite set of interior candidates.

\section{8. Discussion}\label{discussion}

\noindent\textbf{What the fixed-profile theorem proves.} The two endpoint models
determine both the essential lines and the index. The central analytic
estimate controls the full nonlocal operator in the common graph domain
and places the homotopy in the upper-semi-Fredholm class. The local
transport problem can then be solved and counted directly at its two
endpoints; continuation within \(\Phi_+\) yields full Fredholmness and
preserves the index. This route uses primal common-domain estimates and
a direct local anchor.

\noindent\textbf{Actual profiles and remaining family questions.} Reference {[}5,
Theorem 3.11{]} supplies fixed points, and Appendix C transfers every
one of them with \(0<a<\underline a\) to the complete fixed-profile
hypothesis registry after
normalization. It further proves \(F<O\) for every such fixed point at
all sufficiently small positive \(a\). Thus the line separation holds
uniformly over the fixed-point set for all sufficiently small positive
\(a\). Along the Huang--Qin--Wang--Wei small-\(a\) limit, \(\mu\to2\log2-1\), and
therefore the realized band width satisfies
\(O-F=a(1-2\mu)/(1-a\mu)\sim(3-4\log2)a\approx0.2274a\). The band is
narrow to first order, but rigorously nonzero for all sufficiently small
positive \(a\). Family-level uniqueness, continuous selection, and
uniform registry constants remain separate questions.

\noindent\textbf{Interpretation of the two lines.} After \(y=e^t\), the transport
skeleton has two cylindrical ends. The far-field line is the plain-\(L^2\)
dilation threshold at \(t\to+\infty\), while the origin line is shifted
by the second-derivative weight at \(t\to-\infty\). Endpoint
subordination, rather than compact perturbation, controls the Hilbert
coupling. This is the sense in which the argument parallels two-ended
Lockhart--McOwen/Melrose Fredholm theory while remaining self-contained
for the present problem.

\noindent\textbf{Browder interpretation.} The exact Browder essential spectrum is
the closed inter-line strip, while the first three Edmunds--Evans
spectra remain exactly the boundary line or lines. Thus, for a
separated-line profile satisfying both hypothesis packages, the two-line
description posed in {[}6, Section 4.5{]} is correct for
\(\sigma_{e1}\), \(\sigma_{e2}\), and \(\sigma_{e3}\), but not for the
Browder essential spectrum under the convention adopted there; the
latter is \(\overline S\).

This does not conflict with the numerical report in {[}6{]} of no
resolution-stable discrete eigenvalues: in the \(F<O\) orientation the
present \(X\)-eigenvalues fill an open band, so none is discrete in the
isolated finite-algebraic-multiplicity sense, and finite-grid candidates
may drift rather than converge. It does show that the perturbative
program proposed in {[}6, Discussion item (i){]} cannot exclude all
point spectrum on the origin-\(H^2\) realization. The maximal-\(L^2\)
smear and the present \(X\)-band belong to different realizations, but
both reflect the same two-ended indicial geometry with different
endpoint weights.
Proposition 2.8 rigorously gives the orientation \(F<O\) for every
sufficiently small positive-\(a\) Huang--Qin--Wang--Wei fixed point. The transition value
\(a_c\approx0.6890665\), at which \(c_l\) changes sign and the
self-similar spatial extent changes from shrinking (focusing) to
expanding, is reported in {[}4{]}, while the crossing \(F=O\) for the
exact \(a=1/2\) profile is discussed in {[}6{]}. Different profiles at the same
\(a\) may have different scalar quantities \((c_l,q(0),\tilde c)\).

\section{Appendix A. Trace splitting and pointwise profile
data}\label{appendix-a.-trace-splitting-and-pointwise-profile-data}

\subsection{Lemma A.1 (density of the test class in the gauge
subspace)}\label{lemma-a.1-density-of-the-test-class-in-the-gauge-subspace}

Let \(\tau\varphi=\varphi'(0)\) and \(X_w=\ker\tau\). Then
\(C_c^\infty(0,\infty)\) is dense in \(X_w\), while it is not dense in
\(X\). Moreover, every \(\varphi\in X_w\) satisfies
\[
\varphi(y)=o(y^{3/2}),\qquad \varphi'(y)=o(y^{1/2})
\qquad(y\downarrow0).
\]

\emph{Proof.} If \(\varphi\in X_w\), then \(\varphi(0)=\varphi'(0)=0\)
and \[
\varphi(y)=\int_0^y(y-t)\varphi''(t)\,dt,\qquad
\varphi'(y)=\int_0^y\varphi''(t)\,dt.
\] Writing \(\varepsilon(y)=\|\varphi''\|_{L^2(0,y)}\to0\) gives \[
|\varphi(y)|\le3^{-1/2}y^{3/2}\varepsilon(y),\qquad
|\varphi'(y)|\le y^{1/2}\varepsilon(y).
\] For a cutoff that changes from \(0\) to \(1\) on
\([\delta,2\delta]\), the three terms in \[
(\chi_\delta\varphi)''-\varphi''=\chi_\delta''\varphi+2\chi_\delta'\varphi'+(\chi_\delta-1)\varphi''
\] are respectively \(O(\varepsilon(2\delta))\),
\(O(\varepsilon(2\delta))\), and \(o(1)\) in \(L^2\). A far cutoff tends
to the identity because \(\varphi,\varphi',\varphi''\in L^2\), and
mollification away from the endpoints produces \(C_c^\infty(0,\infty)\)
approximants in \(X\). Conversely, \(\tau\) is continuous on \(X\) and
vanishes on every compactly supported test function, so the closure of
that class cannot contain any vector with nonzero derivative trace.
\(\square\)

\subsection{Lemma A.2 (boundary-free Riesz
representation)}\label{lemma-a.2-boundary-free-riesz-representation}

Fix \(e=y\chi_0\in D_a\) with \(\tau e=1\). Every \(\Lambda\in X^*\) has
a unique decomposition \[
\Lambda=c\,\tau+\Lambda_w,\qquad c=\Lambda(e),\qquad \Lambda_w(e)=0,
\] and there are \(v_0,v_2\in L^2(0,\infty)\) such that \[
\Lambda_w(\varphi)=\langle\varphi,v_0\rangle+\langle\varphi'',v_2\rangle
\qquad(\varphi\in X).
\] No boundary value of \(v_2\) is asserted.

\emph{Proof.} Equip \(X\) with the equivalent Hilbert inner product \[
\langle\varphi,\psi\rangle_X=\langle\varphi,\psi\rangle_2+\langle\varphi'',\psi''\rangle_2.
\] Apply Riesz representation to \(\Lambda_w=\Lambda-\Lambda(e)\tau\) on
all of \(X\): for a unique \(\psi_w\in X\), \[
\Lambda_w(\varphi)=\langle\varphi,\psi_w\rangle_2+\langle\varphi'',\psi_w''\rangle_2.
\] Take \(v_0=\psi_w\) and \(v_2=\psi_w''\). If the resulting
distribution \(u=v_0+v_2''\) vanishes on tests, it kills \(X_w\) by
Lemma A.1 and continuity; together with \(\Lambda_w(e)=0\) this kills
all of \(X\). Hence \(\Lambda\mapsto(c,u)\) is injective. \(\square\)

\subsection{A.3. Distributional transposition for range
annihilators}\label{a.3.-distributional-transposition-for-range-annihilators}

Every range annihilator satisfies the following distributional
transposition identity, obtained from the decomposition in Lemma A.2.

Let \(T_\lambda=L_a-\lambda:D_a\to X\), and let
\(\Lambda\in\operatorname{Ann}\operatorname{ran}T_\lambda\). Write \[
\Lambda=c_0\tau+\Lambda_w,\qquad c_0=\Lambda(e),\qquad \Lambda_w(e)=0,
\] as in Lemma A.2, and associate to \(\Lambda_w\) the distribution \[
\langle u,\phi\rangle_{\mathcal D',\mathcal D}
:=\overline{\Lambda_w(\overline\phi)},
\qquad \phi\in C_c^\infty(0,\infty). \tag{A.3.1}
\] The density of \(C_c^\infty(0,\infty)\) in \(X_w\), together with
\(\Lambda_w(e)=0\), shows that \(u\) determines \(\Lambda_w\) uniquely.
For distributions arising in this way define the normalized transposed
action \(\mathcal T_{\lambda,e}\) by \[
\left\langle\mathcal T_{\lambda,e}u,\phi\right\rangle
:=\overline{\Lambda_w\!\left((L_a-\lambda)\overline\phi\right)}. \tag{A.3.2}
\] This is a definition by transposition; in particular, it does not
assign a pointwise meaning to a Volterra integral \(W(u)\) for a general
distribution. Because \(\Lambda_w(e)=0\) is built into this
normalization, no identification with a bare \(L^2\) formal-adjoint
realization is asserted.

Put \(d=\overline{c_0}\) and \(m(y)=-2/(\pi y)\). Since \[
\tau\!\left((L_a-\lambda)\phi\right)
=(1-a)\Omega'(0)\int_0^\infty m(y)\phi(y)\,dy
\qquad(\phi\in C_c^\infty(0,\infty)),
\] every range annihilator necessarily satisfies \[
\mathcal T_{\lambda,e}u=-d(1-a)\Omega'(0)m
\qquad\text{in }\mathcal D'(0,\infty). \tag{A.3.3}
\] It must additionally satisfy \(\Lambda((L_a-\lambda)e)=0\). Thus
every range annihilator gives a distributional solution of (A.3.3).
This forward implication is recorded for completeness; no converse is
claimed, and it is not used in the Fredholm argument above.

\subsection{\texorpdfstring{Registry A.4 (fixed-profile \(K4^+\)
bounds)}{Registry A.4 (fixed-profile K4+ bounds)}}\label{registry-a.4-fixed-profile-k4k4-bounds}

Fix \(0<a<1\) and one profile satisfying Section 2.2. The registry
records the pointwise hypotheses needed by the operator estimates for
this fixed profile. Appendix C proves the full fixed-profile hypothesis registry for every
Huang--Qin--Wang--Wei fixed point with \(0<a<\underline a\), after
normalization. Other profiles require individual verification.

{\def\LTcaptype{none} 
\begin{longtable}[]{@{}
  >{\raggedright\arraybackslash}p{(\linewidth - 4\tabcolsep) * \real{0.3333}}
  >{\raggedright\arraybackslash}p{(\linewidth - 4\tabcolsep) * \real{0.3333}}
  >{\raggedright\arraybackslash}p{(\linewidth - 4\tabcolsep) * \real{0.3333}}@{}}
\toprule\noalign{}
\begin{minipage}[b]{\linewidth}\raggedright
quantity
\end{minipage} & \begin{minipage}[b]{\linewidth}\raggedright
definition
\end{minipage} & \begin{minipage}[b]{\linewidth}\raggedright
status in this paper
\end{minipage} \\
\midrule\noalign{}
\endhead
\bottomrule\noalign{}
\endlastfoot
\(\beta_-,\beta_+\) & \(0<\beta_-\le b(y)/y\le\beta_+\) & assumed in
\((H_{\mathrm{prof}})\) \\
\(r_{\mathrm{org}},M_1,M_2,M_3,M_4\) & fixed origin radius and the local bounds on
\(q',\Omega'',(\Omega-\Omega'(0)y)/y^3,q''\) & assumed in
\((H_{\mathrm{prof}})\) \\
\(\Omega'(0)=-4\) & gauge \(\Omega'(0)\) & normalization of the chosen
profile \\
\(Q_0,Q_1,W_0,W_1,W_2\) & global \(L^\infty\) coefficient and derivative
bounds & assumed in \((K4^+)\) \\
\(I_q\) & \(\int_0^\infty|q|\) & assumed in \((K4^+)\) and compatible
with \((D_\infty)\) \\
\(J_k\) & \((\int_0^\infty y|\Omega^{(k)}|^2)^{1/2}\), \(k=1,2,3\) &
assumed in \((K4^+)\) \\
\(\varepsilon_\Omega,\varepsilon_q,\varepsilon_V\) & tail moduli tending
to zero & defined after \((D_\infty)\); their vanishing follows from
\((K4^+)\) and \((D_\infty)\) \\
\(C_q,C_\Omega,\alpha_{\rm dec},U_\infty\) & far-field constants in Section 2.2 &
assumed in \((D_\infty)\) \\
\end{longtable}
}

Under these hypotheses, every coefficient constant used in Sections 3--5
is a finite algebraic function of the displayed registry. In particular,
Corollary 4.4 uses \[
C_P=2W_0+\sqrt2W_1+W_2
+a\left(J_1+2^{-1/2}W_1+2W_2+J_3\right),
\] which is finite by \((K4^+)\). A compact-uniform family theorem would
additionally require uniform bounds for every registry entry and
suitable continuity of the supplied family.

\noindent\textbf{Relation to {[}5{]}.} Reference {[}5, Theorem 3.11{]} supplies
fixed-point existence at each parameter in its stated regime; Appendix C
supplies the pointwise registry transfer for every fixed point with
\(0<a<\underline a\).

\begin{center}\rule{0.5\linewidth}{0.5pt}\end{center}

\section{Appendix B. High-frequency Hilbert--Volterra
gauge}\label{appendix-b.-high-frequency-hilbertvolterra-gauge}

This appendix proves the plain-\(L^2\) estimate used in Proposition 5.6.
Throughout, \(\lambda=\sigma+i\tau\) with \(\sigma,\tau\in\mathbb R\).
Let \[
D_b:=\{w\in L^2(0,\infty):bw'\in L^2(0,\infty)\},\qquad
\|w\|_{D_b}=\|w\|_2+\|bw'\|_2. \tag{B.1}
\] The standing bounds on \(b/y\) identify this with the maximal
plain-\(L^2\) transport domain. Define \[
\rho(y)=\int_1^y\frac{dr}{b(r)},\qquad
\Gamma(y)=\int_1^y\frac{\Omega(r)}{b(r)}\,dr,
\] \[
P_\tau(y)=\exp\{-i\tau\rho(y)+i\operatorname{sgn}(\tau)\Gamma(y)\},\qquad
A_\sigma=-b\partial_y+q-1-\sigma. \tag{B.2}
\] The integral defining \(\Gamma\) has finite limits at both endpoints,
and multiplication by \(P_\tau\) maps \(D_b\) onto itself for each fixed
\(\tau\).

\noindent\textbf{Lemma B.1 (Hilbert--Volterra plateau gauge).} For \(T=|\tau|\)
sufficiently large, \[
P_\tau^{-1}(L_a-(\sigma+i\tau))P_\tau=A_\sigma+R_\tau, \tag{B.3}
\] where \[
R_\tau w=
\Omega P_\tau^{-1}\bigl(H_{\mathrm{odd}}(P_\tau w)-i\operatorname{sgn}(\tau)P_\tau w\bigr)
-a\Omega'P_\tau^{-1}V(P_\tau w), \tag{B.4}
\] and \[
\|R_\tau w\|_2\le\frac{C}{T}\|w\|_{D_b}. \tag{B.5}
\] If \(\sigma\) lies in a bounded interval \(I\), then \[
\|R_\tau w\|_2\le\frac{C_I}{T}
\bigl(\|w\|_2+\|A_\sigma w\|_2\bigr). \tag{B.6}
\] The constant in (B.5) depends only on the fixed-profile registry;
\(C_I\) in (B.6) may also depend on the bounded interval \(I\).

\emph{Proof.} The identity \[
\frac{P_\tau'}{P_\tau}=
\frac{-i\tau+i\operatorname{sgn}(\tau)\Omega}{b}
\] shows that the transport contribution cancels both the spectral term
\(i\tau\) and the leading Hilbert plateau
\(i\operatorname{sgn}(\tau)\Omega\), giving (B.3)--(B.4).

For the remainder, use the logarithmic unitary
\((U_0f)(t)=e^{t/2}f(e^t)\). We fix \[
D_t=-i\partial_t,\qquad
\widehat g(\xi)=\int_{\mathbb R}e^{-it\xi}g(t)\,dt.
\] The Fourier mode \(e^{i\xi t}\) corresponds to
\(f(y)=y^{-1/2+i\xi}=y^{-s}\) with \(s=1/2-i\xi\). Therefore the Mellin
symbol \(-\cot(\pi s/2)\) gives \[
U_0H_{\mathrm{odd}}U_0^{-1}=m(D_t),\qquad
m(\xi)=-\operatorname{sech}(\pi\xi)-i\tanh(\pi\xi). \tag{B.7}
\] This identity is first obtained on Schwartz functions by the
Fourier--Mellin calculation above. Since \(H_{\mathrm{odd}}\) and
\(m(D_t)\) are bounded on \(L^2\), it extends to all of \(L^2\) by
density. The frequency reversal is the reason for the sign of the
\(\tanh\) term.

Writing \(\kappa(t)=b(e^t)/e^t\), the registry gives
\(0<\beta_-\le \kappa\le\beta_+\) and
\(\dot\kappa=c_l+aq-\kappa\in L^\infty\). Hence \(\dot\rho=1/\kappa\) is
bounded below and \(\ddot\rho\in L^\infty\). For
\(\alpha(t)=e^{i\operatorname{sgn}(\tau)\Gamma(e^t)}U_0w(t)\), \[
\|\alpha\|_{H^1}\le C\bigl(\|w\|_2+\|bw'\|_2\bigr).
\]

We spell out the cutoff estimate used here. Let
\(\epsilon=\operatorname{sgn}(\tau)\), \(h=T^{-1}\), and
\(u=e^{-i\epsilon T\rho}\alpha\). Then \[
(\epsilon hD_t+\dot\rho)u
=\epsilon h e^{-i\epsilon T\rho}D_t\alpha. \tag{B.7a}
\] Choose a fixed smooth semiclassical frequency cutoff
\(0\le\chi_-\le1\), equal to one where
\(\epsilon\xi\le-3/(4\beta_+)\) and supported where
\(\epsilon\xi\le-1/(2\beta_+)\), and put \(\chi_+=1-\chi_-\). On the
whole support of \(\chi_-(h\cdot)\), \[
\|(m(D_t)-i\epsilon)\chi_-(hD_t)u\|_2
\le 4e^{-\pi/(2\beta_+h)}\|u\|_2, \tag{B.7b}
\] because \(|m(\eta)-i\epsilon|\le4e^{-\pi|\eta|}\) whenever
\(\epsilon\eta\le0\). Here \(\eta\) is the unscaled Fourier variable and
the cutoff condition is \(\epsilon h\eta\le-1/(2\beta_+)\).

We estimate the commutator directly. Put \(a_0(t)=\dot\rho(t)\) and
\(g=\chi_+'\). The function \(g\) is fixed, smooth, and compactly
supported, so \(\|\check g\|_1\) is part of the constant below. We use
\(\check g(z)=(2\pi)^{-1}\int e^{iz\xi}g(\xi)\,d\xi\). Distributionally,
\[
\mathcal F^{-1}\chi_+
=c\delta_0+i\,\mathrm{p.v.}\!\frac{\check g(z)}z.
\] The delta term is a multiple of the identity and drops out of the
commutator. For Schwartz functions, the remaining Fourier kernel gives
\[
[a_0,\chi_+(hD_t)]u(t)
=i\int_{\mathbb R}
\frac{a_0(t)-a_0(s)}{t-s}
\check g\!\left(\frac{t-s}{h}\right)u(s)\,ds. \tag{B.7c}
\] The formula extends by density. Since \(a_0'=\ddot\rho\in L^\infty\),
Schur's test yields the explicit bound \[
\|[a_0,\chi_+(hD_t)]u\|_2
\le h\|\ddot\rho\|_\infty\|\check g\|_1\|u\|_2. \tag{B.7d}
\]

On the complementary region put \(v=\chi_+(hD_t)u\). Its Fourier support
and \(\dot\rho\ge1/\beta_+\) give \[
\operatorname{Re}\langle(\epsilon hD_t+\dot\rho)v,v\rangle
\ge\frac{1}{4\beta_+}\|v\|_2^2.
\] Because \(D_t\) commutes with the cutoff, \[
(\epsilon hD_t+\dot\rho)v
=\chi_+(hD_t)(\epsilon hD_t+\dot\rho)u
+[\dot\rho,\chi_+(hD_t)]u.
\] Together with (B.7a), (B.7d), and \(\|u\|_2=\|\alpha\|_2\), this
yields \[
\|\chi_+(hD_t)u\|_2
\le Ch\|\alpha\|_{H^1}.
\] The boundedness of \(m(D_t)-i\epsilon\) on this component, together
with the exponentially small first component, proves \[
\|H_{\mathrm{odd}}(P_\tau w)-i\operatorname{sgn}(\tau)P_\tau w\|_2
\le\frac{C}{T}\|w\|_{D_b}. \tag{B.8}
\]

It remains to control the Volterra channel. Set \(Jf(y)=\int_0^yf\) and
\(r_\tau=H_{\mathrm{odd}}(P_\tau w)-i\operatorname{sgn}(\tau)P_\tau w\).
Then \[
V(P_\tau w)=i\operatorname{sgn}(\tau)J(P_\tau w)+Jr_\tau.
\] Because \(P_\tau'=-i\operatorname{sgn}(\tau)(T-\Omega)P_\tau/b\),
integration by parts gives \[
J(P_\tau w)=i\operatorname{sgn}(\tau)\frac{b}{T-\Omega}P_\tau w
-i\operatorname{sgn}(\tau)J\left[P_\tau\left(\frac{bw}{T-\Omega}\right)'\right]. \tag{B.9}
\] The lower endpoint term vanishes: if \(u(t)=e^{t/2}w(e^t)\), then
\(w,bw'\in L^2\) imply \(u\in H^1(\mathbb R)\subset L^\infty\), and
therefore \(b(y)w(y)\to0\) as \(y\downarrow0\). For
\(T\ge2\|\Omega\|_\infty\), \[
\left(\frac{bw}{T-\Omega}\right)'
=\frac{bw'}{T-\Omega}+\frac{b'w}{T-\Omega}
+\frac{b\Omega'w}{(T-\Omega)^2}.
\] Consequently, \[
\left\|\left(\frac{bw}{T-\Omega}\right)'\right\|_2
\le\frac2T\left(\|bw'\|_2+\|b'\|_\infty\|w\|_2\right)
+\frac4{T^2}\|b\Omega'\|_\infty\|w\|_2. \tag{B.9a}
\] Here \(\|b'\|_\infty\le c_l+aQ_0\) and, by the profile identity
\(b\Omega'=(q-1)\Omega\),
\(\|b\Omega'\|_\infty\le W_0(1+Q_0)\). The Volterra map obeys \[
\|\Omega'Jg\|_2^2
\le\int_0^\infty y|\Omega'(y)|^2\,dy\,\|g\|_2^2
=J_1^2\|g\|_2^2. \tag{B.9b}
\] Multiplying (B.9) by \(\Omega'\), using (B.9a)--(B.9b), and
estimating the first term in (B.9) directly gives \[
\|\Omega'J(P_\tau w)\|_2
\le\frac CT\|w\|_{D_b}. \tag{B.9c}
\] The same Volterra bound and (B.8) give
\(\|\Omega'Jr_\tau\|_2\le CT^{-1}\|w\|_{D_b}\). Therefore \[
\|\Omega'V(P_\tau w)\|_2\le\frac{C}{T}\|w\|_{D_b}. \tag{B.10}
\] Equations (B.8) and (B.10) prove (B.5). Finally, \[
\|bw'\|_2\le\|A_\sigma w\|_2+(1+Q_0+|\sigma|)\|w\|_2,
\] which proves (B.6). \(\square\)

The remaining collar count is recorded separately.

\noindent\textbf{Lemma B.2 (uniform transport collar count).} On every compact
\(J\Subset(-\infty,F)\), the operators \(A_\sigma:D_b\to L^2\) are
invertible and \[
\sup_{\sigma\in J}\|A_\sigma^{-1}\|_{L^2\to D_b}<\infty. \tag{B.11}
\] On every compact \(I\Subset(F,\tfrac32\tilde c)\), they are
surjective with one-dimensional kernel and possess right inverses
\(Q_\sigma\) satisfying \[
\sup_{\sigma\in I}\|Q_\sigma\|_{L^2\to D_b}<\infty. \tag{B.12}
\]

\emph{Proof.} Put \(x=\rho(y)\) and \((Sw)(x)=b(y(x))^{1/2}w(y(x))\).
The registry bounds imply that \(S\) is unitary from \(L^2(dy)\) to
\(L^2(dx)\) and maps \(D_b\) isomorphically onto \(H^1(\mathbb R)\).
Direct differentiation gives \[
SA_\sigma S^{-1}=-\partial_x+r_\sigma(x),\qquad
r_\sigma=q-1-\sigma+\frac{b'}2, \tag{B.13}
\] with \[
r_\sigma(-\infty)=\frac32\tilde c-\sigma,\qquad
r_\sigma(+\infty)=F-\sigma. \tag{B.14}
\] The graph-norm identification is explicit: \[
(Sw)'=S(bw')+\frac12b'(y(x))Sw.
\] Together with \(b'\in L^\infty\), this identity and its inverse prove
\(S(D_b)=H^1(\mathbb R)\) with equivalent graph norms.

Replace \(r_\sigma\) by the step function having the two limits in
(B.14). Half-line variation of constants shows that the resulting
operator \(H^1\to L^2\) is invertible when both limits are positive, and
is surjective with a one-dimensional kernel when the left limit is
positive and the right limit is negative. The difference between
\(r_\sigma\) and this step function is bounded and tends to zero at both
ends; it may jump at \(x=0\). Multiplication by any such coefficient is
compact from \(H^1\) to \(L^2\), by Rellich on a fixed compact interval
and a small \(L^\infty\) tail outside it.
Consequently (B.13) is Fredholm with index zero for \(\sigma<F\) and
index one for \(F<\sigma<\tfrac32\tilde c\).

The dimensions can also be read directly. Set \[
\Phi_\sigma(x)=\exp\!\left(\int_0^x r_\sigma(s)\,ds\right). \tag{B.15}
\] Every homogeneous solution of \((-\partial_x+r_\sigma)u=0\) is a
multiple of \(\Phi_\sigma\), while every homogeneous adjoint solution is
a multiple of \(\Phi_\sigma^{-1}\). If \(\sigma<F\), \(\Phi_\sigma\)
decays at \(-\infty\) and grows at \(+\infty\), whereas its reciprocal
has the opposite behavior. Neither belongs to \(L^2\), so the Fredholm
operator and its adjoint have zero kernel; it is invertible. If
\(F<\sigma<\tfrac32\tilde c\), then \(\Phi_\sigma\) decays at both
ends and \(\Phi_\sigma^{-1}\) grows at both ends. Thus the kernel is
one-dimensional and the adjoint kernel is zero, so the operator is
surjective. Finally, \(\sigma\mapsto A_\sigma\) is norm-continuous from
\(D_b\) to \(L^2\). To see the right-inverse uniformity explicitly, fix
\(\sigma_0\in I\) and one right inverse \(Q_{\sigma_0}\). For nearby
\(\sigma\), \[
A_\sigma Q_{\sigma_0}
=I-(\sigma-\sigma_0)Q_{\sigma_0},
\] whose right-hand side is inverted by a Neumann series. A finite cover
of \(I\) gives (B.12); the same argument with inverses gives (B.11) on
\(J\). \(\square\)

We now apply Lemma B.1. For \(\sigma\in J\), let
\(G_\sigma=A_\sigma^{-1}\). Then \[
A_\sigma+R_\tau=(I+R_\tau G_\sigma)A_\sigma,\qquad
\sup_{\sigma\in J}\|R_\tau G_\sigma\|_{2\to2}
\le\frac{C_J}{|\tau|}. \tag{B.16}
\] Hence the perturbed operator is invertible once \(|\tau|\) exceeds a
single threshold. For \(\sigma\in I\), \[
(A_\sigma+R_\tau)Q_\sigma=I+R_\tau Q_\sigma,\qquad
\sup_{\sigma\in I}\|R_\tau Q_\sigma\|_{2\to2}
\le\frac{C_I}{|\tau|}. \tag{B.17}
\] The right-hand side is then invertible, so \(A_\sigma+R_\tau\) is
surjective. The compact set
\(\{A_\sigma:\sigma\in I\}\subset\mathcal B(D_b,L^2)\) lies in the open
Fredholm set, so it has a uniform Fredholm-stability radius. Equation
(B.5) gives
\(\sup_{\sigma\in I}\|R_\tau\|_{D_b\to L^2}=O(T^{-1})\); choosing \(T\)
once therefore keeps the entire path
\(A_\sigma+sR_\tau\), \(0\le s\le1\), Fredholm for every
\(\sigma\in I\). Its index is one by Lemma B.2; surjectivity therefore
leaves a one-dimensional kernel. Since \(O<0<\tfrac32\tilde c\), every compact collar
interval used in Proposition 5.6 is covered by either (B.16) or (B.17),
with one uniform threshold \(T_\delta\).

\begin{center}\rule{0.5\linewidth}{0.5pt}\end{center}

\section{Appendix C. Huang--Qin--Wang--Wei profiles and the full fixed-profile
hypothesis registry}\label{appendix-c.-hqww-profiles-and-the-full-fixed-profile-registry}

This appendix proves Propositions 2.7--2.8. The argument applies
separately to every fixed point of the Huang--Qin--Wang--Wei map. It
never selects a branch in \(a\).

We first make the imported fixed-point interface explicit. Following
{[}5, Section 3.1{]}, put \[
\rho(x)=(1+|x|)^{-1/2},\qquad
\mathbb V=\{f\in C(\mathbb R):f(-x)=f(x),\quad
\|f\|_{L^\infty_\rho}:=\|\rho f\|_\infty<\infty\}. \tag{C.0a}
\] With \(\eta=(3\cdot2^{20}\sqrt2)^{-1}\), define \[
\begin{aligned}
\mathbb D=\{f\in\mathbb V:\;&f(0)=1,\quad
(1-x^2)_+\le f(x)\le1,\\
&f\text{ is nonincreasing on }[0,\infty),\quad
s\mapsto f(\sqrt s)\text{ is convex},\\
&f'_{-}(1/2)\le-\eta\}.
\end{aligned} \tag{C.0b}
\] For \(f\in\mathbb D\), set \[
\mathbf T(f)(x)=\frac1\pi\int_0^\infty f(y)
\left[\frac yx\log\left|\frac{x+y}{x-y}\right|-2\right]dy,
\qquad
c(f)=\frac2\pi\int_0^\infty\frac{1-f(y)}{y^2}\,dy. \tag{C.0c}
\] The continuous value at the origin is \(\mathbf T(f)(0)=0\). For
\(0<a\le1\), \[
\mathbf T_a(f)=
\left(1+\frac{2a\,\mathbf T(f)}{(1-a/3)c(f)}\right)_+, \tag{C.0d}
\] and \[
R_a(f)(x)=\mathbf T_a(f)(x)^{1/a}
\exp\left[
\frac{1-a}{a}\int_0^x
\frac{\mathbf T_a(f)(y)-1}{y\mathbf T_a(f)(y)}\,dy
\right]. \tag{C.0e}
\] If \(\mathbf T_a(f)\) has a first positive zero, the expression in
(C.0e) is taken on the interval before that zero and extended by zero
after it; the result is extended evenly to \(\mathbb R\). At the
removable endpoint \(a=0\), define the even extension by \[
R_0(f)(x)=
\exp\left[
\frac2{c(f)}
\left(\mathbf T(f)(x)+\int_0^x\frac{\mathbf T(f)(y)}y\,dy\right)
\right], \tag{C.0f}
\]

For a general \(f\in\mathbb D\), monotonicity gives
\(f_\infty=\lim_{x\to\infty}f(x)\), and {[}5, Section 3.1{]} defines \[
b(f)=\frac2\pi\int_0^\infty(f(y)-f_\infty)\,dy\in[0,\infty].
\tag{C.0g}
\] For a fixed point, {[}5, Lemmas 4.1--4.2{]} give \(f_\infty=0\) and
\(b(f)<\infty\), so the integrand may then be replaced by \(f(y)\).
{[}5, Proposition 3.1{]} gives the raw similarity constants \[
c_{l,0}=\frac{1-a/3}{2}c(f)-ab(f),\qquad
c_{\omega,0}=\frac{1-a/3}{2}c(f)-b(f). \tag{C.0h}
\] The scalar used in Section C.5 is \[
\begin{aligned}
Q(f)&=\frac1{\pi^2}\int_0^\infty\int_0^\infty f(x)f(y)
\left[
\left(\frac xy+\frac yx\right)
\log\left|\frac{x+y}{x-y}\right|-2
\right]dx\,dy
=-\frac2\pi\int_0^\infty\mathbf T(f)'(x)\,x f(x)\,dx,\\
\mu(f)&=\frac{2Q(f)}{b(f)^2}.
\end{aligned} \tag{C.0i}
\] The quantities \(b,Q,\mu\) in (C.0g)--(C.0i) are only invoked for
fixed points, where the cited results and Section C.5 justify their
finiteness.

For clarity, we record exactly which facts are imported from {[}5{]}.

\begin{enumerate}
\def\labelenumi{\arabic{enumi}.}
\tightlist
\item
  {[}5, Proposition 3.1 and equations (3.8)--(3.9){]} supply the raw
  profile and similarity constants, while {[}5, Theorem 3.11{]} supplies
  fixed-point existence at each parameter.
\item
  {[}5, Corollary 4.6 and Theorem 4.10{]} give focusing in the stated
  range; {[}5, Lemma 4.9{]} gives the tail input used in (C.18); and
  {[}5, Theorem 4.12{]}
  supplies the regularity transfer.
\item
  {[}5, Theorem 4.3{]} supplies the formula for \(\mu\) and the
  parameter-uniform bound \(\mu\le\overline\mu<1\) used in Section C.5.
\item
  {[}5, Lemma 3.10{]} gives closedness and compactness of the
  fixed-point class in its parameter-independent weighted topology.
  After equation (4.9), the proof of {[}5, Lemma 4.7{]} notes the
  corresponding varying-parameter convergence by modifying the proof of
  {[}5, Theorem 3.7{]}. Lemma C.3 below supplies the details, especially
  at \(a=0\), and deduces sequential closedness of the graph
  \(\{(a,f):R_a(f)=f\}\).
\item
  The uniqueness input from {[}5, Section 5.3{]} is precisely the
  limiting \(a=0\) fixed point used in Section C.5.
\end{enumerate}

\subsection{C.1 Raw profile and
normalization}\label{c.1-raw-profile-and-normalization}

Let \(f\in\mathbb D\) be a fixed point of the map in (C.0d)--(C.0f). In
the gauge of {[}5{]}, write \[
\omega_0(x)=-xf(x),
\qquad
q_0=H\omega_0,
\qquad
U_0'=q_0,
\qquad
U_0(0)=0.
\tag{C.1}
\] The raw profile and similarity constants are produced in {[}5,
Proposition 3.1 and equations (3.8)--(3.9){]}. The integral
representation and finiteness of the scalar \(b(f)\), together with
\(c_{\omega,0}<0\), follow from {[}5, Lemmas 4.1--4.2{]}. For \[
0<a<\underline a=\frac{400}{848-9\pi^2},
\] {[}5, Corollary 4.6 and Theorem 4.10{]} put every such fixed point in
the focusing case, with \(c_{l,0}>0\).

Set \[
\alpha=-\frac1{c_{\omega,0}}>0,
\qquad
\beta=\frac4\alpha>0,
\tag{C.2}
\] and use the scaling {[}5, equation (2.3){]}: \[
\Omega(y)=\alpha\omega_0(\beta y),
\qquad
q(y)=\alpha q_0(\beta y),
\qquad
U(y)=\frac{\alpha}{\beta}U_0(\beta y),
\qquad
c_l=\alpha c_{l,0}.
\tag{C.3}
\] Because \(f(0)=1\), one has \(\omega_0'(0)=-1\). Thus (C.2)--(C.3)
give \[
c_\omega=-1,
\qquad
\Omega'(0)=-4.
\tag{C.4}
\] Equation (4.2) of {[}5{]} gives
\(c_{\omega,0}/c_{l,0}=1-r_{\mathrm{HQWW}}\), where
\(r_{\mathrm{HQWW}}=r_a(f)\) is the decay exponent used there.
Consequently \[
s_{\mathrm{prof}}:=r_{\mathrm{HQWW}}-1=\frac1{c_l}>1.
\tag{C.5}
\] The focusing, tail, and regularity results {[}5, Corollaries 4.5--4.6
and Theorems 4.10 and 4.12{]} give, after the positive scaling, \[
0<c_l<1,
\qquad
|\Omega(x)|\le C(1+x)^{-s_{\mathrm{prof}}},
\qquad
\Omega(x)\sim-C_\Omega x^{-s_{\mathrm{prof}}},
\tag{C.6}
\] and \[
f(x)\sim C_f x^{-1-s_{\mathrm{prof}}},
\qquad
\Omega'\in H^p(\mathbb R)
\quad\hbox{for every }p\ge0.
\tag{C.7}
\] The exponent in the asymptotic formula for \(f\) is stated in the raw
gauge and is unchanged by the fixed positive rescaling.

\subsection{C.2 Weighted derivative
transfer}\label{c.2-weighted-derivative-transfer}

We use a short tail interpolation argument to obtain the three weighted
derivative integrals required by \((K4^+)\).

\noindent\textbf{Lemma C.1.} Suppose \(u\in C^\infty([1,\infty))\) satisfies \[
|u(x)|\le C(1+x)^{-\alpha},
\qquad \alpha>1,
\] and \(u^{(m)}\in L^2(1,\infty)\) for every integer \(m\ge1\). Then,
for every integer \(k\ge1\), \[
\int_1^\infty x|u^{(k)}(x)|^2\,dx<\infty.
\tag{C.8}
\]

\emph{Proof.} Put \(R_n=2^n\), \(I_n=[R_n,2R_n]\), and
\(I_n^*=[R_n/2,4R_n]\). Choose \[
m>\frac{k(2\alpha-1)}{2\alpha-2},
\qquad
\theta=\frac{k}{m}.
\] A rescaled interior interpolation inequality gives \[
\|u^{(k)}\|_{L^2(I_n)}
\le
C\|u\|_{L^2(I_n^*)}^{1-\theta}
 \|u^{(m)}\|_{L^2(I_n^*)}^{\theta}
+C R_n^{-k}\|u\|_{L^2(I_n^*)},
\tag{C.9}
\] with \(C\) independent of \(n\). Since \[
\|u\|_{L^2(I_n^*)}^2\le C R_n^{1-2\alpha},
\] writing \(E_n=\|u^{(m)}\|_{L^2(I_n^*)}^2\) yields \[
\int_{I_n}x|u^{(k)}|^2\,dx
\le
C R_n^\gamma E_n^\theta
+C R_n^{2-2k-2\alpha},
\qquad
\gamma=2-2\alpha+(2\alpha-1)\frac{k}{m}<0.
\tag{C.10}
\] The enlarged annuli have finite overlap, so \(\sum_nE_n<\infty\).
Discrete Hölder gives \[
\sum_nR_n^\gamma E_n^\theta
\le
\left(\sum_nR_n^{\gamma/(1-\theta)}\right)^{1-\theta}
\left(\sum_nE_n\right)^\theta<\infty.
\] The second term in (C.10) is summable as well, and smoothness handles
the remaining compact interval. \(\square\)

Apply Lemma C.1 to \(u=\Omega\) using (C.6)--(C.7). Then \[
J_k^2=\int_0^\infty x|\Omega^{(k)}(x)|^2\,dx<\infty,
\qquad k=1,2,3.
\tag{C.11}
\]

\subsection{C.3 Velocity quotient, integrability, and far-field
regimes}\label{c.3-velocity-quotient-integrability-and-far-field-regimes}

Return first to the raw gauge. By {[}5, equation (3.2){]}, set \[
G_0(x):=\mathbf T(f)(x)+b(f)=\frac{U_0(x)}x.
\tag{C.12}
\] Equivalently, the positive-kernel representation is \[
G_0(x)=\frac1\pi\int_0^\infty
f(y)\frac yx\log\left|\frac{x+y}{x-y}\right|\,dy,
\tag{C.13}
\] so \(G_0\ge0\). Lemma 3.4 of {[}5{]} gives \(G_0'\le0\) and states
that \(z\mapsto G_0(\sqrt z)\) is convex. Since \[
\frac{d^2}{dz^2}G_0(\sqrt z)
=\frac1{4x}\left(\frac{G_0'(x)}x\right)',
\] we have \[
G_0'\le0,
\qquad
\left(\frac{G_0'}x\right)'\ge0.
\tag{C.14}
\]

Choose \(1<\delta<\min\{s_{\mathrm{prof}},2\}\). The asymptotic in (C.7) gives
\(\sup_x x^{1+\delta}f(x)<\infty\), and {[}5, Lemma 4.9{]} yields \[
G_0(x)=O(x^{-\delta}).
\tag{C.15}
\] Thus \(G_0\in L^1(0,\infty)\) and \(xG_0(x)\to0\) at both endpoints.
Since \(q_0=U_0'=(xG_0)'\) and \(G_0\) is decreasing, \[
\begin{aligned}
\|q_0\|_1
&\le\int_0^\infty G_0(x)\,dx
 +\int_0^\infty x(-G_0'(x))\,dx\\
&=2\int_0^\infty G_0(x)\,dx<\infty.
\end{aligned}
\tag{C.16}
\]

Put \(h=-G_0'\ge0\). Equation (C.14) says that \(h(x)/x\) is
nonincreasing. Hence \[
G_0(x/2)-G_0(x)
=\int_{x/2}^x h(t)\,dt
\ge\frac{h(x)}x\int_{x/2}^x t\,dt
=\frac38xh(x),
\] and therefore \[
x|G_0'(x)|\le\frac83G_0(x/2),
\qquad
|q_0(x)|\le G_0(x)+\frac83G_0(x/2).
\tag{C.17}
\]

The tail inputs {[}5, equation (4.2), Lemma 4.9, and Theorems 4.10 and
4.12{]} now give three regimes. If \(1<s_{\mathrm{prof}}<2\), the
endpoint estimate of {[}5, Lemma 4.9{]} with
\(\delta=s_{\mathrm{prof}}\), together with (C.17), gives \[
G_0(x)+|q_0(x)|=O(x^{-s_{\mathrm{prof}}}).
\tag{C.18}
\] This endpoint statement is distinct from the strict interior choice
in (C.15). If \(s_{\mathrm{prof}}>2\), then
\(\int_0^\infty y^2f(y)\,dy<\infty\), and the
second-moment conclusion of Lemma 4.9 gives \[
G_0(x)=O(x^{-2}),
\qquad
G_0'(x)=O(x^{-3}),
\qquad
q_0(x)=O(x^{-2}).
\tag{C.19}
\] In the resonant case \(s_{\mathrm{prof}}=2\), for every fixed
\(1<\widehat\alpha<2\), \[
G_0(x)+|q_0(x)|=O(x^{-\widehat\alpha}).
\tag{C.20}
\] No endpoint \(O(x^{-2})\) estimate is asserted or used in this case.

Define \(\alpha_{\rm reg}=s_{\mathrm{prof}}\) for
\(1<s_{\mathrm{prof}}<2\), set
\(\alpha_{\rm reg}=2\) for \(s_{\mathrm{prof}}>2\), and choose any
\(\alpha_{\rm reg}=\widehat\alpha\in(1,2)\) when
\(s_{\mathrm{prof}}=2\).
Under (C.3), \[
G(y):=\frac{U(y)}y=\alpha G_0(\beta y).
\tag{C.21}
\] Positivity, monotonicity, integrability, and the preceding rates are
preserved. Since \(\alpha_{\rm reg}\le s_{\mathrm{prof}}\), the
\(\Omega\)-bound below is also the weakened consequence of (C.6). Thus
\[
|\Omega(y)|\le C(1+y)^{-\alpha_{\rm reg}},
\qquad
|q(y)|\le C(1+y)^{-\min(\alpha_{\rm reg},2)},
\tag{C.22}
\] and (C.15), with \(\delta>1\), gives \[
U(y)=yG(y)\longrightarrow0.
\tag{C.23}
\]

\subsection{C.4 Transfer of the complete
registry}\label{c.4-transfer-of-the-complete-registry}

The normalized profile satisfies \[
\Omega+(c_ly+aU)\Omega'=q\Omega,
\qquad
q=H\Omega.
\tag{C.24}
\] The source results in C.1 give smoothness, oddness, focusing,
\(0<c_l<1<2\), and \(\Omega\in L^2\). The normalization is (C.4).

By (C.14) and (C.21), the normalized quotient \(G=U/y\) is nonnegative
and decreasing. Smoothness at zero gives \(G(0)=U'(0)=q(0)\), so \[
c_l\le\frac{b(y)}y=c_l+aG(y)
\le c_l+aq(0)=:\tilde c.
\tag{C.25}
\] This proves the focusing transport bounds. Smoothness and oddness
give finite \(M_1,M_2\) on a sufficiently small interval and \[
\Omega(y)-\Omega'(0)y=O(y^3),
\] so \(M_3<\infty\). Comparing the order-\(y\) terms in (C.24) gives
the exact origin identity \[
q(0)-\bigl(c_l+aq(0)\bigr)=1.
\tag{C.26}
\] Moreover, (C.7), boundedness of the Hilbert transform on Sobolev
spaces, and Sobolev embedding give \(q''\in L^\infty\), hence
\(M_4<\infty\). Thus every item of \((H_{\mathrm{prof}})\) holds.

Equations (C.6)--(C.7) imply \(\Omega\in H^m(\mathbb R)\) for every
integer \(m\ge0\). Boundedness of the Hilbert transform and Sobolev
embedding therefore give \[
Q_0,Q_1,W_0,W_1,W_2<\infty.
\] Equation (C.11) gives \(J_1,J_2,J_3<\infty\), and (C.16), transferred
by (C.3), gives \(I_q<\infty\). Hence \((K4^+)\) holds.

The pointwise bounds in \((D_\infty)\) are exactly (C.22), and (C.23)
gives \(U_\infty=0\). Consequently \[
b(y)=c_ly+aU(y)=c_ly+O(1).
\tag{C.27}
\] The weighted tails in \(J_1,J_2,J_3\) tend to zero by absolute
continuity. For \(v=\Omega'\) and \(v=\Omega''\), \[
\sup_{y\ge R}|v(y)|^2
\le
2\|v\|_{L^2(R,\infty)}
 \|v'\|_{L^2(R,\infty)}.
\tag{C.28}
\] Together with (C.6), (C.11), Sobolev regularity, and (C.22), this
proves \[
\varepsilon_\Omega(R)\to0,
\qquad
\varepsilon_q(R)\to0,
\qquad
\varepsilon_V(R)\to0.
\tag{C.29}
\] Every entry of \((H_{\mathrm{prof}})+(K4^+)+(D_\infty)\) has now been
verified for the normalized profile. Since the fixed point was
arbitrary, Proposition 2.7 holds for the entire Huang--Qin--Wang--Wei fixed-point set.

\subsection{\texorpdfstring{C.5 Uniform small-\(a\) separation over the
fixed-point
set}{C.5 Uniform small-a separation over the fixed-point set}}\label{c.5-uniform-small-a-separation-over-the-fixed-point-set}

For a fixed point \(f\in\mathbb D\), {[}5, Theorem 4.3{]} uses \[
b(f)=\frac2\pi\int_0^\infty f(x)\,dx,
\qquad
\mu(f)=\frac{2Q(f)}{b(f)^2},
\tag{C.30}
\] where \[
Q(f)=\frac1{\pi^2}\int_0^\infty\int_0^\infty
f(x)f(y)
\left[
\left(\frac{x}{y}+\frac{y}{x}\right)
\log\left|\frac{x+y}{x-y}\right|-2
\right]dx\,dy.
\tag{C.31}
\] The formula for the normalized similarity exponent in {[}5, Theorem
4.3 and Corollary 4.5{]} is \[
c_l=\frac{1-a(2-\mu)}{1-a\mu}.
\tag{C.32}
\] Combining (C.32) with the origin identity \[
\tilde c=\frac{c_l+a}{1-a}
\] gives the exact scalar relation \[
F-O=\frac{a(2\mu-1)}{1-a\mu}.
\tag{C.33}
\] For \(0<a<1\), the bounds in {[}5, Theorem 4.3{]} give \(0<\mu<1\).
Hence \(F<O\) is equivalent to \(\mu<1/2\).

At \(a=0\), {[}5, Section 5.3{]} proves uniqueness only for the limiting
fixed-point problem: \[
f_0(x)=\frac1{1+x^2}.
\tag{C.34}
\] Here \(b(f_0)=1\), \(\omega_0(x)=-x/(1+x^2)\), and
\(U_0(x)=\arctan x\). The formula used in the proof of {[}5, Theorem
4.3{]} gives \[
Q(f_0)=\frac2\pi\int_0^\infty
\left(\frac{U_0(x)}x\right)'\omega_0(x)\,dx.
\] With \(x=\tan t\), \[
Q(f_0)
=\frac2\pi\int_0^{\pi/2}
\left(t\cot t-\cos^2t\right)dt
=\log2-\frac12,
\] and therefore \[
\mu(f_0)=2\log2-1<\frac12.
\tag{C.35}
\]

We next record the continuity input for \(Q\).

\noindent\textbf{Lemma C.2.} Let \(f_n,f\ge0\) be nonincreasing on
\((0,\infty)\), suppose \(f_n\to f\) in \(L^1(0,\infty)\), and assume
their \(L^1\) norms are uniformly bounded. Then \(Q(f_n)\to Q(f)\).

\emph{Proof.} Set \[
g_n(s)=e^sf_n(e^s),
\qquad
g(s)=e^sf(e^s).
\] Then \[
\|g_n-g\|_{L^1(\mathbb R)}
=\|f_n-f\|_{L^1(0,\infty)}.
\tag{C.36}
\] Monotonicity gives \[
xf_n(x)\le\int_0^x f_n(y)\,dy\le\|f_n\|_1,
\] so \(g_n\) is uniformly bounded in \(L^\infty\); the same holds for
\(g\). The bracket in (C.31) is nonnegative. After the substitutions
\(x=e^s\), \(y=e^t\), Tonelli's theorem rewrites (C.31) as \[
Q(f_n)=\int_{\mathbb R}\int_{\mathbb R}
g_n(s)g_n(t)k(s-t)\,ds\,dt,
\tag{C.37}
\] where \[
k(v)=\frac1{\pi^2}
\left[
(e^v+e^{-v})\log\left|\coth\frac v2\right|-2
\right].
\tag{C.38}
\] This even kernel satisfies \[
k(v)=\frac2{\pi^2}\log\frac2{|v|}+O(1)
\quad(v\to0),
\] and \[
k(v)=\frac8{3\pi^2}e^{-2|v|}
+O(e^{-4|v|})
\quad(|v|\to\infty).
\] Thus \(k\in L^1(\mathbb R)\). If \[
B(g,h)=\int_{\mathbb R}g(s)(k*h)(s)\,ds,
\] then \(B\) is symmetric and Young's inequality gives \[
|B(g,h)|\le\|k\|_1\|g\|_1\|h\|_\infty.
\] Using symmetry in one of the two difference terms, \[
|Q(f_n)-Q(f)|
\le
\|k\|_1\|g_n-g\|_1
\bigl(\|g_n\|_\infty+\|g\|_\infty\bigr)
\longrightarrow0.
\tag{C.39}
\] \(\square\)

\noindent\textbf{Lemma C.3 (sequential closedness of the Huang--Qin--Wang--Wei fixed-point
graph).} Let \(a_n,a_*\in[0,1]\) with \(a_n\to a_*\), and let
\(f_n\in\mathbb D\) satisfy \[
R_{a_n}(f_n)=f_n.
\] Then, after passing to a subsequence, there exists
\(f_*\in\mathbb D\) such that \[
f_n\longrightarrow f_*
\quad\hbox{in }L^\infty_\rho,
\qquad R_{a_*}(f_*)=f_*, \tag{C.40}
\] where \(L^\infty_\rho\) is the weighted topology used to define the
compact class \(\mathbb D\) in {[}5{]}. If \(a_*=0\), then \[
f_*(x)=f_0(x):=\frac1{1+x^2}.
\]

\emph{Proof.} By {[}5, Lemma 3.10{]}, \(\mathbb D\) is closed and
compact in \(L^\infty_\rho\), so a subsequence converges to some
\(f_*\in\mathbb D\). Reference {[}5, Lemma 3.3{]} and the kernel
estimates in the proof of {[}5, Theorem 3.7{]} give \[
c(f_n)\longrightarrow c(f_*),\qquad
\mathbf T(f_n)\longrightarrow\mathbf T(f_*)
\quad\hbox{locally uniformly}. \tag{C.40a}
\]

First suppose \(a_*>0\). Equations (C.0c)--(C.0d), continuity of the
positive-part map, and (C.40a) imply \[
\mathbf T_{a_n}(f_n)\longrightarrow\mathbf T_{a_*}(f_*)
\quad\hbox{locally uniformly}. \tag{C.40b}
\] On a compact subinterval where the limiting function is bounded away
from zero, formula (C.0e) gives uniform convergence of the corresponding
\(R_{a_n}(f_n)\). If \(\mathbf T_{a_*}(f_*)\) has a first zero \(L\),
apply this argument on \([0,L-\delta]\). The remaining collar is
controlled by \[
0\le R_a(f)\le\mathbf T_a(f)\le1,\qquad 0<a\le1,
\] which follows from the proof of {[}5, Lemma 3.6{]}. Letting
\(\delta\downarrow0\) proves local uniform convergence across the zero.
All constants are uniform because \(a_n\) remains in a compact
subinterval of \((0,1]\).

It remains to prove the assertion at \(a_*=0\). Put \[
A_{a,f}(x)=\frac{2\mathbf T(f)(x)}{(1-a/3)c(f)}.
\] Then \(\mathbf T_a(f)=(1+aA_{a,f})_+\), and (C.40a) gives \[
A_{a_n,f_n}\longrightarrow
A_{0,f_*}=\frac{2\mathbf T(f_*)}{c(f_*)}
\quad\hbox{locally uniformly}. \tag{C.40c}
\] The convexity calculation in {[}5, Lemma 3.4{]}, together with \[
\lim_{x\to0}\frac{\mathbf T(f)'(x)}{2x}=-\frac{c(f)}3,
\] gives the uniform origin estimate \[
|\mathbf T(f)(x)|\le Cx^2,\qquad
f\in\mathbb D,\quad 0\le x\le1. \tag{C.40d}
\] Uniformity follows from the bounds for \(c(f)\) in
{[}5, Lemma 3.3{]} and the parameter-independent convexity constraints
defining \(\mathbb D\). Thus \(A_{a_n,f_n}(y)/y\) is dominated near zero
by \(Cy\).

On each fixed compact interval the positive part is inactive for all
large \(n\). For indices with \(a_n>0\), formula (C.0e) becomes \[
\begin{aligned}
\log R_{a_n}(f_n)(x)
&=\frac{\log(1+a_nA_{a_n,f_n}(x))}{a_n}\\
&\quad +(1-a_n)\int_0^x
\frac{A_{a_n,f_n}(y)}
{y(1+a_nA_{a_n,f_n}(y))}\,dy .
\end{aligned} \tag{C.40e}
\] Using (C.40c)--(C.40d) and \(a^{-1}\log(1+az)\to z\) uniformly for
bounded \(z\), dominated convergence yields \[
\begin{aligned}
\log R_{a_n}(f_n)(x)
&\longrightarrow
\frac2{c(f_*)}\left[
\mathbf T(f_*)(x)+\int_0^x\frac{\mathbf T(f_*)(y)}y\,dy
\right]\\
&=\log R_0(f_*)(x)
\end{aligned} \tag{C.40f}
\] locally uniformly. The same conclusion follows directly from (C.0f)
along any indices with \(a_n=0\).

Finally, {[}5, Lemma 3.6{]} gives \(R_a(f)\in\mathbb D\), so
\(0\le R_a(f)\le1\). Hence \[
\sup_{|x|\ge X}\rho(x)
|R_{a_n}(f_n)(x)-R_{a_*}(f_*)(x)|
\le2\rho(X), \tag{C.40g}
\] uniformly in \(n\). Local uniform convergence followed by
\(X\to\infty\) proves convergence in \(L^\infty_\rho\). Passing to the
limit in \(f_n=R_{a_n}(f_n)\) gives \(f_*=R_{a_*}(f_*)\). If \(a_*=0\),
the uniqueness theorem of {[}5, Section 5.3{]} identifies
\(f_*(x)=(1+x^2)^{-1}\). This is the varying-parameter convergence
invoked after equation (4.9) in the proof of {[}5, Lemma 4.7{]}, now
with the endpoint argument written out. \(\square\)

We now prove the all-fixed-point statement by contradiction. If it
failed, there would be \(a_n\downarrow0\) and fixed points
\(f_n\in\mathbb D\) of \(R_{a_n}\) such that \[
\mu(f_n)\ge\frac12.
\tag{C.41}
\] By Lemma C.3, after passing to a subsequence, \[
f_n\longrightarrow f_0=\frac1{1+x^2}
\quad\hbox{in }L^\infty_\rho. \tag{C.42}
\] This convergence supplies the parameter-closure step used below.

The parameter bounds {[}5, Lemma 3.3 and Theorem 4.3{]} give \[
c(f_n)\longrightarrow c(f_0)=1,
\qquad
\frac{b(f_n)}{c(f_n)}
=\frac{1-a_n/3}{1+a_n\mu(f_n)}.
\tag{C.43}
\] The same theorem supplies a uniform bound
\(0\le\mu(f_n)\le\overline\mu<1\), so (C.43) yields \[
b(f_n)\longrightarrow1=b(f_0).
\tag{C.44}
\] The compactness convergence is pointwise, and
\(f_n,f_0\in\mathbb D\) are nonnegative and nonincreasing. Because
\(b(f)\) in (C.30) is exactly a constant multiple of the mass, (C.44)
gives convergence, and in particular uniform boundedness, of their
\(L^1\) norms. Scheffé's lemma therefore gives the required strong
convergence \[
f_n\longrightarrow f_0
\quad\hbox{in }L^1(0,\infty).
\tag{C.45}
\] All hypotheses of Lemma C.2 are now explicit, so it implies
\(Q(f_n)\to Q(f_0)\), and hence \[
\mu(f_n)=\frac{2Q(f_n)}{b(f_n)^2}
\longrightarrow2\log2-1<\frac12,
\] contradicting (C.41).

The arbitrary-sequence argument in fact proves \[
\lim_{a\downarrow0}
\sup_{\substack{f\in\mathbb D\\R_a(f)=f}}
\left|\mu(f)-(2\log2-1)\right|=0.
\tag{C.46}
\] Choose \(a_{\mathrm{sep}}\in(0,1/2]\) so that \(\mu(f)<1/2\) for
every fixed point whenever \(0<a<a_{\mathrm{sep}}\). Such a choice is
possible by (C.46), since \(1/2<\underline a\). Equations (C.33) and
(C.46) prove Proposition 2.8.

The upper cap is forced by a direct witness at \(a=1/2\). Consider \[
f_{1/2}(x)=\frac4{(2+x^2)^2}. \tag{C.47}
\] This function lies in \(\mathbb D\): it is even, \(f_{1/2}(0)=1\),
and \[
f_{1/2}'(x)=-\frac{16x}{(2+x^2)^3}\le0,\qquad
\frac{d^2}{ds^2}\frac4{(2+s)^2}=\frac{24}{(2+s)^4}>0.
\] For \(0\le s\le1\), \[
4-(1-s)(2+s)^2=s^2(3+s)\ge0,
\] so \((1-x^2)_+\le f_{1/2}\le1\). Finally, \[
f_{1/2,-}'(1/2)=-\frac{512}{729}<-\eta.
\]

The fixed-point identity can also be checked directly. Evaluation of
(C.0c) gives \[
\mathbf T(f_{1/2})(x)
=-\frac{x^2}{\sqrt2(2+x^2)},\qquad
c(f_{1/2})=\frac3{2\sqrt2}. \tag{C.48}
\] Thus \[
\mathbf T_{1/2}(f_{1/2})(x)
=1+\frac{\mathbf T(f_{1/2})(x)}{(5/6)c(f_{1/2})}
=\frac{10+x^2}{5(2+x^2)}=:g(x)>0. \tag{C.49}
\] Moreover, \[
\int_0^x\frac{g(y)-1}{yg(y)}\,dy
=-4\int_0^x\frac{y}{y^2+10}\,dy
=-2\log\frac{x^2+10}{10}. \tag{C.50}
\] Since \((1-a)/a=1\) at \(a=1/2\), equations (C.0e) and (C.50) give \[
\begin{aligned}
R_{1/2}(f_{1/2})(x)
&=g(x)^2\exp\left(\int_0^x\frac{g(y)-1}{yg(y)}\,dy\right)\\
&=\frac{(10+x^2)^2}{25(2+x^2)^2}
\frac{100}{(10+x^2)^2}
=f_{1/2}(x).
\end{aligned} \tag{C.51}
\] This independently closes the membership and fixed-point statement
also recorded in {[}5, Section 5.2{]}.

The relevant scalars follow without using (C.32): \[
b(f_{1/2})
=\frac2\pi\int_0^\infty\frac4{(2+x^2)^2}\,dx
=\frac1{\sqrt2}, \tag{C.52}
\] and, by (C.0i), \[
\begin{aligned}
Q(f_{1/2})
&=-\frac2\pi\int_0^\infty
\mathbf T(f_{1/2})'(x)\,x f_{1/2}(x)\,dx\\
&=\frac{16\sqrt2}{\pi}
\int_0^\infty\frac{x^2}{(2+x^2)^4}\,dx
=\frac18.
\end{aligned} \tag{C.53}
\] Therefore \[
\mu(f_{1/2})
=\frac{2Q(f_{1/2})}{b(f_{1/2})^2}
=\frac12. \tag{C.54}
\] Equation (C.0h) gives the raw constants \(c_{l,0}=\sqrt2/16\) and
\(c_{\omega,0}=-3\sqrt2/16\). After the normalization \(c_\omega=-1\)
used in the main text, \(c_l=1/3\), and direct scaling gives
\[
\Omega(y)=-\frac{1024y}{(16+9y^2)^2},\qquad
q(y)=\frac{128}{3}\frac{16-9y^2}{(16+9y^2)^2},\qquad
U(y)=\frac{128y}{3(16+9y^2)}.
\]
These rational formulas directly verify all four requirements of
\(\mathrm{Adm}(1/2)\) in {[}6, Definition 4.2{]}, including the
pointwise derivative rates in (H3). For the transport weight
\(b(y)=c_ly+aU(y)\), they give
\[
\frac{b(y)}y=\frac{80+9y^2}{3(16+9y^2)}\in\left[\frac13,\frac53\right],
\]
and the remaining derivative and endpoint bounds in (H4) follow by
direct differentiation. Thus the two packages have nonempty
intersection. However, (C.33), together with \(c_l=1/3\), gives
\(F=O=-5/6\) for this profile. No
all-fixed-point separation interval of the form \(0<a<A\) can therefore
have \(A>1/2\). This does not identify an optimal positive threshold or
assert separation for every \(0<a<1/2\).

Finally, {[}5, Theorem 3.11{]} supplies at least one fixed point at
every parameter in the stated range. The conclusion concerns the entire
fixed-point set and uses uniqueness only for the limiting \(a=0\)
problem.

\section*{Acknowledgements}\label{acknowledgements}

The author used OpenAI Codex (GPT-5.6 Sol), Anthropic Claude (Opus 5),
and Moonshot AI Kimi K3 for
drafting and revising mathematical exposition, literature
triage, symbolic cross-checks, and internal critical review. The author
critically revised and checked all AI-assisted material and takes full
responsibility for the manuscript.

\section*{Funding}\label{funding}

No external funding was received specifically for this work.

\section*{Competing interests}\label{competing-interests}

The author declares no competing interests.

\section*{Data availability}\label{data-availability}

No new data were created or analysed in this study.

\section*{References}\label{references}

\begingroup
\newcommand{\refitem}{\par\hangindent=2em\hangafter=1\noindent}

\refitem {[}1{]} P. Constantin, P. D. Lax, A. Majda. \emph{A simple
one-dimensional model for the three-dimensional vorticity equation.}
Comm. Pure Appl. Math. 38 (1985) 715--724.

\refitem {[}2{]} S. De Gregorio. \emph{On a one-dimensional model for the
three-dimensional vorticity equation.} J. Stat. Phys. 59 (1990)
1251--1263.

\refitem {[}3{]} H. Okamoto, T. Sakajo, M. Wunsch. \emph{On a generalization of
the Constantin-Lax-Majda equation.} Nonlinearity 21 (2008) 2447--2461.

\refitem {[}4{]} P. M. Lushnikov, D. A. Silantyev, M. Siegel. \emph{Collapse
versus blow up and global existence in the generalized
Constantin-Lax-Majda equation.} J. Nonlinear Sci. 31 (2021), art. 82.
doi:10.1007/s00332-021-09737-x; arXiv:2010.01201.

\refitem {[}5{]} D. Huang, X. Qin, X. Wang, D. Wei. \emph{Self-similar
finite-time blowups with smooth profiles of the generalized
Constantin--Lax--Majda model.} Arch. Ration. Mech. Anal. 248 (2024),
art. 22. doi:10.1007/s00205-024-01971-3; arXiv:2305.05895v5.

\refitem {[}6{]} J. Xu. \emph{The spectral picture of self-similar collapse in
the Constantin-Lax-Majda equation.} arXiv:2607.19762 (2026).

\refitem {[}7{]} V. S. Rabinovich, S. Roch, B. Silbermann. \emph{Limit Operators
and Their Applications in Operator Theory.} Operator Theory: Advances
and Applications 150, Birkhäuser, 2004.

\refitem {[}8{]} R. B. Lockhart, R. C. McOwen. \emph{Elliptic differential
operators on noncompact manifolds.} Ann. Scuola Norm. Sup. Pisa Cl. Sci.
(4) 12 (1985) 409--447.

\refitem {[}9{]} R. B. Melrose. \emph{The Atiyah-Patodi-Singer Index Theorem.}
Research Notes in Mathematics 4, A K Peters, 1993.

\refitem {[}10{]} M. Lesch. \emph{Operators of Fuchs Type, Conical Singularities,
and Asymptotic Methods.} Teubner-Texte zur Mathematik 136, Teubner,
1997.

\refitem {[}11{]} V. A. Kondrat'ev. \emph{Boundary value problems for elliptic
equations in domains with conical or angular points.} Trans. Moscow
Math. Soc. 16 (1967) 227--313.

\refitem {[}12{]} R. Duduchava. \emph{Integral Equations with Fixed
Singularities.} Teubner-Texte zur Mathematik 24, BSB B. G. Teubner,
Leipzig, 1979.

\refitem {[}13{]} D. E. Edmunds, W. D. Evans. \emph{Spectral Theory and
Differential Operators.} Oxford Mathematical Monographs, Clarendon
Press, Oxford, 1987.

\refitem {[}14{]} T. M. Elgindi, T.-E. Ghoul, N. Masmoudi. \emph{Stable
self-similar blow-up for a family of nonlocal transport equations.}
Anal. PDE 14 (2021), 891--908. doi:10.2140/apde.2021.14.891.

\refitem {[}15{]} D. Huang, X. Qin, X. Wang. \emph{Multiscale self-similar
finite-time blowups of the Constantin--Lax--Majda model for the
three-dimensional Euler equations.} SIAM J. Math. Anal. 57 (2025),
4068--4096. doi:10.1137/24M168845X.

\refitem {[}16{]} D. Huang, J. Tong, X. Wang. \emph{Self-similar finite-time
blowups with singular profiles of the generalized Constantin--Lax--Majda
model: theoretical and numerical investigations.} arXiv:2603.25104
(2026).

\refitem {[}17{]} T. Kato. \emph{Perturbation Theory for Linear Operators.}
Reprint of the 1980 edition, Classics in Mathematics, Springer-Verlag,
Berlin, 1995.

\par
\endgroup

\end{document}